\documentclass[a4paper,11pt]{article}
\usepackage[utf8]{inputenc}
\usepackage{amsfonts}
\usepackage{amsmath}
\usepackage{amssymb}
\usepackage{amsthm}
\usepackage{color,a4wide}
\usepackage{algorithm,algpseudocode, algorithmicx}
\usepackage{graphicx}
\usepackage[caption=false]{subfig}
\usepackage{tikz}
\usepackage{mathtools}
\usetikzlibrary{arrows}
\tikzstyle{block}=[draw opacity=0.7,line width=1.4cm]
\usepackage{pgfplots}
\usepackage{hyperref}
\hypersetup{
 colorlinks=true,
 linkcolor=blue,
 filecolor=magenta,  
 urlcolor=cyan,
}

\usepackage{array,multirow}
\newcolumntype{B}[1]{>{\centering\arraybackslash}m{#1}}
\usepackage{colortbl}
\usepackage[sorting=none]{biblatex} %
\usepackage[T1]{fontenc}
\usepackage[normalem]{ulem}

\theoremstyle{plain}

\theoremstyle{plain}

\theoremstyle{plain}

\theoremstyle{definition}

\theoremstyle{plain}

\theoremstyle{remark}

\theoremstyle{remark}

\theoremstyle{definition}

\DeclareMathOperator*{\argmin}{\arg \min}
\newcommand{\R}{\mathbb{R}}
\newcommand\inv[1]{#1\raisebox{1.15ex}{$\scriptscriptstyle-\!1$}}

\title{Regularization of Inverse Problems: Deep Equilibrium Models versus Bilevel Learning}
\author{Danilo Riccio$^1$ \and Matthias J. Ehrhardt$^2$ \and Martin Benning$^3$}
\date{%
    $^1$Queen Mary University of London, United Kingdom\\%
    $^2$University of Bath, United Kingdom\\%
    $^3$University College London, United Kingdom%
}

\begin{document}

\maketitle

\begin{abstract}
 Variational regularization methods are commonly used to approximate solutions of inverse problems. In recent years, model-based variational regularization methods have often been replaced with data-driven ones such as the fields-of-expert model \cite{roth2009fields}. Training the parameters of such data-driven methods can be formulated as a bilevel optimization problem. In this paper, we compare the framework of bilevel learning for the training of data-driven variational regularization models with the novel framework of deep equilibrium models \cite{Bai2019} that has recently been introduced in the context of inverse problems \cite{Gilton2021}.
 We show that computing the lower-level optimization problem within the bilevel formulation with a fixed point iteration is a special case of the deep equilibrium framework. We compare both approaches computationally, with a variety of numerical examples for the inverse problems of denoising, inpainting and deconvolution.
\end{abstract}

\section{Introduction}
In inverse problems, a desired quantity can only be recovered indirectly from measurements through the inversion of a so-called forward operator. Inverse problems are found in many diverse research areas like medical and industrial imaging, signal and image processing, computer vision, or geophysics, with various applications utilizing inverse problems such as Computerized Tomography, Magnetic Resonance Imaging, or Transmission Electron Tomography, see e.g., \cite{natterer2001mathematical,scherzer2009variational} and references therein. 
For the majority of relevant inverse problems, an inverse of the forward operator is discontinuous and usually not unique. If we focus on inverse problems with linear forward operators in finite dimensions, an inverse (if it exists) is no longer discontinuous, but can still be very sensitive to small variations in the argument if the condition number of the linear forward operator (which is then a matrix) is large. 

Regularizations are parameterized operators that approximate inverses of forward operators in a continuous fashion \cite{tikhonov1987ill,engl1996regularization,benning2018modern}. Variational regularizations \cite{scherzer2009variational} are a special class of regularizations and usually assign the solution of a convex optimization problem to the operator input. Regularizations are equipped with hyper-parameters that are usually referred to as regularization parameters, and these regularization parameters have to be chosen dependent on the data error to guarantee that a regularization convergences to a (generalized) inverse when the error goes to zero.

A heuristic regularization parameter choice strategy that has been popularized in the past decade is the identification of optimal regularization parameters in a data-driven way with bilevel optimization \cite{DeLosReyes2013, Kunisch2013bilevel, Ochs2015, DelosReyes2017, Ehrhardt2021bileveloptimization, Crockett2021bilevel}. In a bilevel optimization problem, an upper-level objective that is constrained by a lower-level objective is optimized. In the context of variational regularization parameter estimation, one can for example, minimize an empirical risk between the regularization operator outputs and a set of desired outputs as the upper-level optimization problem, subject to the lower-level optimization problem that the regularization operator outputs have to solve the optimization problem associated with the variational regularization. 

In a development parallel to the use of bilevel optimization for the estimation of regularization parameters, more and more inverse problem solutions have been approximated with operators stemming from deep neural networks, with great empirical success \cite{Gregor2010, Jin2017, Zhu2018automap, Adler2018primaldual, Adler2017deepinversion, Kobler2017, Yang2018, Kobler2020tdv}. A recent development is the application of so-called deep equilibrium models \cite{Bai2019} in the context of inverse problems \cite{Gilton2021,heaton2021feasibility}. Deep equilibrium models are deep learning frameworks where a neural network is trained such that its fixed point approximates some ground truth solution.\\

None of the aforementioned papers have investigated the obvious link between deep equilibrium models and bilevel optimization. Our contributions are the framing of bilevel optimization problems with strongly convex lower-level problem as deep equilibrium models, and an empirical comparison of deep equilibrium models and bilevel optimization for three inverse problems.

The paper is organized as follows. In Section~\ref{sec:deep-equib-vs-bilevel}, we give a brief overview of the concepts of inverse problems, regularizations, bilevel optimization and deep equilibrium models. In Section~\ref{sec:arch_design_and_implementation}, we choose different bilevel and deep equilibrium models for comparison and explain their architectures and other design choices. In Section~\ref{sec:applications}, we present the three inverse problems denoising, inpainting and deblurring that we use for the numerical deep equilibrium vs bilevel optimization comparison in Section~\ref{sec:numerical-results}. We then wrap up with conclusions and give an outlook for open research questions in Section~\ref{sec:conclusions-and-outlook}.

\section{Deep equilibrium models and bilevel learning}\label{sec:deep-equib-vs-bilevel}
In this section, we briefly recall the fundamentals of inverse problems and variational regularization. Subsequently, we describe the concepts of bilevel optimization and deep equilibrium models in this context and highlight how they are linked mathematically. We conclude this section with a comment on why learning fixed points in a na\"{i}ve way does not work.

\subsection{Inverse problems}\label{sec:inverse-problems}

Many practical applications, such as Computed Tomography (CT) and Magnetic Resonance Imaging (MRI), require the reconstruction of quantities from indirect measurements. 
This task can be modeled as an inverse problem 
\begin{align}
 Ku=f,\label{eq:general_inverse_problem}
\end{align}
where the goal is to retrieve an unknown quantity $u$ given the operator $K$ and the \textit{ideal} measured data $f$, with $u\in\mathbb{R}^n, f\in\mathbb{R}^m, K\in\mathbb{R}^{m\times n}$. %
In this paper, we limit our analysis to linear and finite-dimensional operators $K$, so we can think of $K$ as a $m\times n$ matrix.

If the operator $K$ is invertible, the unique solution of \eqref{eq:general_inverse_problem} is $u=\inv{K}f$. 
In case the operator $K$ is not invertible, the Moore--Penrose inverse $K^+$ can be used to approximate the inverse instead. The Moore--Penrose inverse is the minimum norm solution of the normal equation $K^\top Ku = K^\top f$, which means that it minimizes the least squares function $J(u)=\frac{1}{2}\lVert Ku-f \rVert^2$, where $\lVert \cdot \rVert$ is the Euclidean norm, with minimal Euclidean norm $\| u \|$. By using the Moore--Penrose inverse to solve \eqref{eq:general_inverse_problem} for $u$, we solve
\begin{align}
 \bar u = K^{+} f.\label{eq:LSE_no_regularizer}
\end{align}
In most applications, this problem is usually ill-conditioned. In other words, as soon as we consider the presence of noise affecting the measurement, the worst-case error between retrieved data and desired data is strongly amplified. In this paper we will consider any measurement error to be additive, i.e., %
\begin{align}
 Ku+\delta=f^\delta,\label{eq:general_inverse_problem__with_additive_noise}
\end{align}
$\delta, f^\delta\in {\mathbb{R}^m}$, where subscript $\delta$ is used to remind that the measured data $f^\delta$ is affected by (additive) noise.

To avoid getting an ill-conditioned matrix like in \eqref{eq:LSE_no_regularizer}, a standard procedure known as \emph{variational regularization} is to define a new cost function $J_R$ by adding a regularizer $\mathcal{R}$ to the cost function $\mathcal{D}$, i.e.,
\begin{align}
 J_{\mathcal{R}}(u) = \lambda \mathcal{D}(Ku, f^\delta) + \mathcal{R}(u) \,,\label{eq:general_regularizer}
\end{align}
where $\mathcal{D}$ is a proper, convex and continuously differentiable function in the first argument, $\mathcal{R}$ is a proper, convex and lower semi-continuous function, and $\lambda > 0$ is the so-called regularization parameter that balances the influence of the data term $J$ and the regularizer $\mathcal{R}$.
One of the most famous examples in the literature is the Tikhonov regularizer $\mathcal{R}(u)=\frac{1}{2}\lVert u \rVert^2$ (cf. \cite{tikhonov1963solution}). For example, Tikhonov regularization with the least squares cost function in \eqref{eq:general_regularizer} leads to
\begin{align}
 u_\lambda = \inv{(\lambda K^\top K + I_n)}K^\top f^\delta,\label{eq:LSE_Tikhonov_regularizer}
\end{align}
where $I_n\in\mathbb{R}^{n\times n}$ is the identity matrix.

More recently, model-based regularizers such as Tikhonov-type regularization have been replaced by data-driven regularizers of the form $\mathcal{R}_{\Psi}(x)$, where the subscript $\Psi$ indicates that $\mathcal{R}_{\Psi}:\mathbb{R}^n \to \mathbb{R} \cup \{\infty\}$ is a parameterized function with parameters $\Psi\in\mathbb{R}^p$. Examples for such data-driven variational regularizations are Markov Random Field priors like the fields-of-expert regularizer \cite{roth2009fields}, which in combination with \eqref{eq:general_regularizer} was popularized for approximating inverse problems solutions in \cite{chen2014insights}:
\begin{align}
 u_{\lambda, \Psi} \in \argmin_{u \in \mathbb{R}^n} \left\{ \frac{\lambda}{2} \| Ku - f^\delta \|^2 + \sum_{i = 1}^r \phi_i(A_i u) \right\} \,.\label{eq:var-reg-with-mrf}
\end{align}
Here $\{ \phi_i \}_{i = 1}^r$ are proper, convex and lower semi-continuous functions and the parameters are $\Psi = \{A_i \}_{i = 1}^r$. Optimal parameters $\Psi$ can be found by training the model for given pairs $(u, f^\delta)$ with different techniques. In this paper we compare two of them, namely \textit{bilevel learning} and \textit{deep equilibrium models}.

\subsection{Deep equilibrium models}\label{sec:deep-equib}
Deep equilibrium models \cite{Bai2019} are deep neural networks that have fixed points, and those fixed points are trained to match data samples from a training data set. Suppose we denote by $G(u, \Psi):\mathbb{R}^n \times \mathbb{R}^p \rightarrow \mathbb{R}^n$ our deep neural network with parameters $\Psi \in \mathbb{R}^p$, then training a deep equilibrium model can be formulated as the constrained optimization problem
\begin{align}
 \min_{\Psi} J(u) \qquad \textnormal{subject to} \qquad u = G(u, \Psi) \,,\label{eq:deep-equilibrium-constrained}
\end{align}
where $J$ is the data fidelity term between the argument $u$ and the corresponding ground truth $u^\dagger$ (e.g., mean squared error, or cross entropy function).
We now want to re-formulate \eqref{eq:deep-equilibrium-constrained} with the help of a Lagrange multiplier $\mu$ to the saddle-point problem $\min_{u, \Psi} \max_{\mu} \mathcal{L}(u, \Psi, \mu)$ with 
$$$$
\begin{align}
 \mathcal{L}(u, \Psi, \mu) = J(u) + \left\langle \mu, u - G(u, \Psi) \right\rangle \,,\label{eq:deep-equilibirum-lagrange}
\end{align}
where $\langle \cdot,\cdot \rangle$ denotes the inner product.
Computing the optimality system of \eqref{eq:deep-equilibirum-lagrange}, i.e., computing the partial derivatives of $\mathcal{L}$ with respect to the individual arguments and setting them to zero, yields the nonlinear system of equations
\begin{subequations}
\begin{align}
 u^\ast &= G(u^\ast, \Psi^\ast) \,, \label{subeq:fixed-point}\\ %
 0 &= \nabla J(u^\ast) + \left(I - (\partial_{u^\ast} G(u^\ast, \Psi^\ast))^\top \right) \mu^\ast \label{subeq:lagrange-fixed-point} \,, \\
 0 &= \partial_{\Psi^\ast} G(u^\ast, \Psi^\ast)^\top \mu^\ast \,, 
\end{align}
\end{subequations}
where $\partial_{u^\ast} G(u^\ast, \Psi^\ast)$ and $\partial_{\Psi^\ast} G(u^\ast, \Psi^\ast)$ denote the Jacobian matrices of $G$ with respect to $u^\ast$ and $\Psi^\ast$, respectively. Note that in order to compute a solution $u^\ast$ of \eqref{subeq:fixed-point} we are required to solve a fixed point problem and further require that such a fixed point exists. In order to compute a solution $\mu^\ast$ to \eqref{subeq:lagrange-fixed-point}, we need to solve the linear system 
\[\mu^\ast = - \left(I - (\partial_{u^\ast} G(u^\ast, \Psi^\ast))^\top \right)^{-1} \nabla J(u^\ast).\] 
In \cite{Bai2019} it was proposed to also formulate this problem as a fixed point problem, i.e., we aim to find $\mu^\ast$ that satisfies
\begin{align*}
 \mu^\ast = (\partial_{u^\ast} G(u^\ast, \Psi^\ast))^\top\mu^\ast -\nabla J(u^\ast) \,. 
\end{align*}
It is important to emphasize that, as for the fixed point $u^\ast$, we obviously require existence of $\mu^\ast$ in order to be able to compute it. Assuming that we can compute both $u^\ast$ and $\mu^\ast$, we can then compute the gradient of $\mathcal{L}$ with respect to the network parameters $\Psi^\ast$, for instance with a gradient-based iterative algorithm like gradient descent,
\begin{align*}
 \Psi^\ast_{j + 1} = \Psi^\ast_{j} - \tau \, \partial_{\Psi^\ast} G(u^\ast_{j}, \Psi^\ast_{j})^\top \mu_{j}^\ast \,,
\end{align*}
where $\{ \Psi^\ast_{j} \}_{j = 1}^\infty$, $\{ u_{j}^\ast \}_{j = 1}^\infty$ and $\{ \mu_{j}^\ast \}_{j = 1}^\infty$ are sequences to approximate $\Psi^\ast$, $u^\ast$ and $\mu^\ast$, and where every $u_{}^\ast$ satisfies $u_{j}^\ast = G(u^\ast_{j}, \Psi^\ast_{j})$, while every $\mu_{j}^\ast$ satisfies $\mu_{j}^\ast = (\partial_{u^\ast} G(u^\ast_{j}, \Psi^\ast_{j}))^\top\mu^\ast_{j} -\nabla J(u^\ast_{j})$. Here, the parameter $\tau$ is a positive step-size parameter that controls the length of the step in the direction of the negative gradient. In practice, any first-order optimization method other than gradient descent can also be used to find optimal parameters $\Psi^\ast$.

In the context of inverse problems of the form \eqref{eq:general_inverse_problem}, we can construct many meaningful deep equilibrium models. 
In \cite{Gilton2021}, one of the methods of consideration is a so-called DeProx-type method, where a neural network is composed with a gradient descent step on the mean-squared error of the forward model output and the measurement data (which in this context is also known as a step of Landweber regularization), i.e.,
\begin{align}
  u^{k + 1} = \mathcal{N}_\Psi\left( u^k - \tau \lambda \, K^\top \left( Ku^k - f^\delta \right) \right) \,,\label{eq:neural-iteration}
\end{align}
for a neural network $\mathcal{N}$ with parameters $\Psi$. Assuming the neural network is $1$-Lipschitz, i.e.,
\begin{align*}
  \| \mathcal{N}_\Psi(u) - \mathcal{N}_\Psi(v) \| \leq \| u - v \| \,,
\end{align*}
for all $u, v \in \mathbb{R}^n$ and choosing $\tau$ such that $\tau < 2/(\lambda \| K \|^2 )$ is satisfied, then we observe that $u^\ast$ with 
\begin{align}
  u^\ast = \mathcal{N}_\Psi\left( u^\ast - \tau \lambda K^\top \left( Ku^\ast - f^\delta \right) \right) \label{eq:fixed_point_equation}
\end{align}
is a fixed point of iteration \eqref{eq:neural-iteration}, which we can conclude from Banach's fixed point theorem if $K^\top K$ has full rank (in which case the fixed point is also unique); if $K^\top K$ does not have full rank, the mapping is only nonexpansive and non-emptiness of the fixed point set has to be verified (cf. \cite{combettes2021fixed}). Hence, we can guarantee convergence of the sequence. The Lipschitz continuity of the mapping $\mathcal{N}_\Psi\left(u - \tau\lambda K^\top \left( Ku - f^\delta \right)\right)$ implies continuity, so we can characterize the fixed point via
\begin{align*}
  u^\ast &= \lim_{k \rightarrow \infty} u^k = \lim_{k \rightarrow \infty} \mathcal{N}_\Psi\left(u^{k - 1} - \tau\lambda K^\top \left( Ku^{k - 1} - f^\delta \right)\right)\\
  &= \mathcal{N}_\Psi\left( \lim_{k \rightarrow \infty } \left(I - \tau \lambda K^\top K\right) u^{k - 1} + \tau\lambda K^\top f^\delta \right) = \mathcal{N}_\Psi\left( \left(I - \tau \lambda K^\top K\right) u^\ast + \tau\lambda K^\top f^\delta \right) \,.
\end{align*}
This means that for the operator $G(u, \Psi):= \mathcal{N}_\Psi\left( u - \tau \lambda K^\top \left( Ku - f^\delta \right) \right)$, the deep equilibrium problem \eqref{eq:deep-equilibrium-constrained} is well-defined if the conditions outlined above are met.

Another method that was considered in \cite{Gilton2021} is the so-called DeGrad-type method, which is motivated by gradient descent where a neural network is used to replace the gradient of a regularization term, i.e.,
\begin{align}
  u^{k + 1} = u^k - \tau \left(\lambda \, K^\top \left( Ku^k - f^\delta \right) + \mathcal{N}_\Psi(u^k) \right) \,.\label{eq:neural-iteration_de_grad}
\end{align}
Making similar assumptions as in the DeProx case, we can assume that a (unique) fixed point %
exists, so that the deep equilibrium problem \eqref{eq:deep-equilibrium-constrained} is well-defined again.

\subsection{Bilevel learning}\label{subsec:bilevel-learning}
Following the description of data-driven variational regularization models such as \eqref{eq:var-reg-with-mrf} in Section~\ref{sec:inverse-problems}, we want to briefly recall how the parameters $\Psi$ of the regularizer can be trained with the help of bilevel learning. In this context, we can formulate the bilevel optimization problem as
\begin{subequations}
\begin{align}
 &\min_{\Psi} J(u^\ast) \,,\label{eq:bilevel-upper-level}
\intertext{subject to}
u^\ast \in &\argmin_{u} \left\{ \frac{\lambda}{2} \| Ku - f^\delta \|^2 + \mathcal{R}_{\Psi}(u) \right\} \,,\label{eq:bilevel-lower-level}
\end{align}\label{eq:bilevel}%
\end{subequations}
where \eqref{eq:bilevel-upper-level} and \eqref{eq:bilevel-lower-level} are the upper-level and the lower-level problems, respectively. Here, $J$ denotes a convex and continuously differentiable loss function, and $\mathcal{R}_\Psi:\R^n \to \R \cup \{ \infty \}$ is a proper, convex and lower semi-continuous function that is parameterized by parameters $\Psi$. Note that due to the convexity of the lower-level problem, $u^\ast$ will always be a global minimizer. 
General existence results of \eqref{eq:bilevel} can for example be found in \cite{DelosReyes2016}. Traditionally, \eqref{eq:bilevel} is solved with the help of gradient-based optimization methods. Historically, in order to compute the gradient of the upper-level optimization problem with regards to $\Psi$ it is common to make use of the implicit function theorem (cf. \cite[Theorem 4.E]{zeidler2012applied} and \cite[Corollary 4.34]{mordukhovich2006variational}). If we assume that $\mathcal{R}_\Psi$ is differentiable with respect to its argument, we can characterize the solution of \eqref{eq:bilevel-lower-level} via the optimality condition
\begin{align*}
 0 = S(u, \Psi) = \lambda K^\top \left( Ku - f^\delta\right) + \nabla \mathcal{R}_\Psi(u) \,.
\end{align*}
If we further assume that $S$ is strictly differentiable and that $\nabla_u S(u, \Psi)$ is invertible, we can compute $\nabla J(u^\ast(\Psi))$ via the implicit function theorem, i.e.
\begin{align*}
 \nabla J(u^\ast(\Psi)) = \left( \nabla_u S(u, \Psi) \right)^{-1} \nabla_\Psi S(u, \Psi) \,.
\end{align*}
Please note that for $S$ to be strictly differentiable, we require $\mathcal{R}_\Psi$ to be twice differentiable with regards to its argument. The invertibility of $\mathcal{R}_\Psi$ can usually be achieved if we make $\mathcal{R}_\Psi$ strictly convex, for example with an additional elliptic regularization term as in \cite{DeLosReyes2013}. 

In the following section, we want to highlight that bilevel optimization can be viewed as a special case of the deep equilbirium model if we replace the lower-level optimization problem with an equivalent fixed-point problem.

\subsection{Bilevel learning as a deep equilibrium model}
Under suitable conditions on $\mathcal{R}_\Psi$ (for example strong convexity), we can conclude uniqueness of $u^\ast$ in \eqref{eq:bilevel-lower-level}. If the solution $u^\ast$ is unique, its computation can also be replaced by a suitable fixed point iteration. To give an example, if $\mathcal{R}_\Psi$ is considered $L$-smooth, i.e. $\mathcal{R}_{\Psi}$ is continuously differentiable and $\nabla \mathcal{R}_\Psi$ is Lipschitz-continuous with Lipschitz constant $L$, one can minimize the lower-level problem \eqref{eq:bilevel-lower-level} via an iterative strategy such as gradient descent, i.e.,
\begin{align}
 u^{k + 1} = u^{k} - \tau \left( \lambda K^\top (K u^{k} - f^\delta ) + \nabla \mathcal{R}_\Psi(u^{k}) \right) \,,\label{eq:lower-level-gradient-descent}
\end{align}
for a sequence $\{ u^{k} \}_{k = 0}^\infty$ with arbitrary initial value $u^{0}$ and a positive step-size parameter~$\tau$. If $\tau$ is chosen appropriately, then convergence of \eqref{eq:lower-level-gradient-descent} to $u^\ast$ can be guaranteed. As a consequence, we can replace \eqref{eq:bilevel-lower-level} with the fixed point constraint
\begin{align}
 u^\ast = u^\ast - \tau \left( \lambda K^\top (K u^\ast - f^\delta ) + \nabla \mathcal{R}_\Psi(u^\ast) \right) \label{eq:lower-level-gradient-descent-fixed-point}
\end{align}
and therefore reformulate \eqref{eq:bilevel} to a saddle-point problem of the form of \eqref{eq:deep-equilibirum-lagrange} with $G(u^\ast, \Psi) = u^\ast - \tau \left( \lambda K^\top (K u^\ast - f^\delta ) + \nabla \mathcal{R}_\Psi(u^\ast) \right)$. Note that, depending on the choice of $\mathcal{R}_\Psi$, instead of gradient descent many other optimization algorithms can be chosen, which leads to a variety of different fixed point equations $u^\ast = G(u^\ast, \Psi)$ with the same fixed point $u^\ast$ that can be chosen in \eqref{eq:deep-equilibirum-lagrange}.

We want to emphasize that every minimization problem $\min_{u} F(u)$ with proper, convex and lower semi-continuous $F$ can na\"{i}vely be turned into a fixed-point iteration of the form
\begin{align}
  u^{k + 1} = \argmin_{u} \left\{ F(u) + \frac{1}{2\tau} \| u - u^k \|^2 \right\} \,,\label{eq:naive-fixed-point}
\end{align}
for a positive step-size parameter $\tau$. Minimizing \eqref{eq:naive-fixed-point} is obviously as difficult as minimizing $F$ itself, so other fixed-point iterations such as \eqref{eq:lower-level-gradient-descent-fixed-point} are usually more suitable in practice if applicable. However, if $F$ has a minimum, then $u^k$ as defined in \eqref{eq:naive-fixed-point} converges to the set of minimizers of $F$ and $F(u^k)$ converges to its optimal value \cite{bauschke2011convex}. This implies that bilevel learning problems with proper, convex and lower semi-continuous lower level problems form a subset of deep equilibrium methods.

Vice versa, every fixed-point problem can na\"{i}vely be converted into a (potentially non-convex) optimization problem by converting a fixed-point equation like $u^\ast = G(u^\ast)$ into $\min_{u} H(u - G(u))$, for a non-negative continuous function $H$ with $H(0) = 0$. An example for $H$ is the squared Euclidean norm, i.e. $H(v) = \| v \|^2 / 2$. In this setting it is obvious that $u^\ast = G(u^\ast)$ is a global minimizer of $H(u - G(u))$, assuming that the fixed-point exists. However, a more interesting question is whether the class of fixed-point operators is strictly larger than the class of fixed-point operators arising from the computational minimization of optimization problems. A thorough mathematical discussion of this question is beyond the scope of this paper, but it is clear that certain vector-fields, such as $C^\top f( A x)$ for a function $f$, a vector $x$ and two matrices $A$ and $C$, cannot be characterized as gradients of functions if $C \neq A$ because the curl of the vector-fields is non-zero, which means the vector fields are not conservative.

\subsection{Why na\"{i}vely learning fixed points does not work}
\label{sec:naively-learning-fixed-point-does-not-work}
We conclude this section by briefly addressing the problem of na\"{i}vely learning fixed points by minimizing empirical risks that measure the deviation between model output and desired output data $u^\dagger$. Suppose we choose the De-Prox architecture \eqref{eq:neural-iteration} from Section~\ref{sec:deep-equib}. Instead of solving \eqref{eq:deep-equilibrium-constrained} with $G$ defined as $G(u, \Psi) = \mathcal{N}_\Psi\left( u - \tau \lambda K^\top \left( Ku - f^\delta \right) \right)$ for $J(u):= \frac12 \| u - u^\dagger \|^2$, we could na\"{i}vely train the parameters $\Psi$ by minimizing $L$ defined as
\begin{align}
  L(u^\dagger, \Psi):= \frac12 \left\| \mathcal{N}_\Psi\left( u^\dagger - \tau \lambda K^\top \left( Ku^\dagger - f^\delta \right) \right) - u^\dagger \right\|^2 \label{eq:naive}
\end{align}
with respect to $\Psi$. Suppose $f^\delta = Ku^\dagger + \delta$, for some perturbation $\delta$, then Problem \eqref{eq:naive} is the same as minimizing an empirical risk for a denoising autoencoder, i.e., 
\begin{align}
 \min_{\Psi} \frac12 \left\| \mathcal{N}_\Psi\left( u^\dagger + \tau \lambda K^\top \delta \right) - u^\dagger \right\|^2 \,. \label{eq:denoising-aec}
\end{align}
The problem with \eqref{eq:denoising-aec} is that the network $\mathcal{N}_\Psi$ is only trained to remove the term $\tau \lambda K^\top \delta$. However, if $K$ is underdetermined (like in the inpainting application later), then $\mathcal{N}_\Psi$ in iteration \eqref{eq:neural-iteration} also has to be able to combat any lack of information that stems from the underdeterminedness of $K$ (for initial values other than $u^\dagger$). For that reason, solving \eqref{eq:deep-equilibrium-constrained} is superior over \eqref{eq:denoising-aec} in this particular context.

Similar considerations can be drawn if we choose the De-Grad architecture \eqref{eq:neural-iteration_de_grad}.
If we assume that the sequence $\{u^k\}_{k=0}^{\infty}$ converges to the fixed point $u^*$, and that the reconstruction is perfect (i.e., $u^* = u^\dagger$), then we can train the parameters $\Psi$ by minimizing $L$ defined as
\begin{align}
 L(u^\dagger, \Psi):= \frac12 \left\| \lambda K^\top ( Ku^\dagger - f^\delta) + \mathcal{N}_\Psi(u^\dagger) \right\|^2 \label{eq:naive-degrad}
\end{align}
with respect to $\Psi$.
If we again assume $f^\delta = Ku^\dagger + \delta$ for some perturbation $\delta$, then Problem \eqref{eq:naive-degrad} is equivalent to 
\begin{align}
 \min_{\Psi} \frac12 \left\| \mathcal{N}_\Psi( u^\dagger) - \lambda K^\top \delta \right\|^2 \,, \label{eq:denoising-aec-degrad}
\end{align}
which means that in this case the network $\mathcal{N}_\Psi$ is trained to approximate $\lambda K^\top \delta$. Again, if $K$ is underdetermined, then iteration \eqref{eq:neural-iteration_de_grad} has to compensate for the lack of information that stems from the underdeterminedness of $K$, which is not the case when \eqref{eq:denoising-aec-degrad} is solved instead. In Section~\ref{subsec:naively-learning-fixed-point} we will train a model by minimizing \eqref{eq:denoising-aec-degrad} and empirically verify that the results obtained with this model are not as good as the ones we get when we minimize \eqref{eq:deep-equilibrium-constrained} instead.

\section{Architecture design and implementation}\label{sec:arch_design_and_implementation}
We compare deep equilibrium and bilevel learning methods empirically for three different inverse problems. In general, it is very difficult to have a fair comparison between different architectures, as their expressivity is not only determined by the number of parameters or the choice of activation functions, but the complex interplay of all architecture design choices as well as other factors like the optimization method that is used for training the models, or the choice of loss and regularization functions, to name only a few. 
We have therefore decided to choose one deep equilibrium model that belongs to the class of deep equilibrium gradient descent, cf. \cite[Section 3.1]{Gilton2021}. We compare this model with a bilevel approach where the lower-level problem is computed via gradient descent, and both architectures are chosen to be as similar as possible. The individual models are described in detail in the following sections. 

\subsection{Deep equilibrium models}
In this section we describe the specific deep equilibrium models that we use for the comparison of numerical results in Section~\ref{sec:numerical-results}. For a fairer comparison with bilevel optimization methods, we focus on models of type DeGrad as described in Section~\ref{sec:deep-equib}. 

\subsubsection{Deep equilibrium gradient descent}\label{subsec:de-grad}
Following \cite[Section 3.1]{Gilton2021}, we define the deep equilibrium gradient descent method for the approximation of inverse problems solutions as
\begin{align}
 u^{k+1} = u^{k} - \tau \left( \lambda K^\top (K u^{k} - f^\delta ) + \mathcal{N}_\Psi(u^{k}) \right) \label{eq:deq-gd} \,,
\end{align}
where $\mathcal{N}_\Psi:\mathbb{R}^n \rightarrow \mathbb{R}^n$ is a neural network with parameters $\Psi$. For the comparison with bilevel optimization methods, we will choose $\mathcal{N}_\Psi$ as
\begin{align}
 \mathcal{N}_\Psi(u) = \gamma \, C^\top \sigma( A u + b ) \,,\label{eq:deq-non-gradient}
\end{align}
for parameters $\Psi = (A, C, b)$ with matrices 
$A, C \in \R^{q \times r}$ and bias vector 
$b \in \mathbb{R}^q$, activation function 
$\sigma:\mathbb{R}^q \rightarrow \mathbb{R}^q$ and positive constant $\gamma$. Please note that \eqref{eq:deq-non-gradient} is not necessarily a gradient of a function with argument $u$, unless $C = A$ and suitable choices of $\sigma$. For consistency in the numerical examples below, we will choose the same activation functions discussed in Section~\ref{subsec:bilevel-optimization-models}.

\subsection{Bilevel learning models}
\label{subsec:bilevel-optimization-models}
In order to effectively compare bilevel optimization with deep equilibrium models, we replace the solution of \eqref{eq:bilevel-lower-level} with the fixed point of \eqref{eq:lower-level-gradient-descent-fixed-point}. We choose $\mathcal{R}_\Psi(u)$ to be of the form
\begin{align}
 \mathcal{R}_\Psi(u) = \gamma \inf_v \left( \frac12 \| v - Au - b \|^2 + \mathcal{R}(v) \right) \,,
\end{align}
which is the Moreau--Yosida regularization \cite{moreau1962fonctions,yosida1964functional} of the proper, convex and lower semi-continuous function $\mathcal{R}:\mathbb{R}^{s} \rightarrow \mathbb{R} \cup \{ \infty \}$ composed with the affine-linear transformation $A \cdot + b$. Here $\Psi = (A, b)$ denote the model parameters, which are a matrix $A \in \mathbb{R}^{s \times n}$ and a bias vector $b \in \mathbb{R}^s$, and $\gamma$ is a positive parameter. The gradient of $\mathcal{R}_\Psi$ with respect to argument $u$ reads
\begin{align}
 \nabla \mathcal{R}_\Psi(u) = \gamma A^\top \left( Au + b - \text{prox}_\mathcal{R}\left(Au + b\right) \right) = \gamma A^\top \text{prox}_{\mathcal{R}^\ast}\left( Au + b \right) \,,\label{eq:gradient-of-regularization-function}
\end{align}
cf.\ \cite[Proposition 12.30]{bauschke2011convex}, where $\text{prox}_\mathcal{R}:\mathbb{R}^s \rightarrow \mathbb{R}^s$ denotes the proximal map \cite{moreau1965proximite} of $\mathcal{R}$, i.e.,
\begin{align*}
 \text{prox}_\mathcal{R}(w) = \argmin_{v \in \mathbb{R}^s} \left\{ \frac12 \| v - w \|^2 + \mathcal{R}(v) \right\} \,,
\end{align*}
and $\text{prox}_{\mathcal{R}^\ast}$ denotes the proximal map of the convex conjugate $\mathcal{R}^\ast$ of $\mathcal{R}$ (cf. \cite{beck2017first}) that is defined as
\begin{align*}
 \mathcal{R}^\ast(p):= \sup_{v} \left( \langle v, p \rangle - \mathcal{R}(v) \right) \,.
\end{align*}
Note that \eqref{eq:gradient-of-regularization-function} is equivalent to \eqref{eq:deq-non-gradient} if $C = A$ and if $\sigma = \text{prox}_{\mathcal{R}^\ast}$ for some proper, convex and lower semi-continuous $\mathcal{R}$ with conjugate $\mathcal{R}^\ast$. 

\noindent Suppose we define $\chi_C$ as the characteristic function over the convex set $C$, i.e.,
\begin{align*}
 \chi_C(u):= \begin{cases}
  0 & u \in C \\ \infty & u \not\in C
 \end{cases} \,.
\end{align*}
If we choose $\mathcal{R}(u) = \chi_{C_1}(u)$ with $C_1 = \{ 0 \}$, we observe $\mathcal{R}_\Psi(u) = \frac12 \| Au + b \|^2$ and $\mathcal{R}^\ast(p) = 0$. The proximal map of $\mathcal{R}^\ast$ then simply reduces to the identity, i.e., $\text{prox}_{\mathcal{R}^\ast}(w) = w$, and \eqref{eq:gradient-of-regularization-function} reduces to $\nabla \mathcal{R}_\Psi(u) = \gamma A^\top (Au + b)$. 

Another interesting example is $\mathcal{R}(u) = \chi_{C_2}(u)$ with $C_2 = (-\infty, 0]^s$, where we recover the elementwise Rectified Linear Unit (ReLU) \cite{nair2010rectified} activation function $\sigma(x)=\text{ReLU}(x)$ with
\begin{align*}
 (\text{ReLU}(x))_i:= 
 \begin{cases} 
 x_i & x_i \ge 0 \\ 
 0 & x_i < 0
 \end{cases} \,, 
\end{align*}
for the proximal map of the conjugate. In this case, \eqref{eq:gradient-of-regularization-function} reduces to $\nabla \mathcal{R}_\Psi(u) = \gamma A^\top \text{ReLU}(A u + b)$.

If we choose $\mathcal{R}(v) = \varepsilon \| v \|_1 = \varepsilon \sum_{j = 1}^m | v_j |$, we have $\mathcal{R}^\ast(p) = \chi_{C_3}(u)$ with $C_3 = \left\{ u \, | \, \| u \|_\infty \leq \epsilon\right\}$ and
\begin{align*}
 \text{prox}_{\mathcal{R}^\ast}(w)_j = \begin{cases} \varepsilon & w_j > \varepsilon \\ w_j & |w_j| \leq \varepsilon \\ -\varepsilon & w_j < -\varepsilon \end{cases} \,,
\end{align*}
for all $j \in \{1, \ldots, s\}$. Then \eqref{eq:gradient-of-regularization-function} reduces to
\begin{align*}
 (\nabla \mathcal{R}_\Psi(u))_j = \gamma \left( A^\top v\right)_j \quad \text{with} \quad v_i:= \begin{cases} \varepsilon & (A u + b)_i > \varepsilon \\ (Au + b)_i & |(Au + b)_i| \leq \varepsilon \\ -\varepsilon & (Au + b)_i < -\varepsilon \end{cases} \,, 
\end{align*}
for $i \in \{1, \ldots, s\}$ and $j \in \{ 1, \ldots, n\}$.
If, in return, we choose $\mathcal{R}(u) = \chi_{C_3}(u)$, we recover the elementwise soft-shrinkage activation function $\text{prox}_{\mathcal{R}^\ast}(x)=\text{Softshrink}_\varepsilon(x)$ with
\begin{align*}
 (\text{Softshrink}_\varepsilon(x))_i:= 
 \begin{cases} 
  x_i+\varepsilon & x_i < -\varepsilon \\ 
  0  & |x_i| \le \varepsilon \\ 
 x_i-\varepsilon & x_i > \varepsilon
 \end{cases} \,,
\end{align*}
in which case \eqref{eq:gradient-of-regularization-function} reduces to $\nabla \mathcal{R}_\Psi(u) = \gamma A^\top \text{Softshrink}_{\varepsilon}(Au +b)$.
Finally, if we choose $\mathcal{R}$ such that $\mathcal{R}^\ast$ is the following
\begin{align*}
 \mathcal{R}^\ast(w) =
 \begin{cases} 
 w \, \inv{\text{tanh}}(w) + \frac{1}{2} \left( \text{log}(1-w^2)-w^2 \right) & |w|<1 \\ 
 \infty & |w|\ge 1 
 \end{cases} \,,
\end{align*}
we get $\text{prox}_{\mathcal{R}^\ast}(x)=\text{tanh}(x)$ with
$$(\text{tanh}(x))_i = \frac{e^{x_i}-e^{-x_i}}{e^{x_i}+e^{-x_i}} \,.$$
In this case, \eqref{eq:gradient-of-regularization-function} reduces to $\nabla \mathcal{R}_\Psi(u) = \gamma A^\top \text{tanh}(Au +b)$.

\subsection{Inverse problems}\label{sec:applications}

We compare deep equilibrium and bilevel learning for the three different inverse problems of denoising, inpainting and deblurring, which we briefly describe in the following subsections.

\subsubsection{Denoising}
In the denoising task we assume measured data $f$ is affected by additive noise $\delta$ only. This means that $K=I_n$, where $I_n \in \{0, 1\}^{n \times n}$ is the identity matrix of size $n$, so equation \eqref{eq:general_inverse_problem} reads
\begin{align*}
 u+\delta=f^\delta \,.
\end{align*}\label{eq:denoising}%
It is also clear that $u$ and $f^\delta$ have the same dimension in this setting. The main goal is to remove the noise $\delta$ to retrieve $u$ from $f^\delta$.

\subsubsection{Inpainting}
In (image) inpainting we have a partial observation $f$ of unknown, complete data $u$. For example, an image can have missing pixels when a user may want to remove an undesired object from said image. The task of filling missing pixels is called inpainting, and it can be modeled as an inverse problem \cite{schonlieb2015partial}. The matrix $K$ can be modeled by removing every row that corresponds to a missing pixel from an identity matrix $I_n$. If we want $u$ and $f^\delta$ to have the same dimension (i.e., $m=n$), it suffices to substitute zero instead of one in each row that corresponds to a missing pixel. 

\subsubsection{Deblurring}
In (image) deblurring the goal is to recover an unknown image $u$ from a blurred and usually noisy image $f^\delta$. This process can mathematically be described with the inversion of a convolution operator, which without regularization is highly ill-conditioned. 
If we want to apply convolutions to images, we can consider $U\in\mathbb{R}^{n_{\text{row}}\times n_{\text{col}}}$, where 
$n_{\text{row}} * n_{\text{col}} = n$, and 
$\: n_{\text{row}}, n_{\text{col}}\in\mathbb{N}$ 
. We consider $u\in\mathbb{R}^n$ in our notation, which we obtain by turning $U$ into a column vector consisting of all columns of $U$ stacked in consecutive order.

A two-dimensional convolution operator $K_{\Omega}$ is then defined as $K_{\Omega} U = F^\delta_{\text{img}}$, with
$[F^\delta_{\text{img}}]_{h,k} = \sum_{i=-a}^a \sum_{j=-b}^b [\Omega]_{i,j} [U]_{h-i,k-j}$ for all $h,k$, where $\Omega$ is a fixed \emph{convolution kernel}. The components of $\Omega$ are usually normalized to guarantee $\sum_{i} \sum_{j} [\Omega]_{i,j} = 1$.

\subsection{Implementation}
\label{sec:implementation}
For the bilevel learning, we want to obtain the fixed point solution $u^*$ defined in \eqref{eq:lower-level-gradient-descent-fixed-point}. We do this by initializing $u_0$ as the zero vector without loss of generality, and then we iteratively use \eqref{eq:lower-level-gradient-descent}. Note that any initial value for $u_0$ can be chosen because of convergence guarantees to the unique fixed point $u^*$. %

 The computation of the fixed point through \eqref{eq:lower-level-gradient-descent-fixed-point} is computed until at least one of two stopping conditions is satisfied: either the relative tolerance $\lVert u^{k}-u^{k-1} \rVert/ \lVert u^{k}\rVert$ is below a threshold, or the maximum number of iterations is reached. We choose $10^{-3}$ as threshold, and $500$ as maximum number of iterations for the MNIST dataset because we experimentally discover that we can get a good trade-off between the quality of the reconstructed image and the computation time with this setting.
For similar reasons, we choose $10^{-14}$ as threshold, and $1,000$ as maximum number of iterations for the CelebA dataset.
To accelerate the convergence towards the fixed point solution, Anderson acceleration \cite{walker2011anderson} can be implemented following the approach shown in \cite{Gilton2021}. However, in this paper we decided not to use it as it could lead to an ill-conditioned problem if the parameters are not carefully chosen.

For deep equilibrium regularizers, we solve \eqref{eq:deq-gd}, where $\mathcal{N}_\Psi$ works as a parametrized regularizer defined in \eqref{eq:deq-non-gradient}. This means that the neural network we are using is composed of one affine input layer with activation function $\sigma$, and one linear operator as output layer with no bias nor activation function, scaled by a scalar $\gamma$ which we consider as a hyperparameter chosen a-priori (i.e., it is not trained).

For bilevel learning models, we use the same setting with the additional constraint $C=A$ for the reasons explained in Section~\ref{subsec:bilevel-optimization-models}. To have a fair comparison, parameters in $C$ are initialized to be equal to $A$ for both the bilevel model and the deep equilibrium one, although during training only the bilevel model is constrained to satisfy $C=A$.

We want to remark that we choose $\nabla\mathcal{R}_\Psi$ to be the composition of one fully connected affine layer, an activation function, and a linear operator. This is done following Section~\ref{subsec:bilevel-optimization-models} to guarantee a fair comparison between the performance of bilevel learning and DEQ models. More expressive architectures can be used to achieve better results, as shown in \cite{Bai2019, Gilton2021, heaton2021feasibility}.

The chosen optimizer for training parameters $\Psi$ is Adam \cite{kingma2014adam} with initial learning rate $10^{-3}$. 

\section{Numerical results}\label{sec:numerical-results}

We compare bilevel optimization method and deep equilibrium model in our numerical experiments\footnote{Our code is available at \url{https://github.com/reacho/deep-equilibrium-vs-bilevel}}. In particular, we want to experimentally verify whether one method works better than the other in terms of test error and robustness to tuning parameters. %
The model used for the bilevel optimization examples is defined in \eqref{eq:lower-level-gradient-descent}. For deep equilibrium models, we use \eqref{eq:deq-gd}-\eqref{eq:deq-non-gradient}.

We use the MNIST dataset \cite{lecun1998mnist}, a dataset of digits that contains 70,000 grayscale images (60,000 for training, and 10,000 testing) of size $28\times28$, and CelebA \cite{liu2015faceattributes}, a dataset of 202,599 RGB face images of size $178\times218$. Given the larger amount of images in the CelebA dataset and their larger size with respect to MNIST, we randomly select 10 images for training and 10 for testing. For simplicity, we also decide to convert CelebA images to grayscale. 
We rescale pixel values to be within the range $[-1,1]$.
Those images are the unknown $u$ we want to retrieve in \eqref{eq:general_inverse_problem__with_additive_noise}. 

We now describe the settings for the MNIST dataset, for which a broad analysis of hyperparameter choices has been performed. The selection of settings for the CelebA dataset is detailed in subsection~\ref{subsec:comparison-celeba}.
To generate the input $f^\delta$, operator $K$ (either denoising, inpainting, or deblurring described in Section~\ref{sec:applications}) is applied to $u$. To avoid committing an inverse crime, instances $\delta$ of Gaussian random variables are added, where $\delta_i \sim\mathcal{N}(0, \alpha^2)$, and $\alpha$ is the noise level hyperparameter. We use 
$\alpha=0.0,\,0.05,\,0.5$ for noise levels in training, and we keep the same value of $\alpha$ for the test dataset. 
We note that instances $\delta$ of Gaussian random variables are not generated once for each image; instead we generate new values of $\delta$ at each training epoch. 

For the regularizers, the hyperparameters are set to the values %
$\tau=0.01,\,0.1,\,0.9,\,1.1,\,2.1$, and $\gamma=0.1,\,0.5,\,1.0$.
We initialize $A$ and $C$ as square matrices with dimension $784\times784$. It follows that the bias $b$ has dimension $784$. This is not the only possible choice since it is sufficient that $A$ and $C$ have the same dimension (e.g., they can be fully-connected affine layers with rectangular matrices). 
We choose the activation function $\sigma$ as $\text{ReLU}$, $\text{Softshrink}$, or $\text{identity}$, as discussed in Section~\ref{subsec:bilevel-optimization-models}. The value for the threshold of the $\text{Softshrink}$ activation function is chosen as $\tau$, although in principle they are different parameters and can be chosen differently.

We choose the mean-squared error (MSE) as the loss function to minimize the Euclidean distance between the original images $u^\dagger$ and the reconstructed images $u^*$.

The simulations are coded using Python (v3.10.7), and in particular the PyTorch (v1.12.1) library. All the simulations run on the GPU, which allows remarkable speed-up in terms of computation time compared to running the same experiments on the CPU. 
The code has been written and tested on a Windows 10 laptop, and the simulations were run on the High Performance Cluster system Apocrita~\cite{king_thomas_2017_438045}.

\subsection{Denoising}

\begin{figure}[!h]
 \begin{minipage}{1.0\textwidth}
  \centering
  \includegraphics[width=0.5\linewidth]{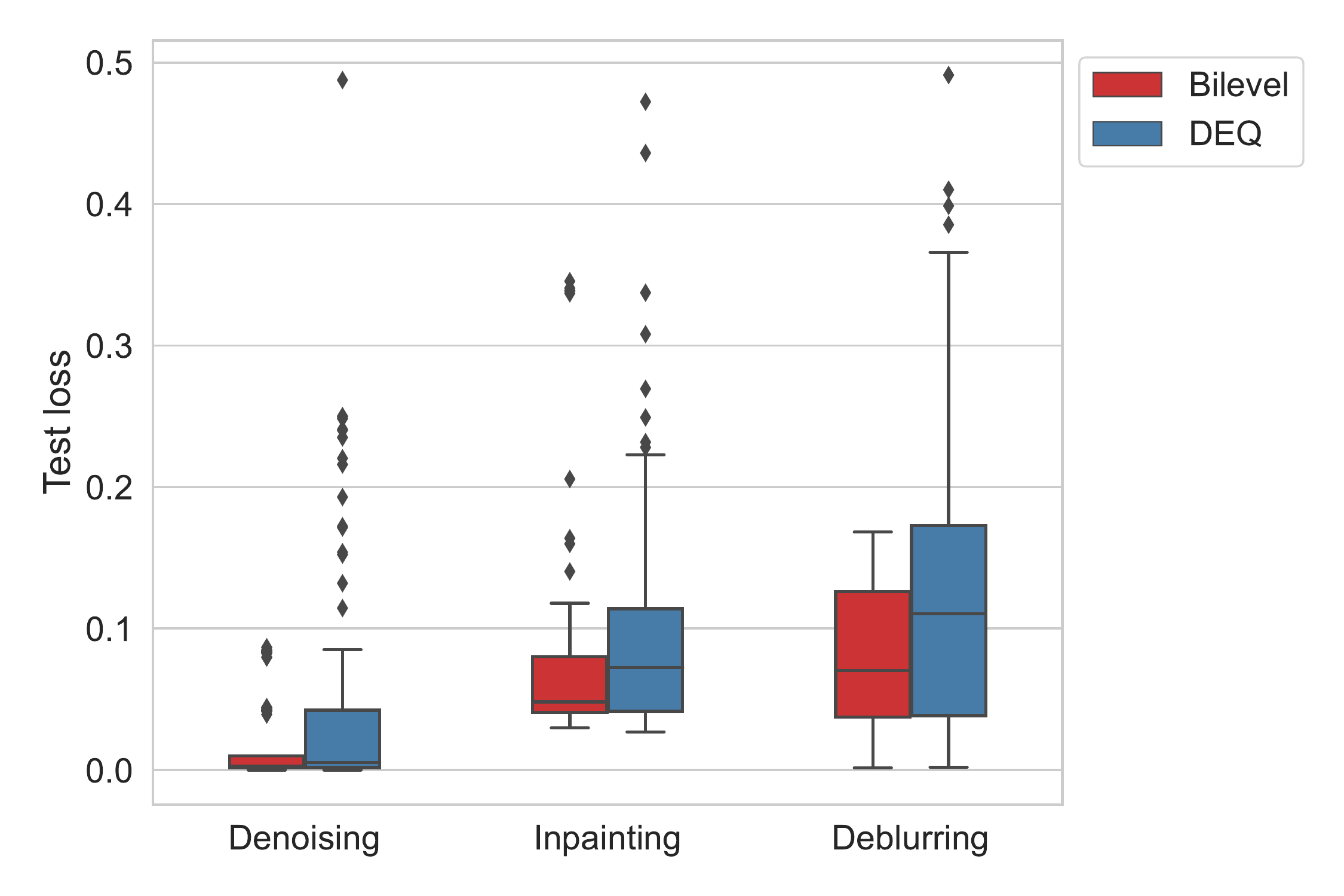}
 \end{minipage}
 \caption{
 Comparison between bilevel optimization and deep equilibrium models for each of the three considered inverse problems, namely denoising, inpainting, and deblurring, over all the range of possible parameters. These boxplots consider the loss of the trained models evaluated on the test dataset. We removed all the configurations with a final loss larger than $0.5$, a value we arbitrarily chose by looking for an empirical relation between the loss and the image quality.
 }
 \label{fig:boxplot-comparison-tasks}
\end{figure}

We start analyzing the denoising task. As it can be expected, our experiments suggest that the denoising task is the easiest among the three tasks considered. Indeed, we see from Fig.\ref{fig:boxplot-comparison-tasks} that the average loss of the trained model is lower in the denoising task, and the interquantile range is smaller too.

\begin{figure}[!h]
 \begin{minipage}{0.49\textwidth}
  \centering
  \includegraphics[width=0.9\linewidth]{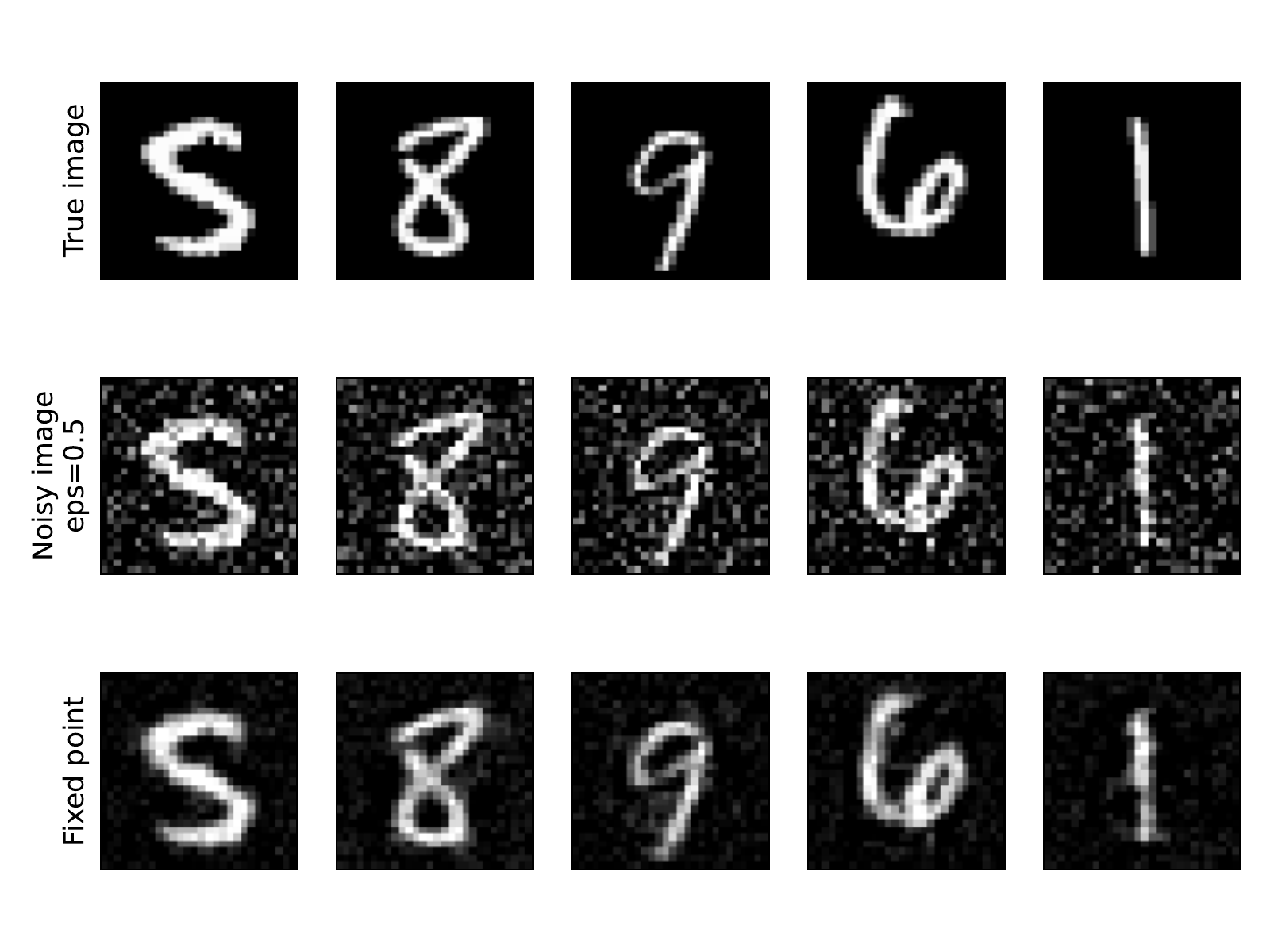}
 \end{minipage}\hfill
 \begin{minipage}{0.49\textwidth}
  \centering
  \includegraphics[width=0.9\linewidth]{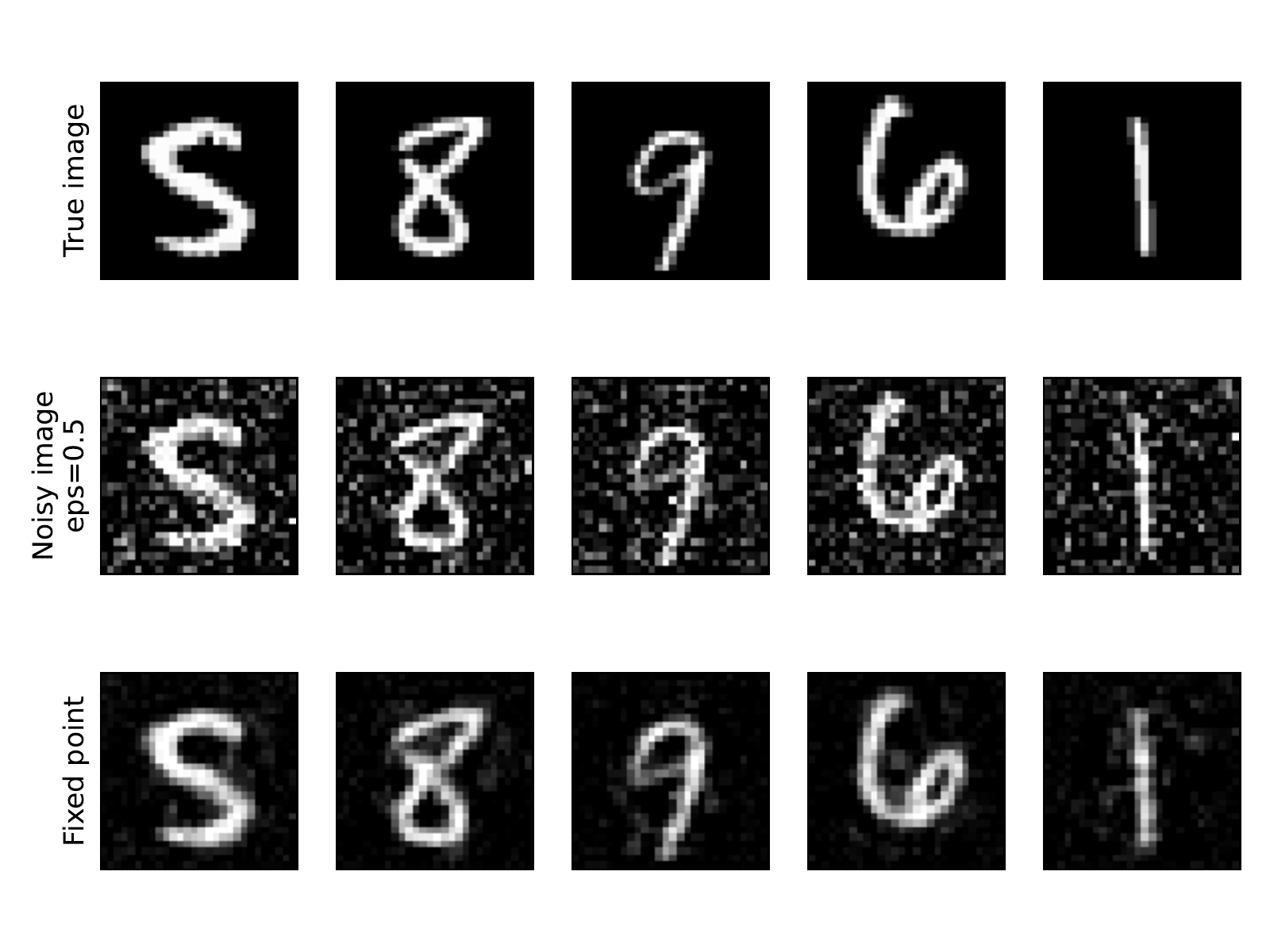}
 \end{minipage}
 \caption{Denoising the MNIST dataset. Visual comparison between bilevel method (left) and deep equilibrium model (right), with parameters $\tau=0.5$, $\gamma=0.1$, and $\sigma=\text{(ReLU)}$. Images are taken from the test dataset. The first row shows the original images; the second row is the model input. The last row is the output of the trained models.}
 \label{fig:denoising__bilevel_vs_DEQ_results_final}
\end{figure}

The visual comparison between bilevel optimization method and deep equilibrium model is shown in Fig.\ref{fig:denoising__bilevel_vs_DEQ_results_final}. The original image (without noise) can be recovered by both methods, and there are no major differences between the reconstructions of the two methods. All the models that achieve a loss smaller than $0.5$ achieve similar visual results.

\subsection{Inpainting}
For inpainting, we choose to mask one third of the image rows, starting from the top and rounding up to the nearest integer. For MNIST images, this corresponds to masking $10$ rows out of $28$.
Note that this task is harder than recovering the same number of pixels if they were randomly selected. In this latter scenario, the values of each missing pixel can be reasonably guessed by looking at the values of its neighbors, which is not true in our setting because most of the missing pixels are surrounded by other missing pixels.
Only $117$ hyperparameter-settings out of $225$ worked for the bilevel model. Deep equilibrium models are more sensitive to the hyperparameter selection than bilevel methods for the inpainting task, since only $73$ hyperparameters settings out of $225$ lead to satisfactory results. We see from Fig.\ref{fig:boxplot-comparison-tasks} that the chosen bilevel methods seem to perform better than their deep equilibrium counterparts in terms of loss minimization.

\begin{figure}[!h]
 \begin{minipage}{0.49\textwidth}
  \centering
  \includegraphics[width=0.9\linewidth]{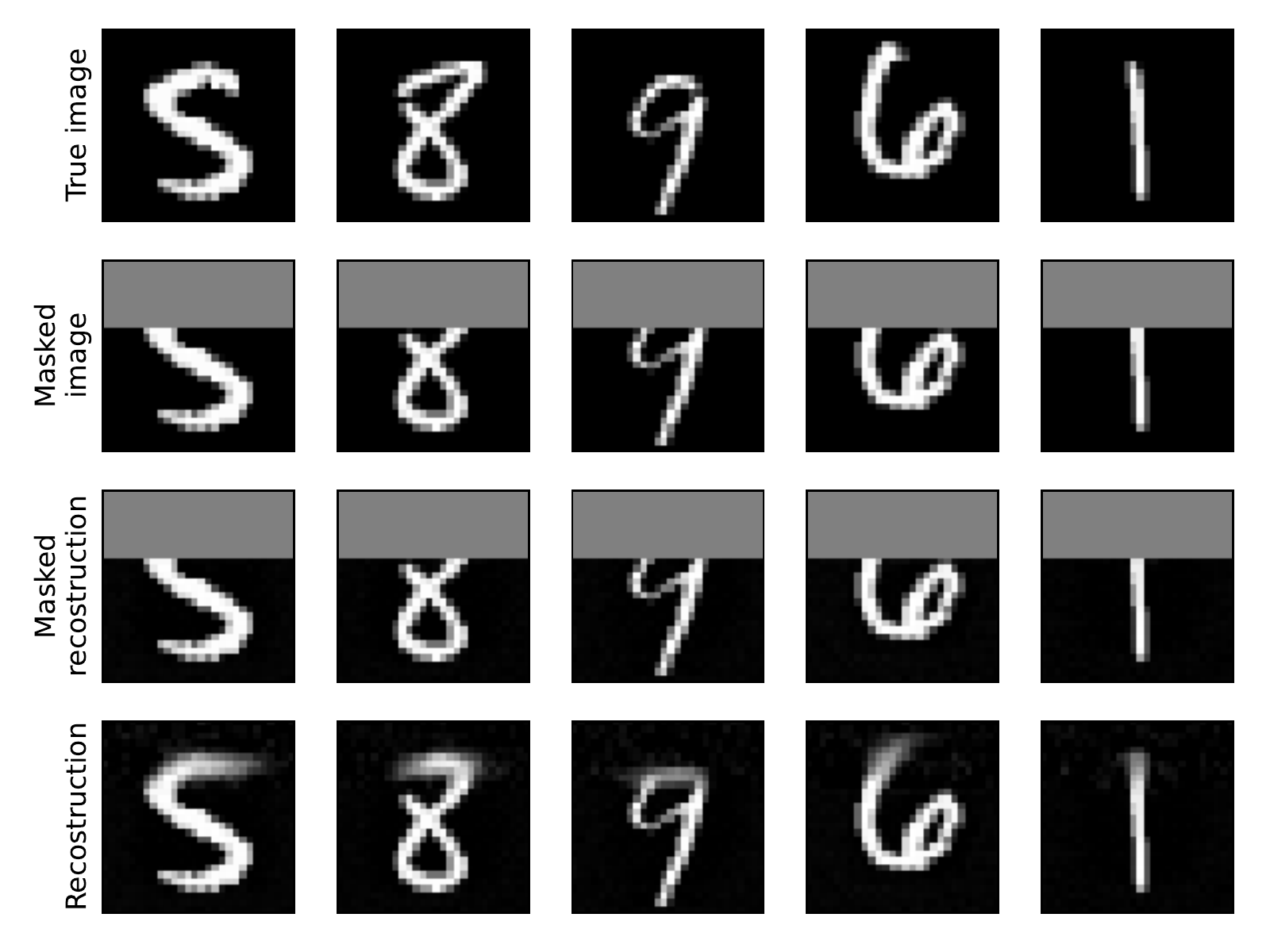}
 \end{minipage}\hfill
 \begin{minipage}{0.49\textwidth}
  \centering
  \includegraphics[width=0.9\linewidth]{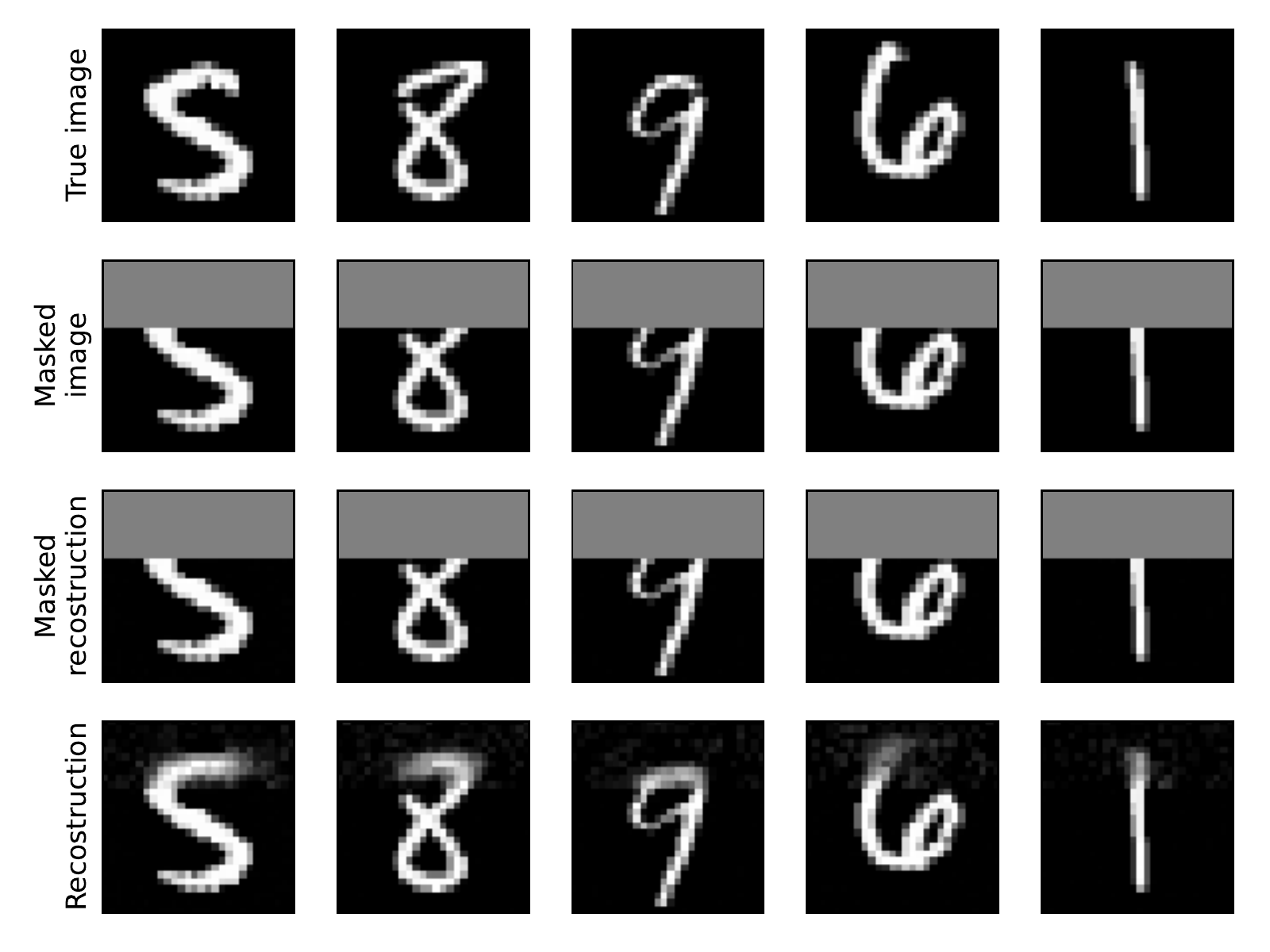}
 \end{minipage}
 \caption{Inpainting MNIST. Comparison between bilevel method (left) and deep equilibrium model (right), with parameters $\tau=0.5$, $\gamma=1.0$, and $\sigma=\text{(Softshrink)}$. Images are taken from the test dataset. The first row shows the original image, the second row is the masked image, i.e., the input of the algorithm. The fourth row is the output of the trained models. Finally, the third row shows what happens when we apply the inpainting operator on the output. The fourth row is the output of the trained deep equilibrium optimization problem. 
 Ideally, the difference between the second and third row should be small.}
 \label{fig:inpainting__train_and_test_results_final}
\end{figure}

Visual results are shown in
Fig.\ref{fig:inpainting__train_and_test_results_final}. Interestingly, both methods achieve good images reconstructions. 

\subsection{Deblurring}

\begin{figure}[!h]
 \begin{minipage}{0.49\textwidth}
  \centering
  \includegraphics[width=0.9\linewidth]{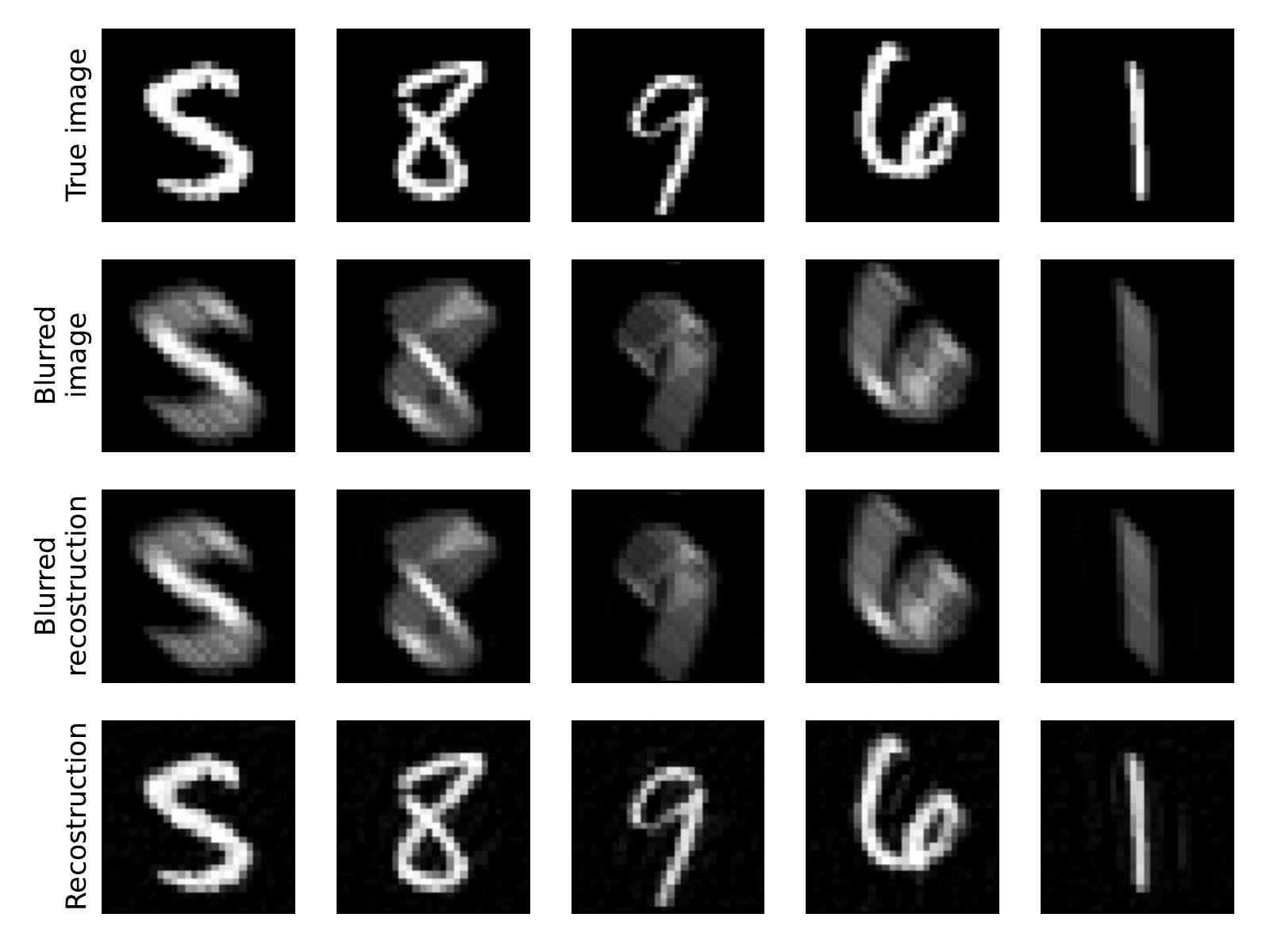}
 \end{minipage}\hfill
 \begin{minipage}{0.49\textwidth}
  \centering
  \includegraphics[width=0.9\linewidth]{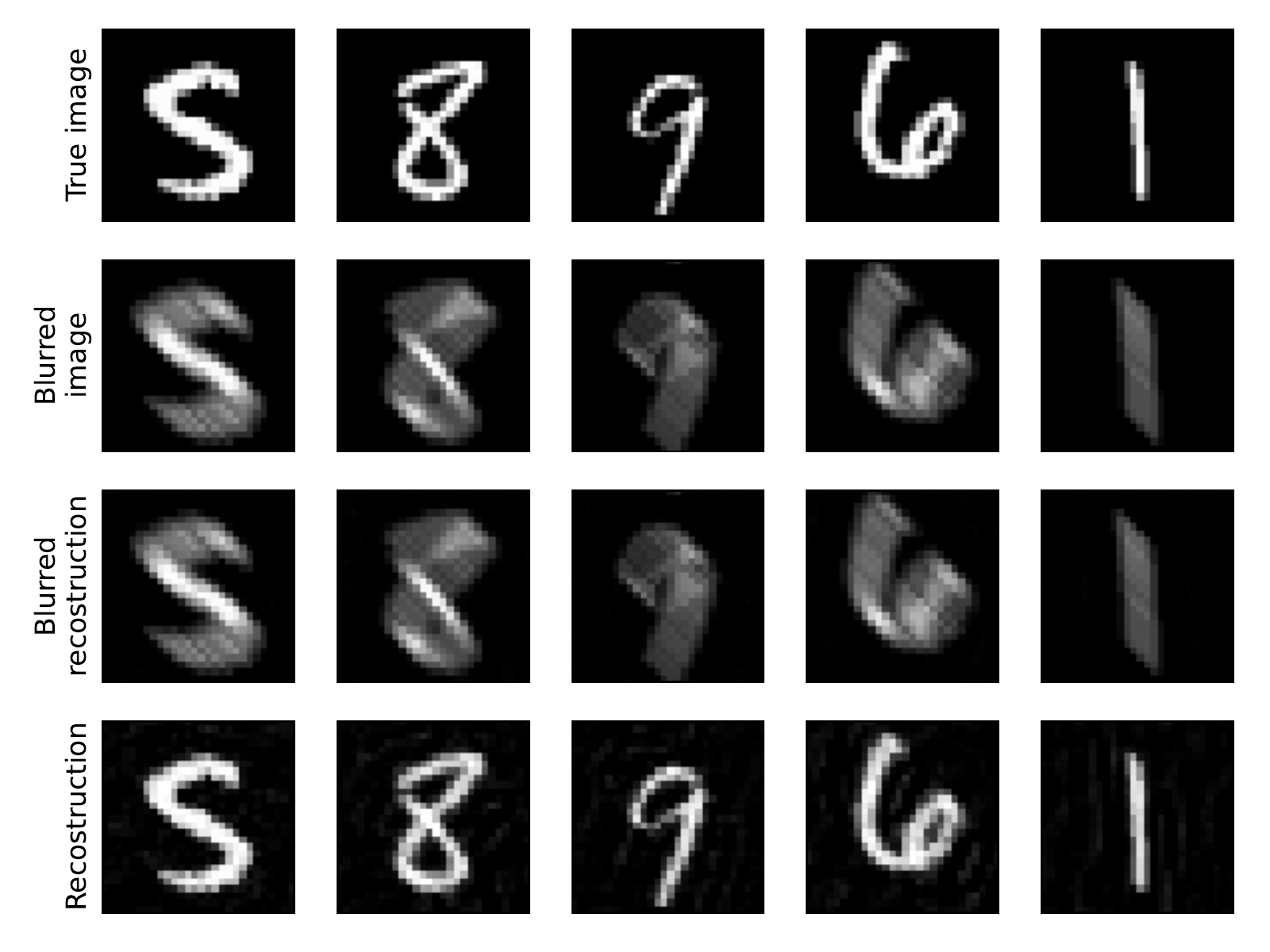}
 \end{minipage}
 \caption{Deblurring MNIST. Comparison between bilevel method (left) and deep equilibrium model (right), with parameters $\tau=0.5$, $\gamma=0.5$, and $\sigma=\text{(Softshrink)}$. Images are taken from the test dataset. The first row shows the original images; the second row is the model input. The last row is the output of the trained models. The third row shows the model output after we apply the convolution kernel to it. Ideally, the difference between the second and the third rows should be small.}
 \label{fig:deblurring__train_and_test_results}
\end{figure}

The last task we consider is deblurring.
We model the convolution operator to mimic a diagonal motion blur, for which we choose the convolution kernel $\Omega$ as $\Omega=\frac{1}{5}I_5$, where $I_5$ denotes the $5\times5$ identity matrix. 
This task has the largest variability in terms of loss results for trained models, as shown in Fig.\ref{fig:boxplot-comparison-tasks}. We also see that bilevel methods achieve a smaller average loss in comparison to deep equilibrium models, with also a smaller interquantile range.
The quality of the retrieved images is similar for the two methods as shown in Fig.\ref{fig:deblurring__train_and_test_results}, even though the average loss seems to suggest a different result (this is because the mean-squared error is not always a good indicator of \textit{similarity} between a pair of images).

\subsection{Sensitivity against noise level}
\begin{figure}[!h]
 \begin{minipage}{1.0\textwidth}
  \centering
  \includegraphics[width=0.32\linewidth]{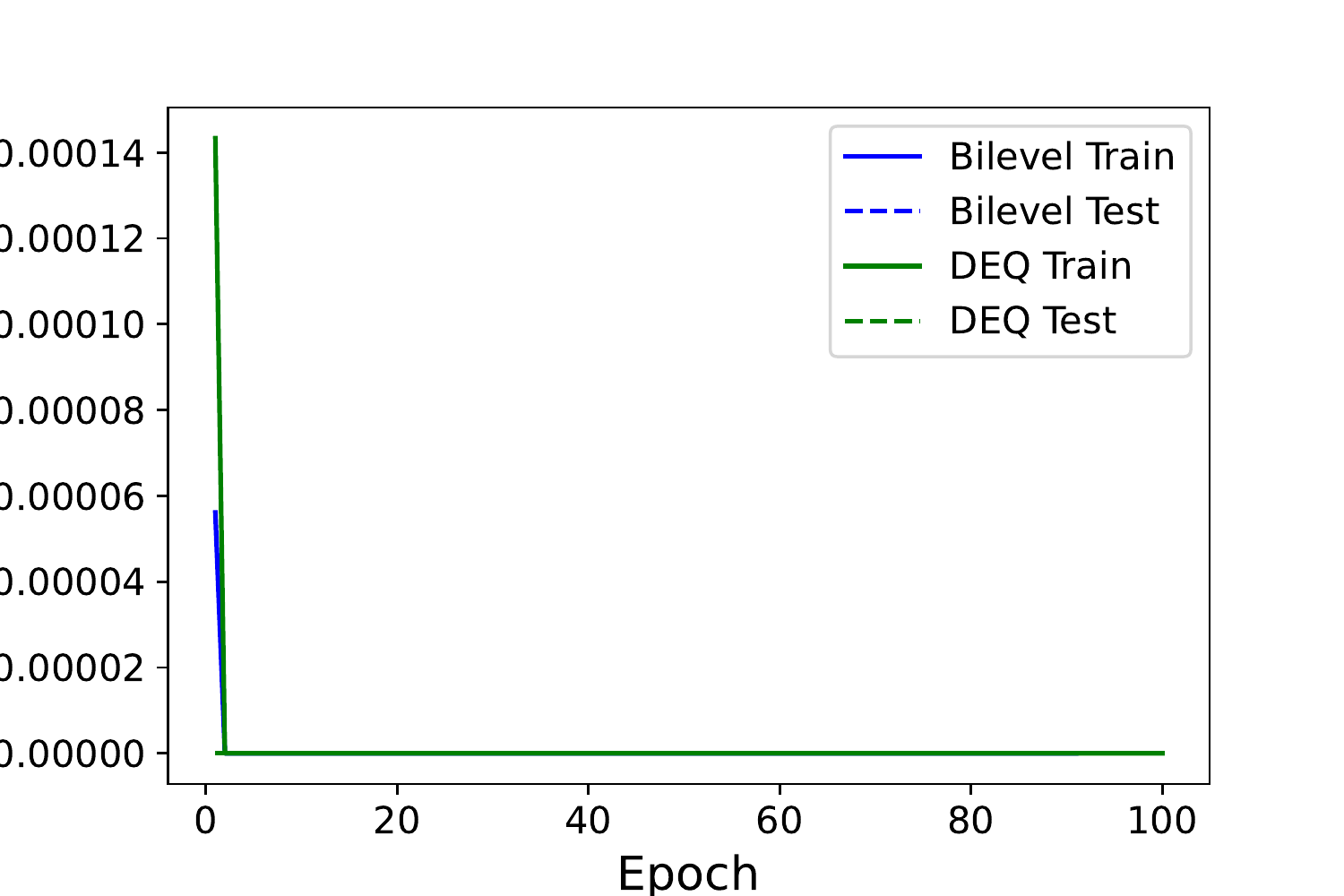}
  \includegraphics[width=0.32\linewidth]{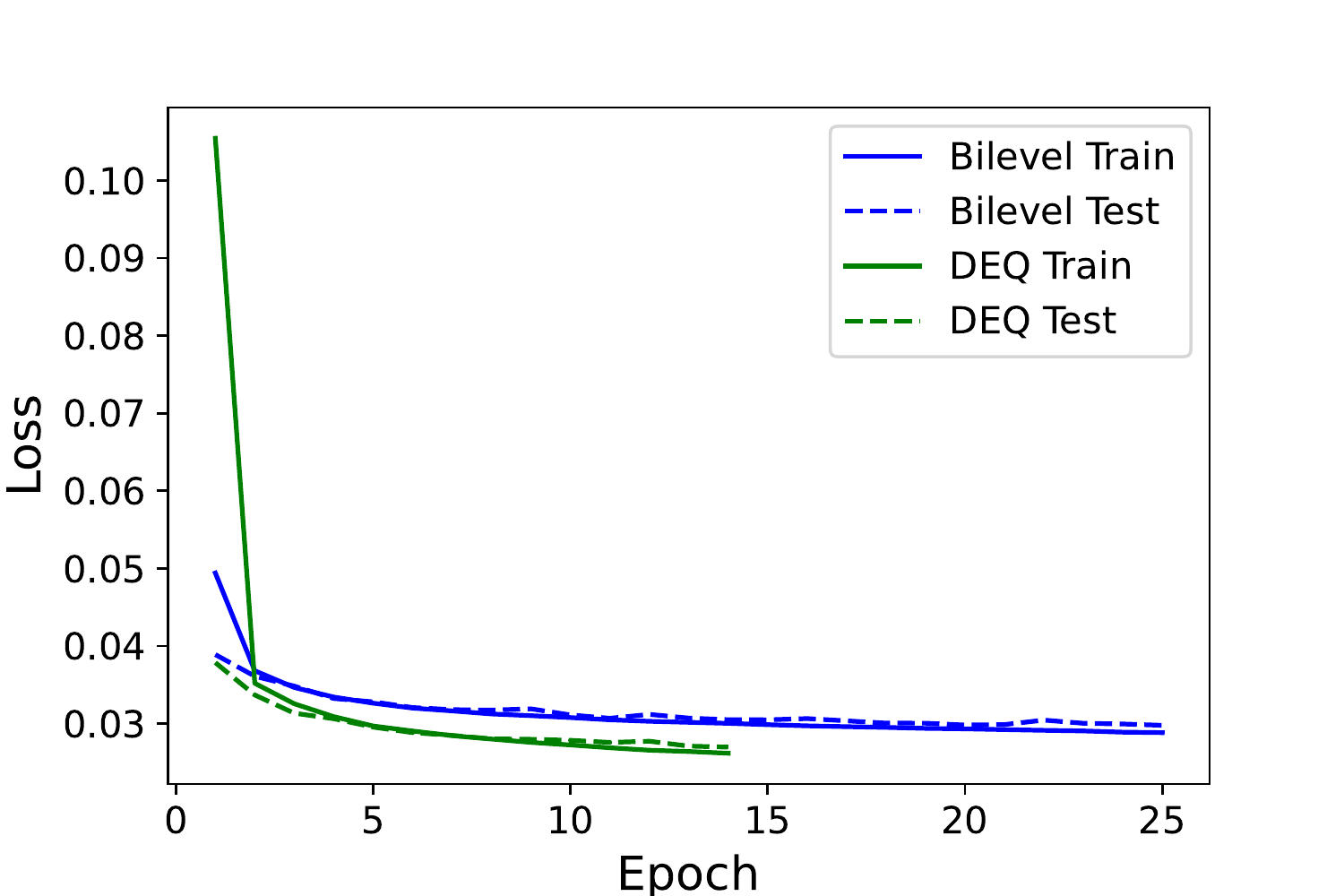}
  \includegraphics[width=0.32\linewidth]{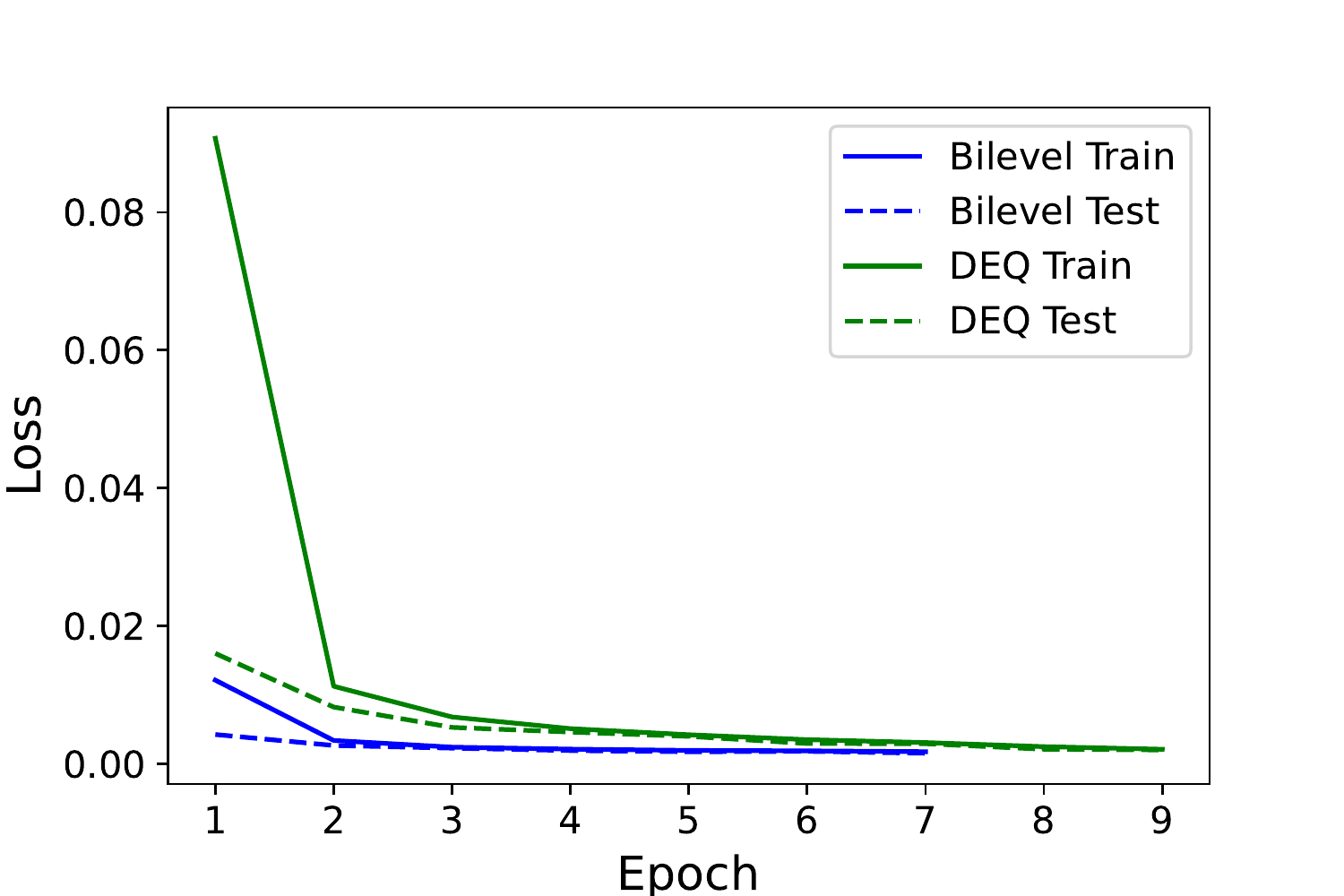}
 \end{minipage}
 \begin{minipage}{1.0\textwidth}
  \centering
  \includegraphics[width=0.32\linewidth]{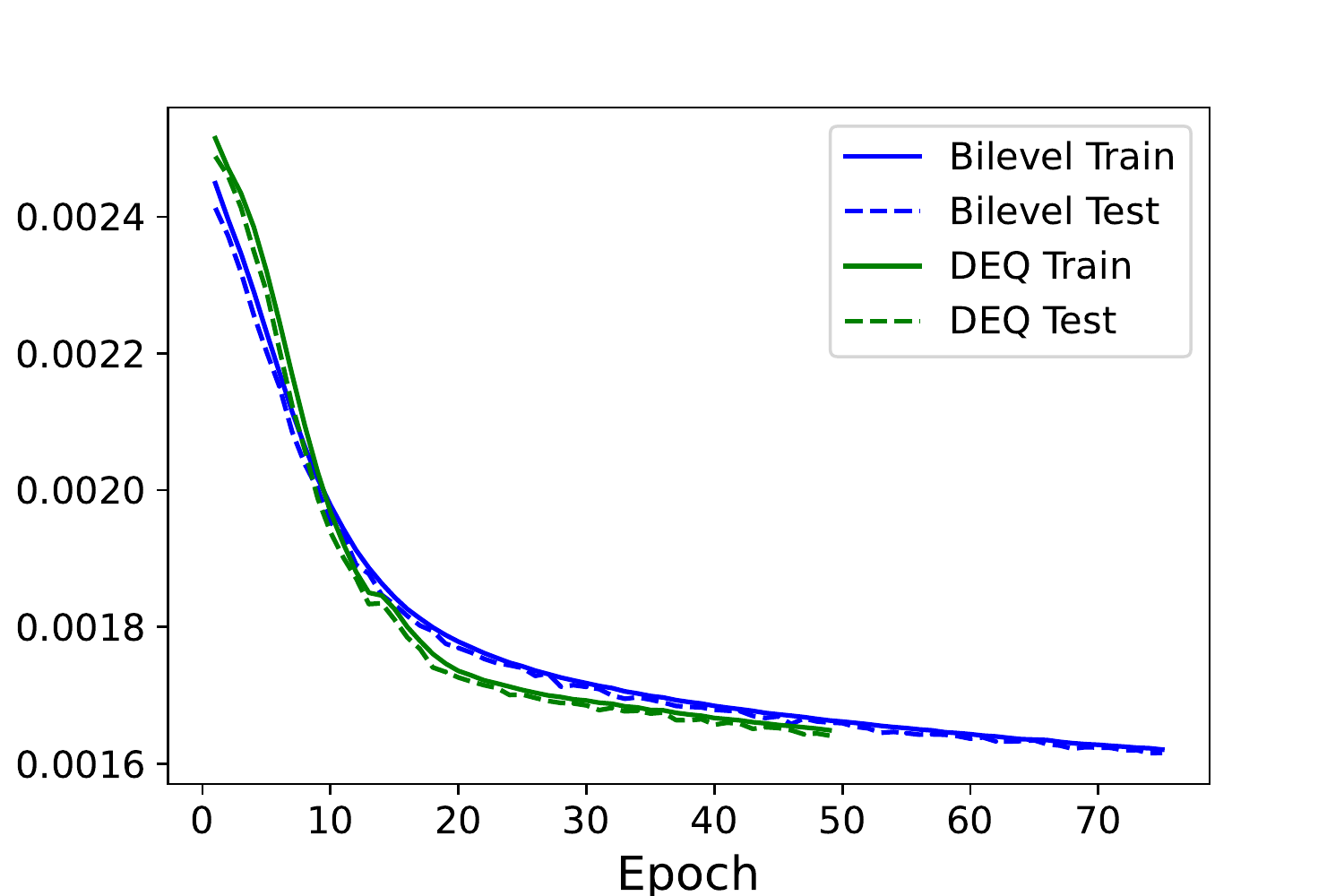}
  \includegraphics[width=0.32\linewidth]{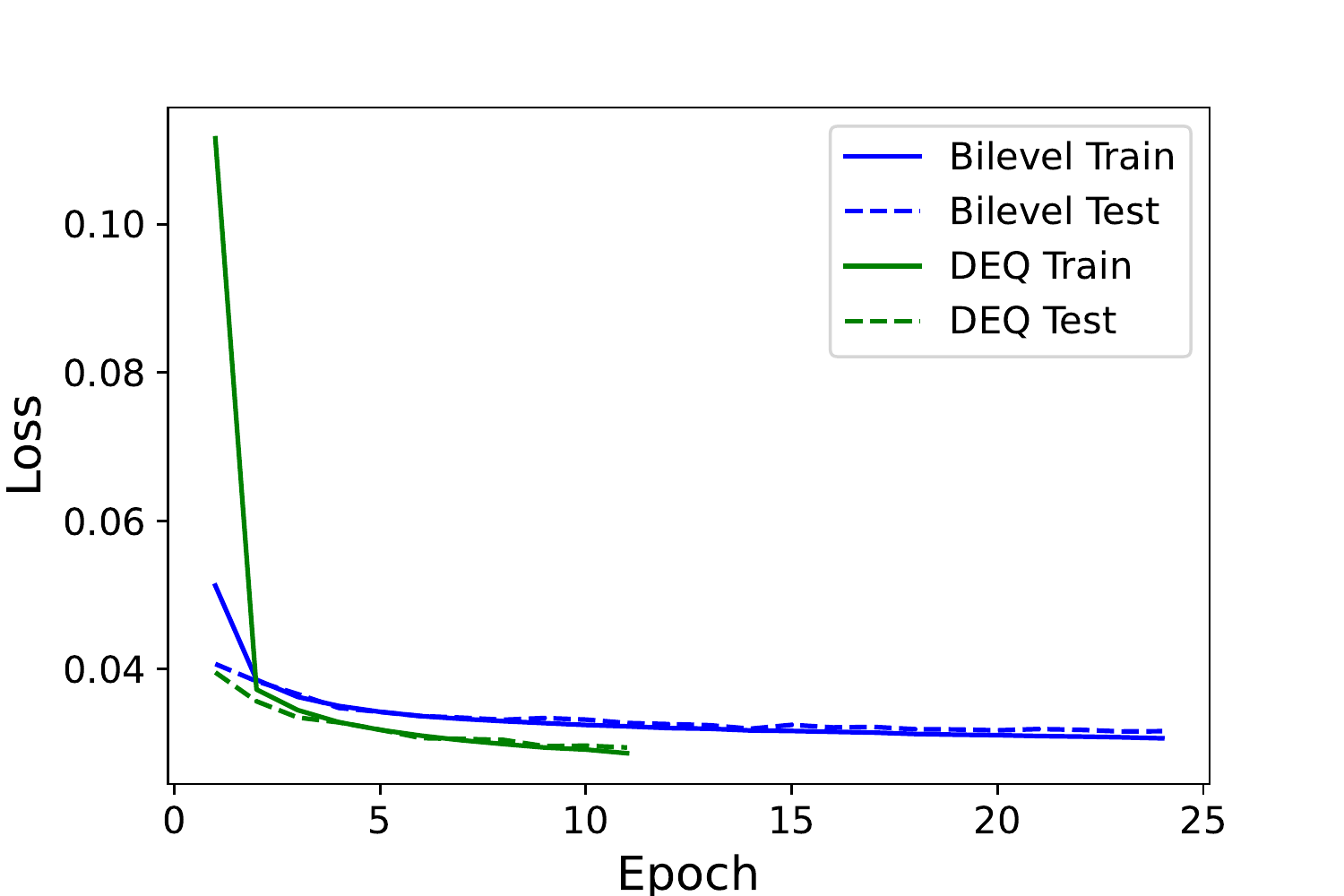}
  \includegraphics[width=0.32\linewidth]{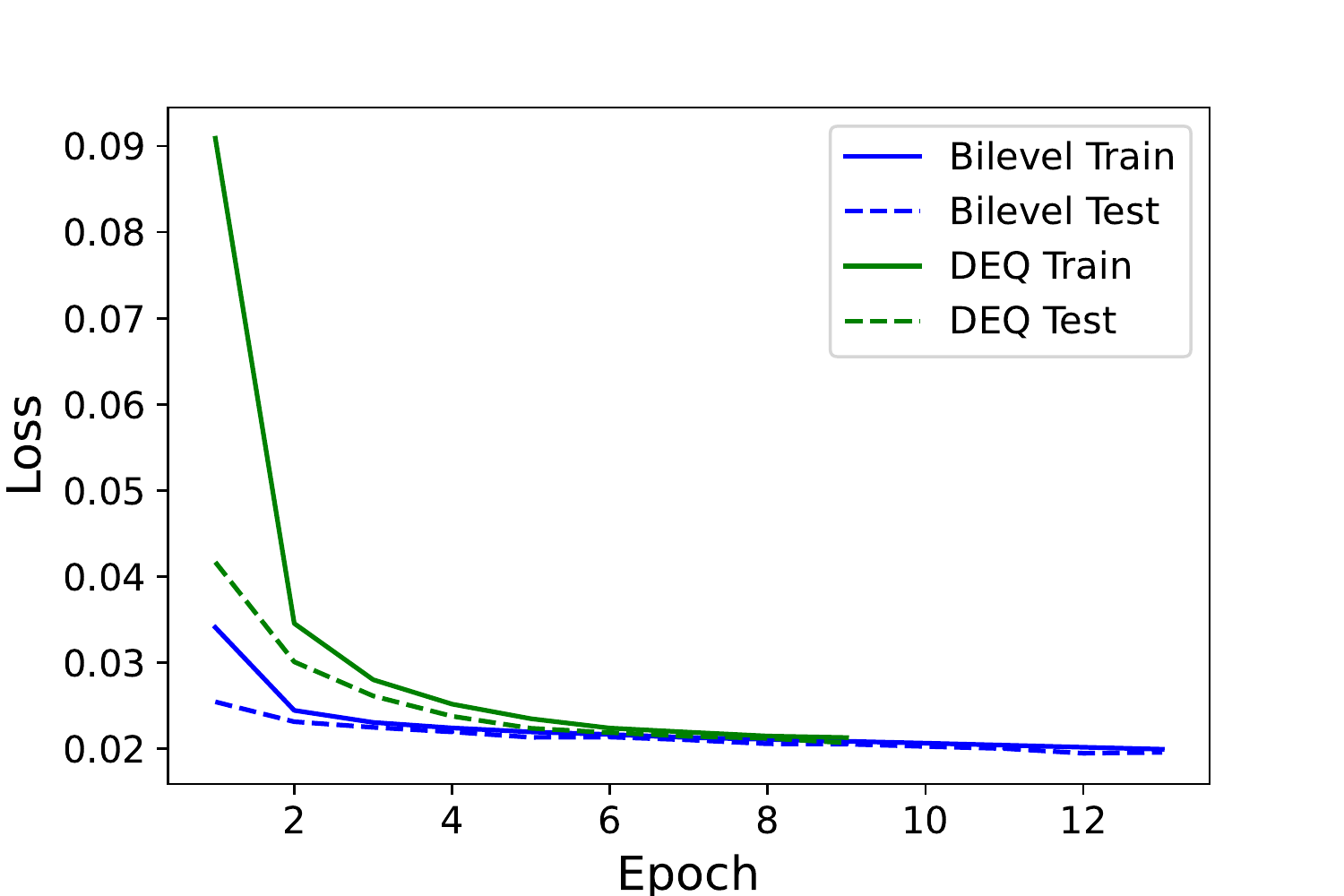}
 \end{minipage}
 \begin{minipage}{1.0\textwidth}
  \centering
  \includegraphics[width=0.32\linewidth]{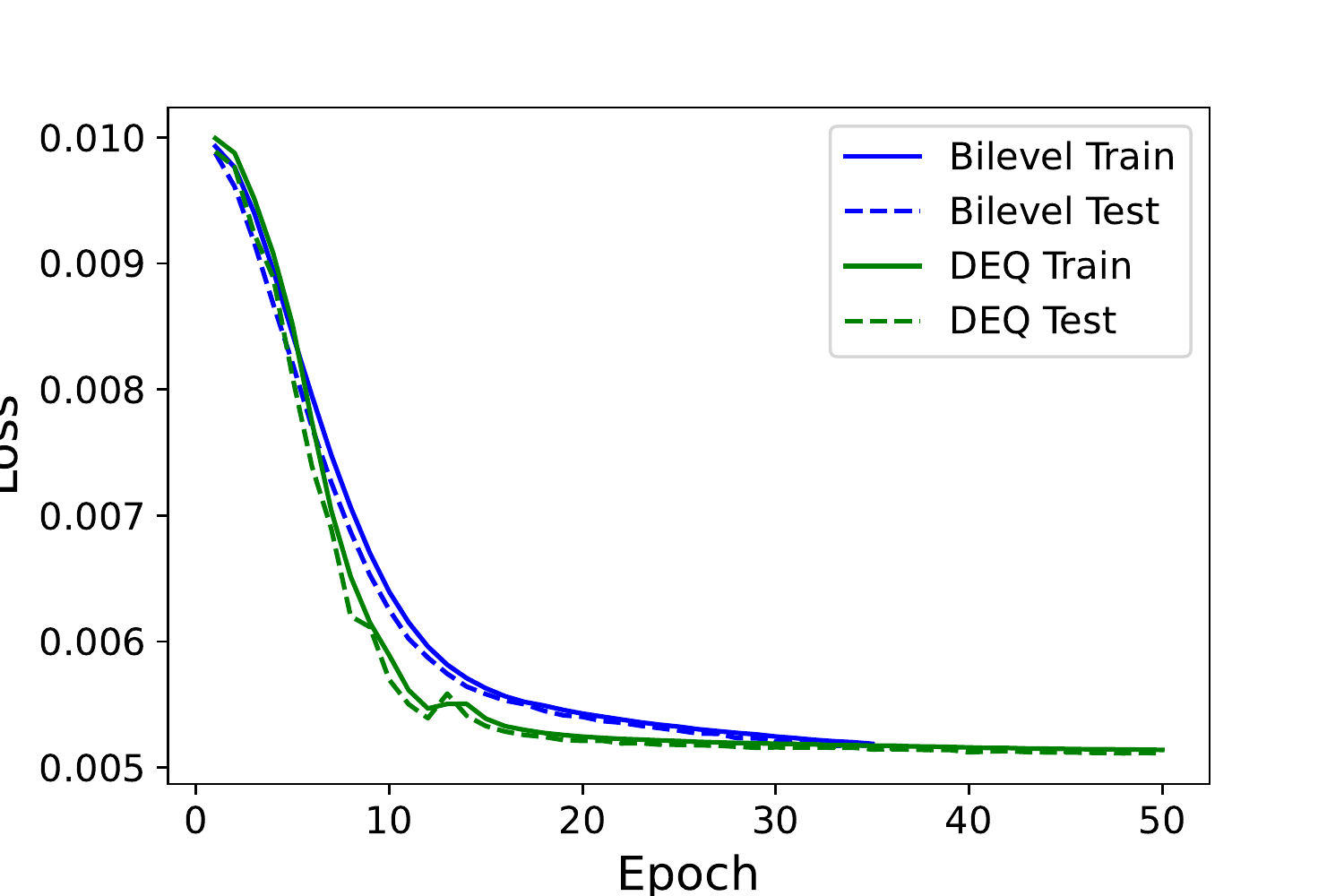}
  \includegraphics[width=0.32\linewidth]{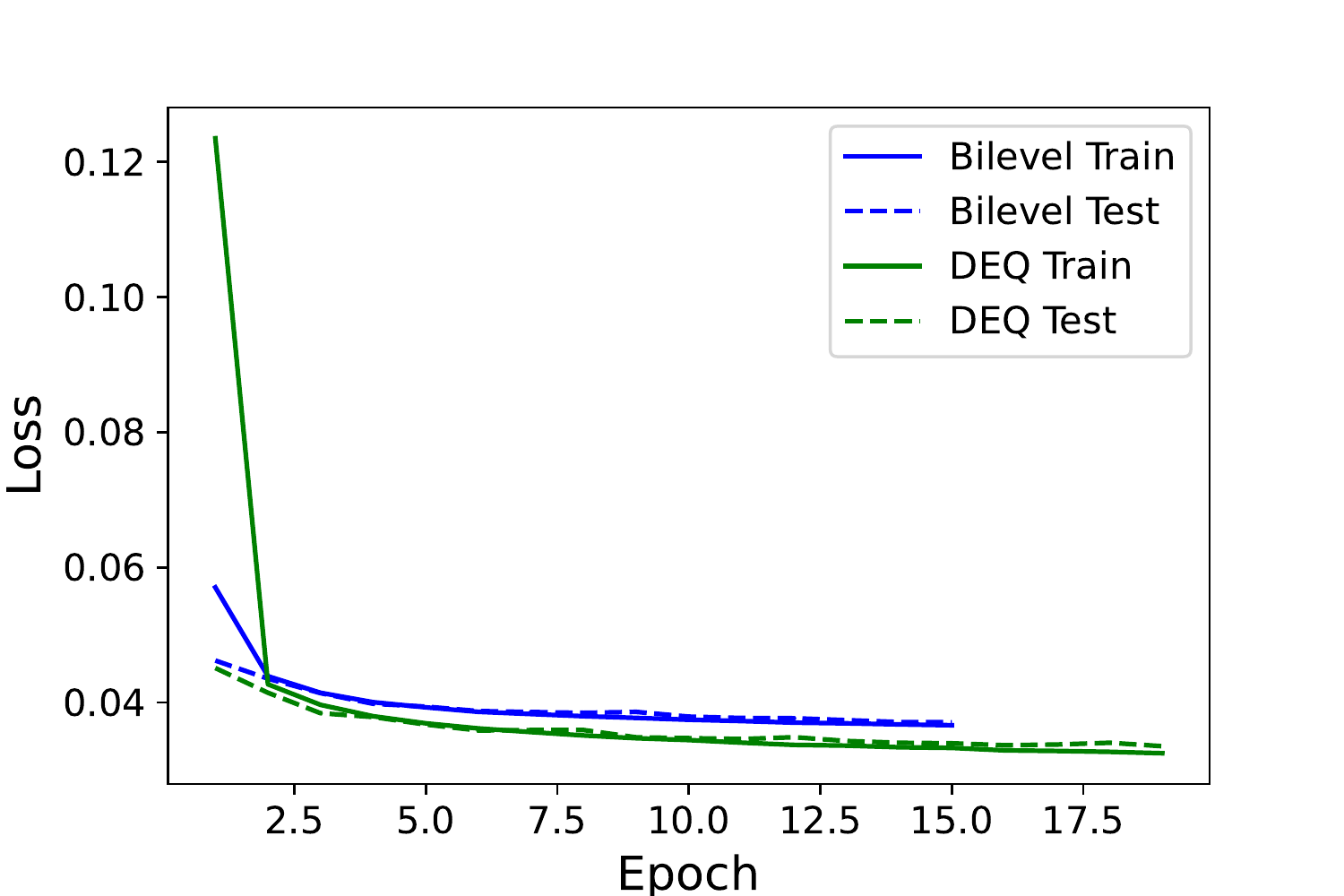}
  \includegraphics[width=0.32\linewidth]{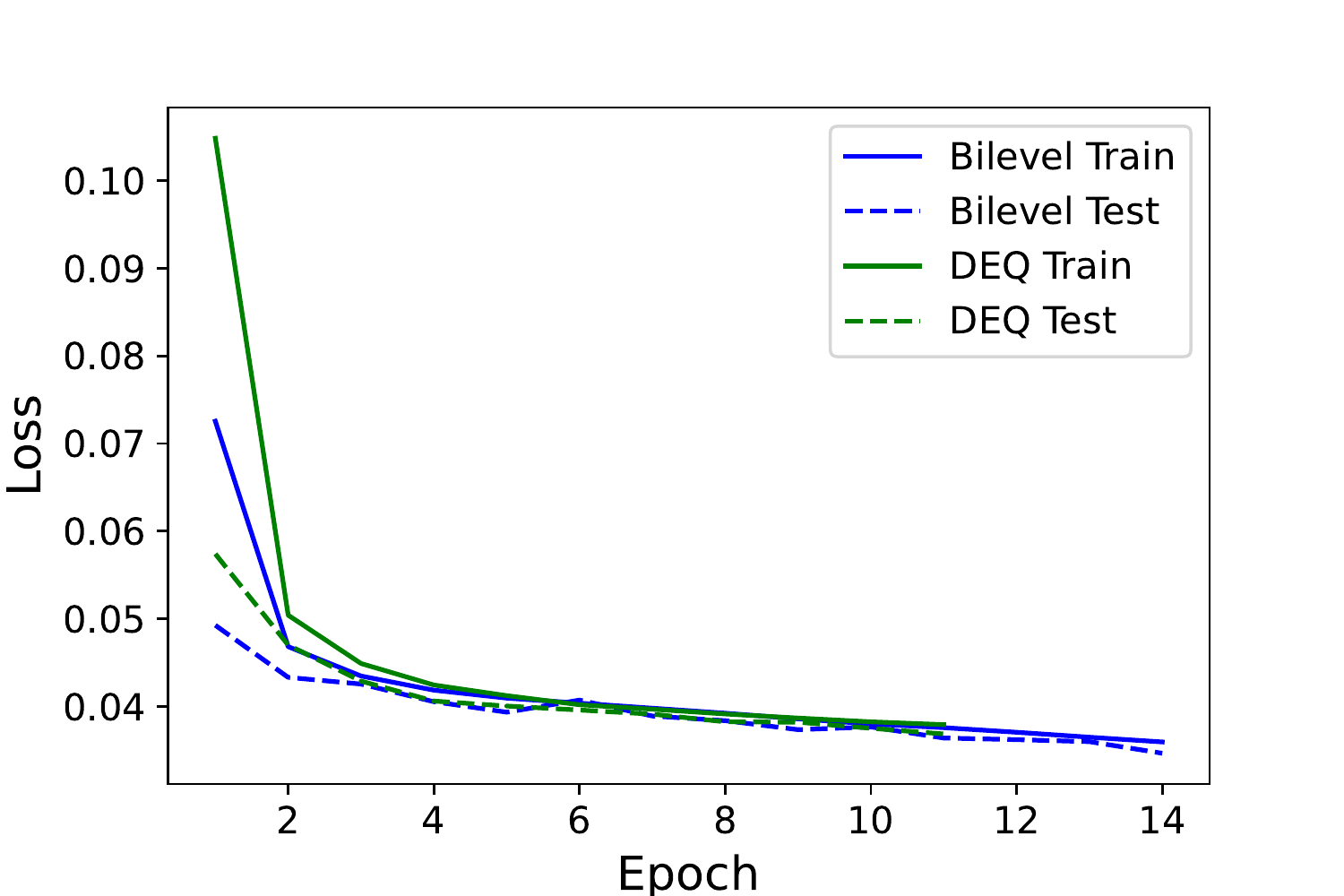}
 \end{minipage}
 \begin{minipage}{1.0\textwidth}
  \centering
  \includegraphics[width=0.32\linewidth]{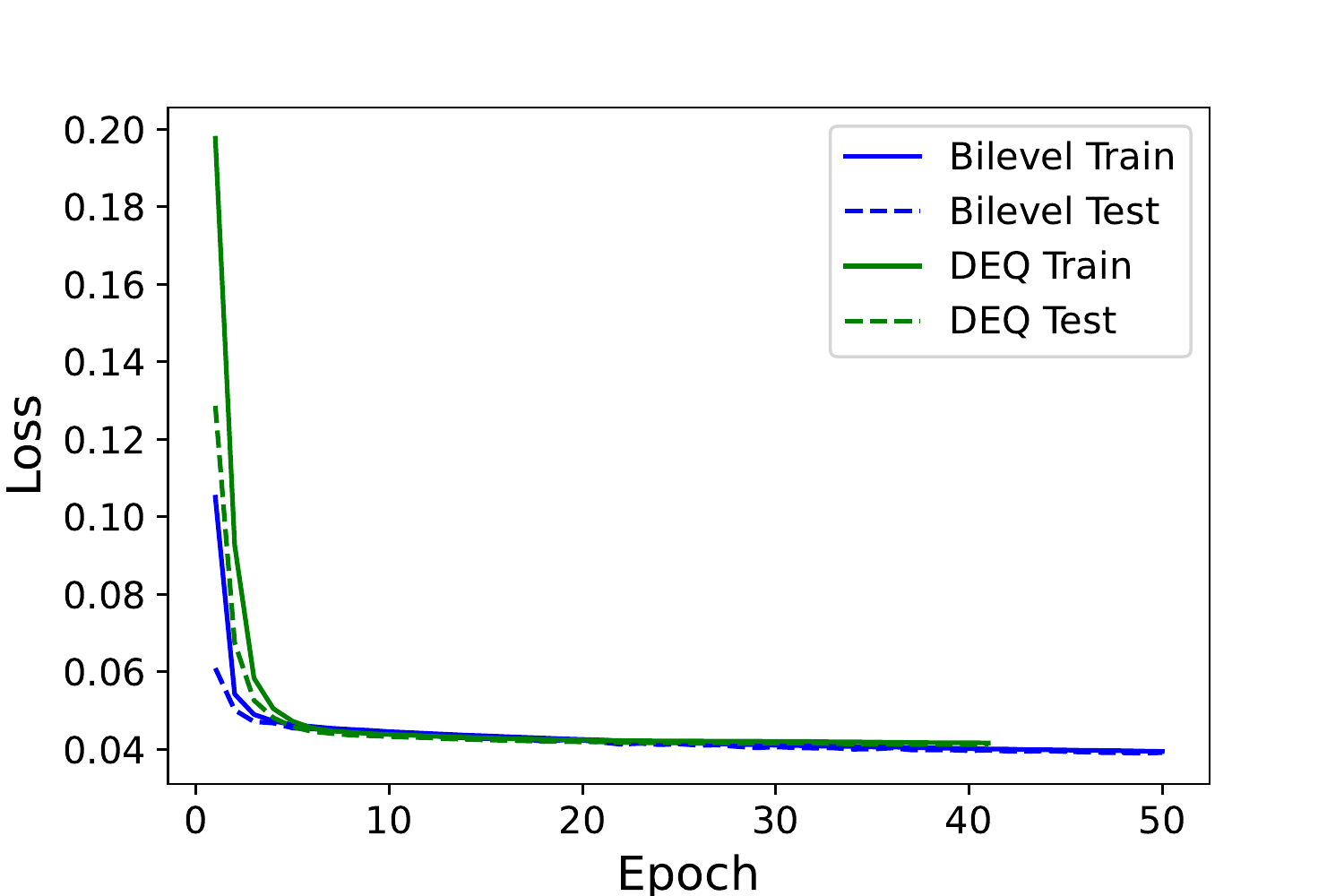}
  \includegraphics[width=0.32\linewidth]{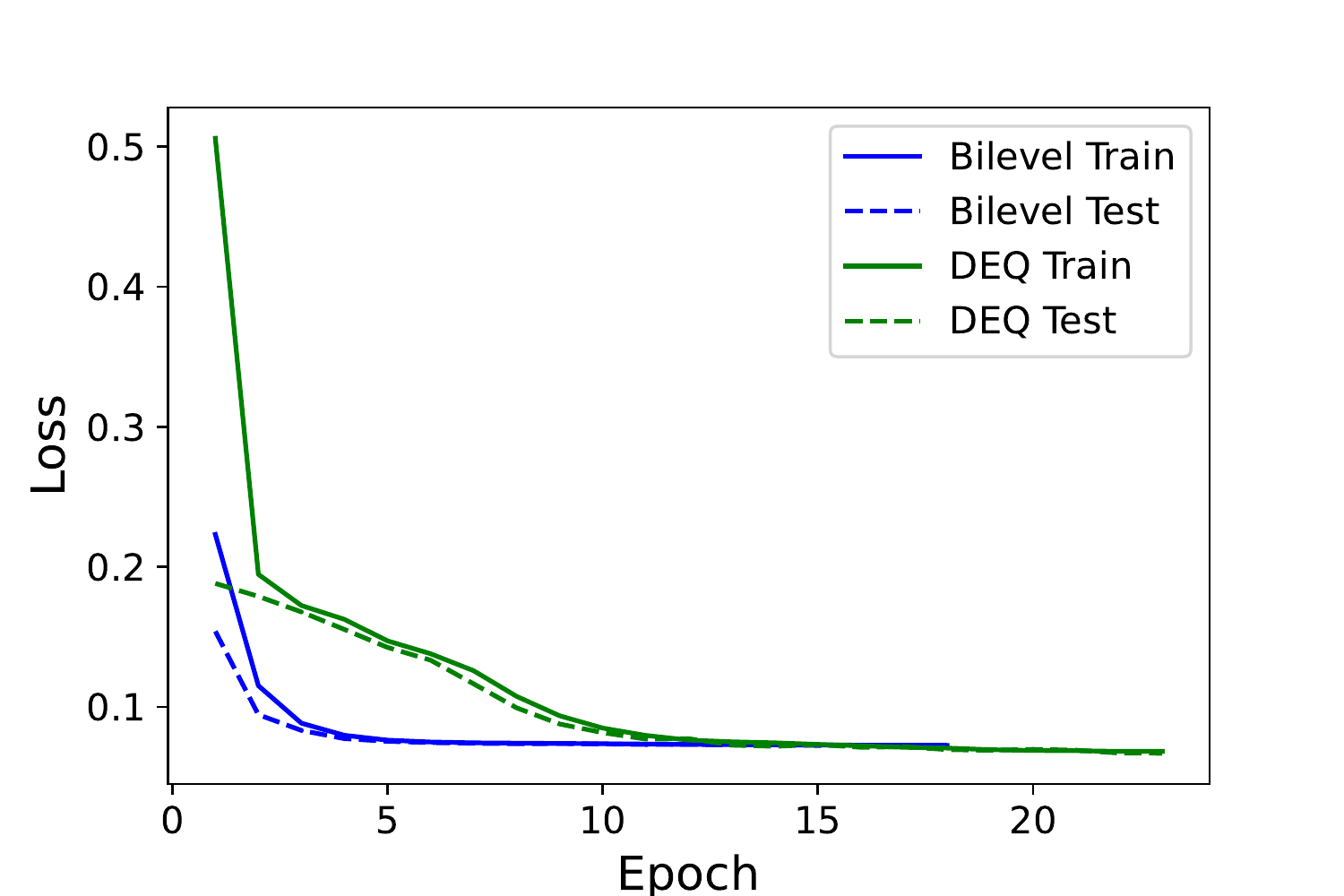}
  \includegraphics[width=0.32\linewidth]{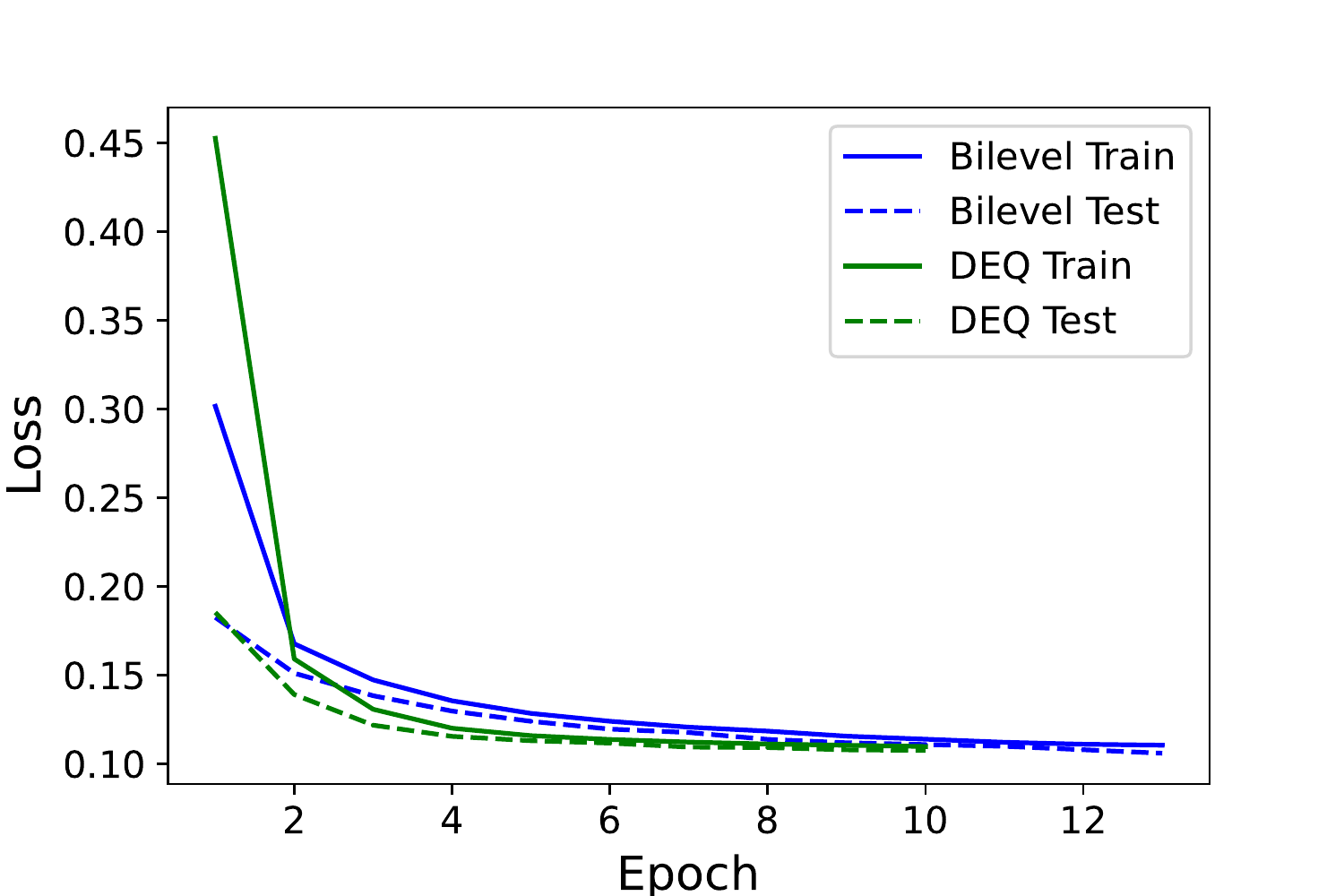}
 \end{minipage}
 \begin{minipage}{1.0\textwidth}
  \centering
  \includegraphics[width=0.32\linewidth]{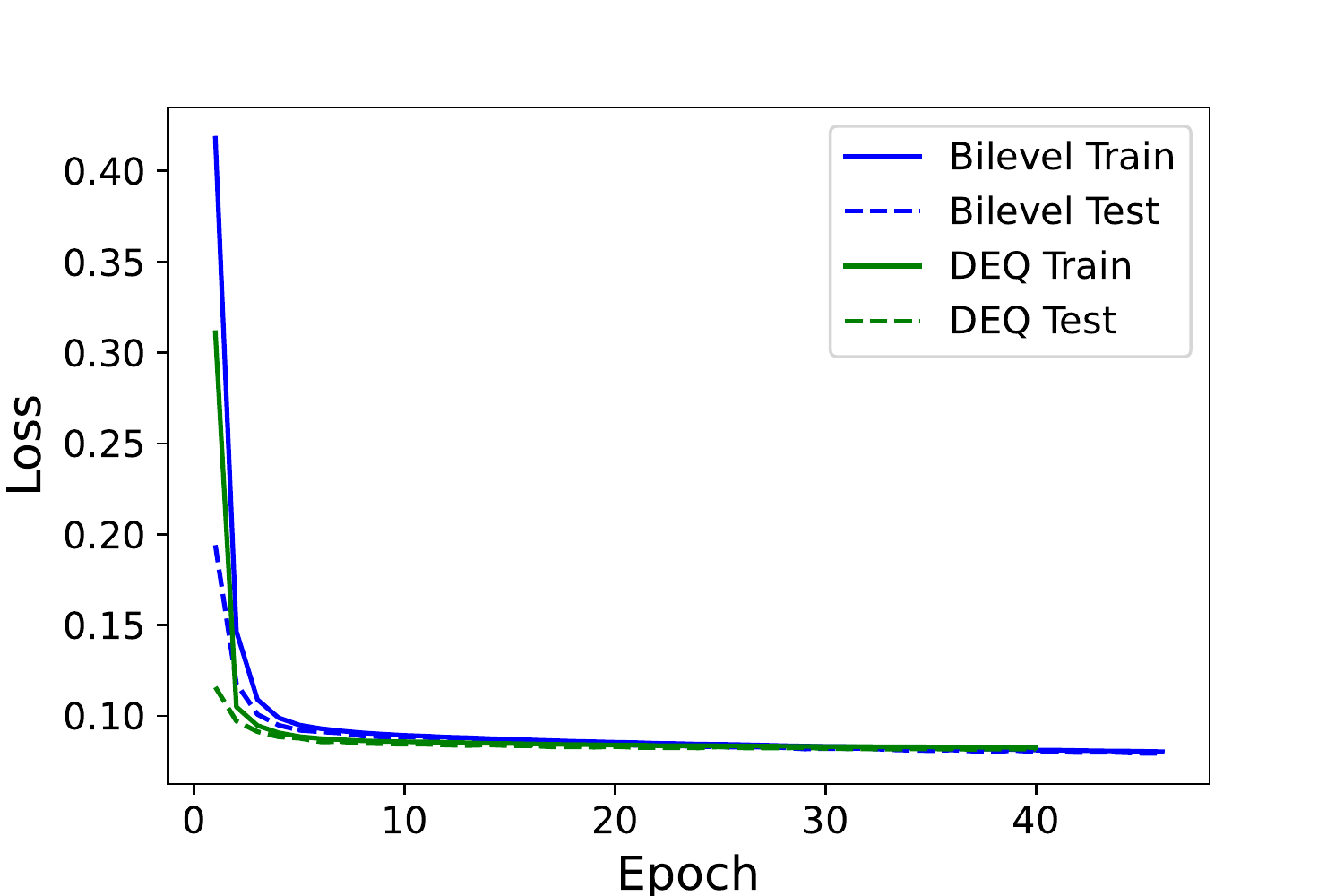}
  \includegraphics[width=0.32\linewidth]{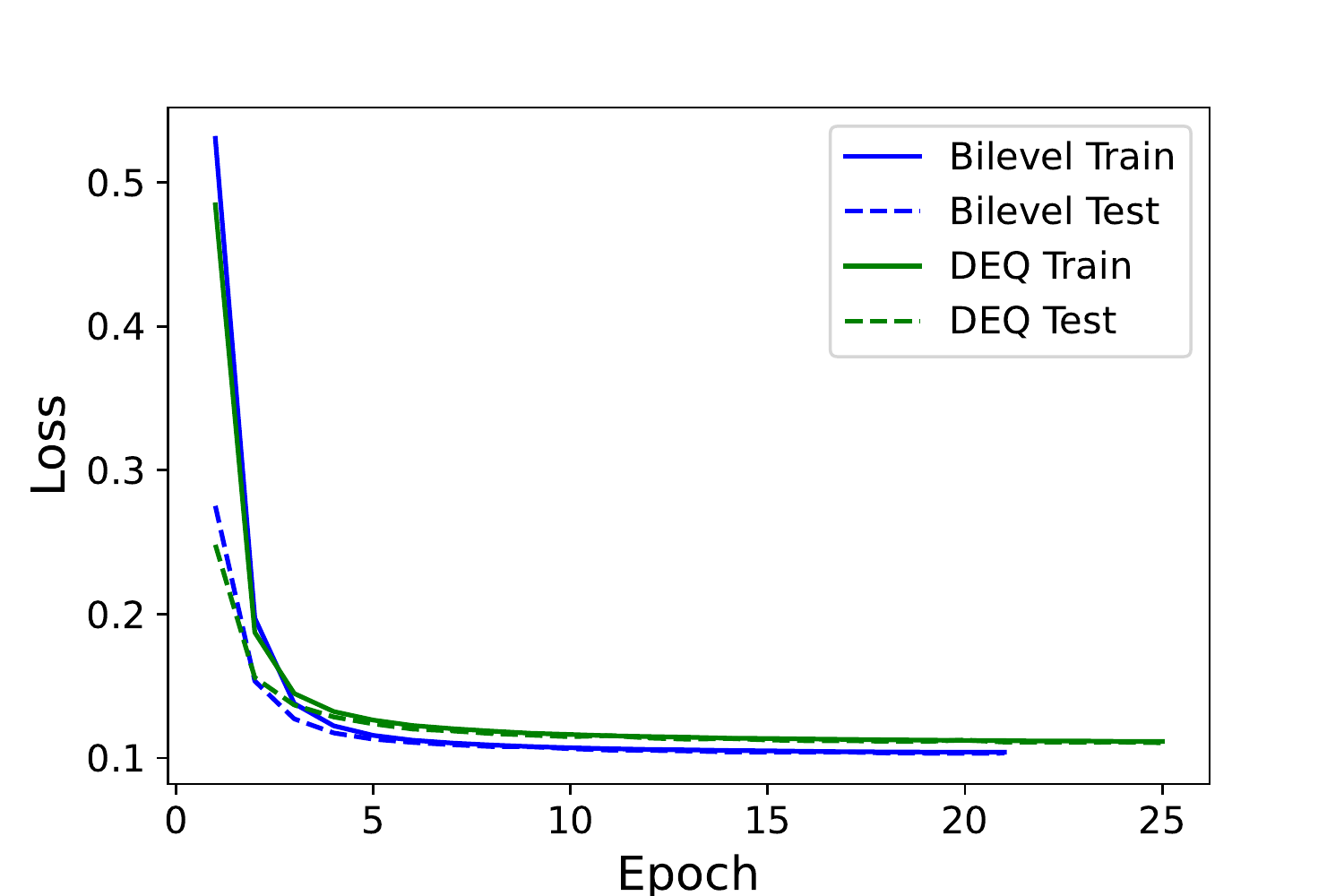}
  \includegraphics[width=0.32\linewidth]{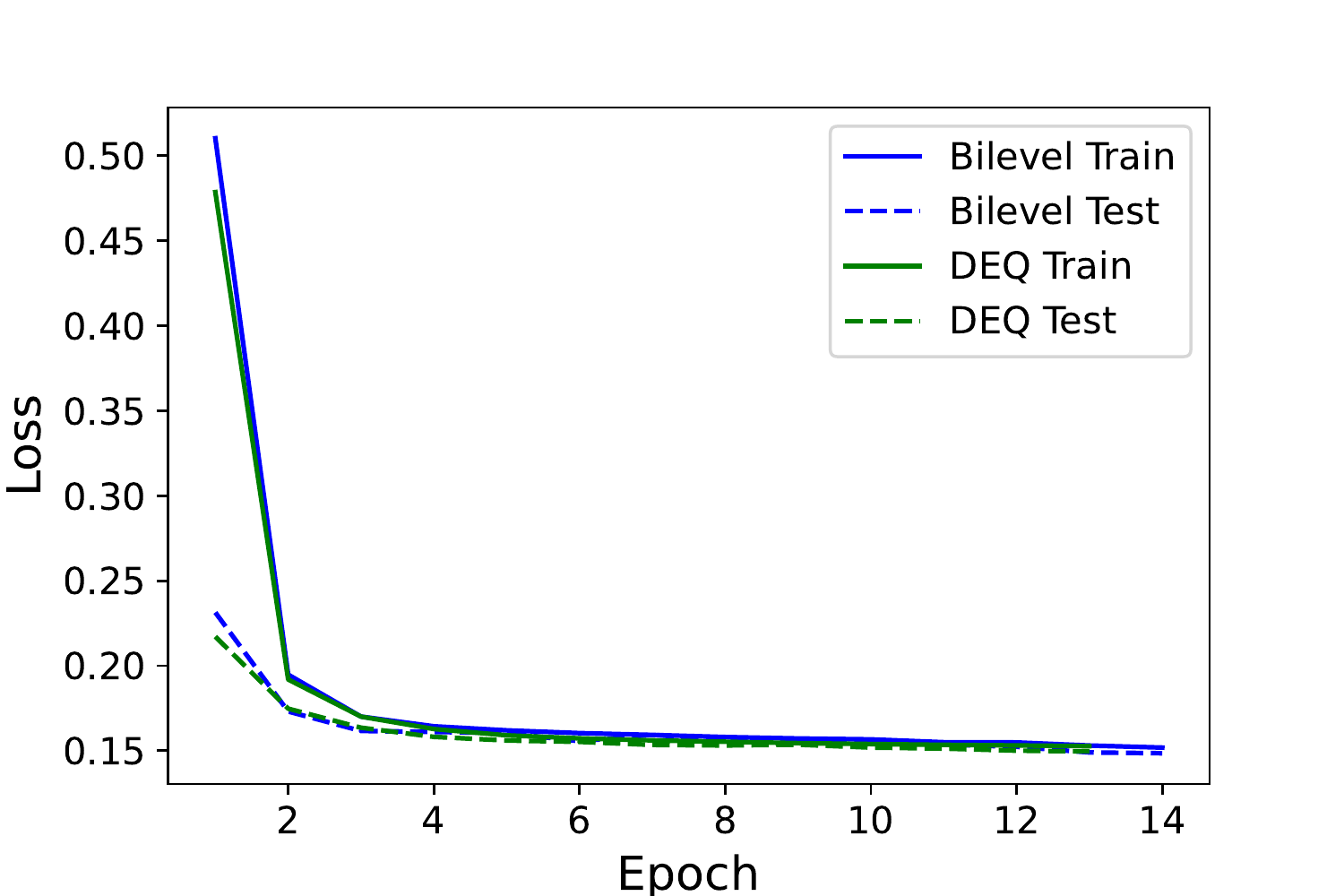}
 \end{minipage}
 \caption{Comparison of the loss error for the test dataset evaluated after each training epoch, for increasing values of noise levels in training (noise levels from top to bottom row: $0, 0.05, 0.1, 0.5, 1$). 
 Simulations are grouped by the tasks, namely denoising, inpainting, and deblurring (left, center, right columns).
 Each plot shows the simulation with the configurations that achieve the lowest final test loss.}
 \label{fig:loss-all-eps}
\end{figure}

The number of parameters in deep equilibrium models is larger than the one for bilevel methods in the framework we considered. 
Therefore, we want to assess whether this has an impact on the loss for different noise levels.
To do so, we visualize the loss for both training and test datasets. 
Results are shown in Fig.\ref{fig:loss-all-eps}.
Both methods seem to have the same behavior for increasing noise level values. Training takes longer for larger noise levels, and the loss becomes larger and larger as expected.

\subsection{Sensitivity against hyperparameter selection}
\begin{figure}[!h]
 \begin{minipage}{1.0\textwidth}
  \centering
  \includegraphics[width=0.5\linewidth]{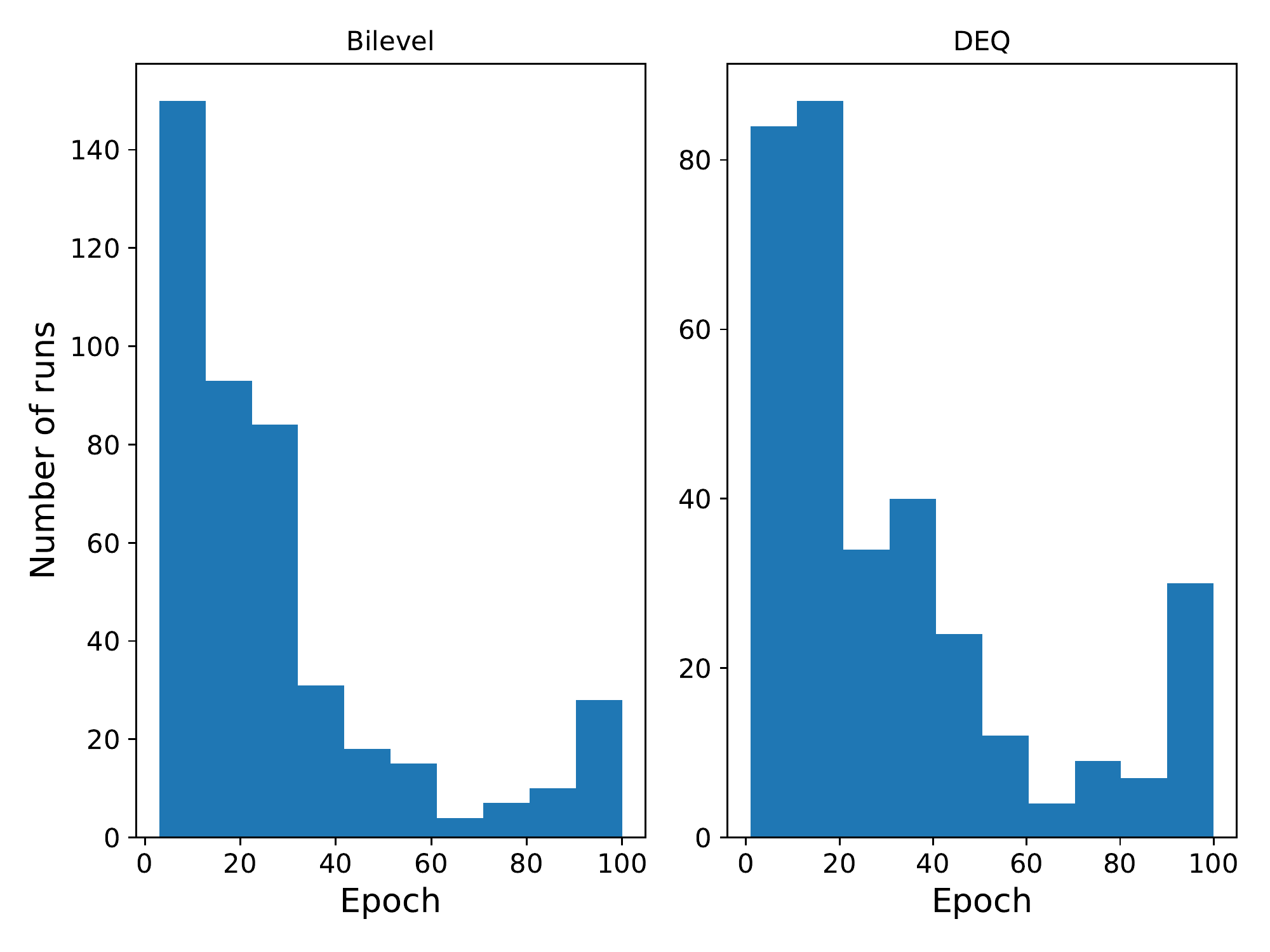}
 \end{minipage}
 \caption{These histograms show how many simulations were finished within an hour as a function of the number of epochs. 
 Each simulation is a different configuration of hyperparameters. 
 We consider only those runs where the loss on the test dataset is smaller than $0.5$.}
 \label{fig:hist-epochs--overall}
\end{figure}

Deep equilibrium models have twice as many parameters as bilevel learning methods in our setting. For our computations, we capped the running time to at most one hour, which is why we need to evaluate the impact of this choice on the results.
To do so, we consider all the simulations, keeping only those with a final test loss smaller than $0.5$, and collect the number of epochs in bins in a histogram. The resulting histograms are shown in Fig.\ref{fig:hist-epochs--overall}.
The shapes of bilevel and deep equilibirum histograms are similar, suggesting that capping the running time to one hour for the simulations does not have an impact on the quality of the results.
The different scale of the y axes for the two plots reflects that more simulations for the deep equilibrium models have exceded the loss threshold of $0.5$ in comparison to the bilevel methods. This further supports that the chosen bilevel methods are less sensitive to hyperparameter selection than their deep equilibrium counterparts. 
 
\subsection{Na\"{i}vely learning the fixed-point}
\label{subsec:naively-learning-fixed-point}
In this section we show the results we get if we choose to train the model by na\"{i}vely learning the fixed-points, as explained in Section~\ref{sec:naively-learning-fixed-point-does-not-work}. 
Here we train the parameters by minimizing \eqref{eq:denoising-aec-degrad}. To easily impose convergence to the fixed point, we perform this test with the bilevel optimization model.
In order to compare the reconstruction with the one we show in Fig.\ref{fig:inpainting__train_and_test_results_final}, we choose the same settings (namely, $\tau=0.5$, $\gamma=1.0$, and $\sigma=\text{Softshrink}$).

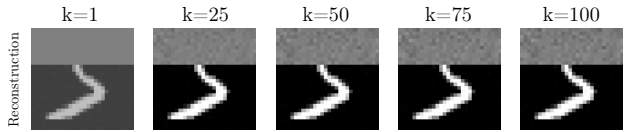
\begin{figure}[!h]
 \begin{minipage}{1.0\textwidth}
  \centering
  \begin{picture}(0,0)
   \put(15,41){{\scalebox{0.7}{k=1}}} %
   \put(60,41){{\scalebox{0.7}{k=25}}} %
   \put(105,41){{\scalebox{0.7}{k=50}}} %
   \put(151,41){{\scalebox{0.7}{k=75}}} %
   \put(195,41){{\scalebox{0.7}{k=100}}} %
   \put(-5,1){\rotatebox{90}{\scalebox{0.5}{Reconstruction}}} %
  \end{picture}
  \includegraphics[angle=0, width=0.5\linewidth]{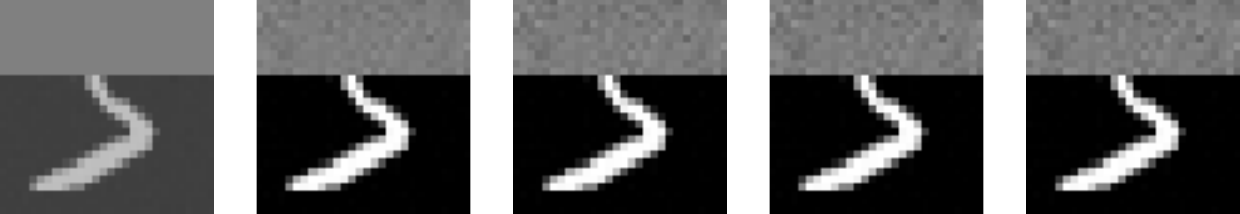}
 \end{minipage}
 \caption{Reconstruction for the inpainting task for a bilevel optimization model whose parameters have been trained by minimizing the error of the reconstruction $u^*$ w.r.t. the true image $u^\dagger$ (na\"{i}ve approach). We show how the reconstruction $\{u^k\}$ changes for different values of the iteration $k$. As we can see, the model trained with the na\"{i}ve approach is not able to inpaint the masked area.}
 \label{fig:naive_training_reconstruction_2}
\end{figure}

After training the parameters using \eqref{eq:denoising-aec-degrad}, we check whether the trained model is able to fill the missing pixels from a masked image.
To do so, we initialize $u_0$ in the same way explained at the beginning of Section~\ref{sec:implementation} and we run 
100 iterations of \eqref{eq:lower-level-gradient-descent-fixed-point}. 
The results are shown in Fig.\ref{fig:naive_training_reconstruction_2}, where it is clear that missing regions are not filled. We again emphasize that the model and the hyperparameters are the same ones used in Fig.\ref{fig:inpainting__train_and_test_results_final}; the reconstruction is satisfactory in the latter. The only difference between the two models is the way in which the parameters have been learned. This empirically supports the claims made in Section \ref{sec:naively-learning-fixed-point-does-not-work}.

\subsection{Comparison on a higher-resolution dataset}
\label{subsec:comparison-celeba}
We now compare how bilevel learning and DEQ model perform on CelebA, a higher-resolution dataset.
For simplicity, the images are converted to grayscale, and pixel values are normalized in the range $[-1, 1]$.
Because of the dimension of the images, using fully-connected dense matrices for $A$ and $C^\top$ is not a viable option, as it would lead to memory issues. Instead, we use 2D convolutional layers, allowing us to have fewer parameters.
We choose to omit the inpainting task for this dataset, as it necessitates non-local information that convolutional layers cannot provide, in contrast to fully-connected dense layers.

\subsubsection{Denoising}
Motivated by the Rudin Osher Fatemi (ROF) model \cite{rudin1992nonlinear}, we modify 
\eqref{eq:deq-gd} and \eqref{eq:deq-non-gradient} to
$$u^{k + 1} = u^k - \tau \left(\lambda \, K^\top \left( Ku^k - f^\delta \right) + \gamma C^{\top} \sigma (\xi A u^k) \right) \, $$
with $\sigma(\cdot) = \tanh(\cdot)$, and $\xi \in \R$ `large enough' to approximate 
the \textit{sign} activation function.
Choosing $\xi=100$, $\gamma=1$, the value of $\lambda$ that empirically minimizes the MSE for noise with standard deviation $0.1$ is $\lambda=18.156$. Please note that for the anisotropic ROF model, $\mathcal{R}(u)$ reads $\mathcal{R}(u) = \| \nabla u \|_1$, where $\| \cdot \|_1$ denotes the one-norm and $\nabla$ a discretization of the gradient operator. Hence $\nabla^\top \textit{sign}(\nabla u)$ is a subgradient of this function, and replacing $\nabla$ and $\nabla^\top$ with operators $A$ and $C$ and \textit{sign} with the hyperbolic tangent yields a trainable model in the vein of ROF denoising.

We consider different initializations for the denoising task, with either 2 or 30 output channels and kernels of size $11\times 11$ or $3\times 3$ kernel, respectively. This is done to test two extreme scenarios: either having a low number of output channels with relatively large kernels, compared against many output channels with small kernels. The number of input channels is dataset-dependent, and for grayscale images, it is one. 
The kernel weights of $A$ and $C$ are initialized with the same values for both Deep Equilibrium models and bilevel learning, whereas in the latter we add the constraint that 
$C=A$ through all the training.

We use Adam as the optimizer, together with a scheduler which makes the learning rate decrease linearly over the epochs, starting in the range $\left[3.2\cdot 10^{-3},10^{-2}\right]$ and ending with a value which is $100$ times smaller. This is done to help escaping local minima in the first epochs, and to converge faster towards the end of the training. 
All simulations are run over $500$ epochs.
For the computation of the fixed point we use a na\"{i}ve forward iteration scheme that we stop when either $1,000$ iterations have been computed, or a relative norm difference between iterations below $10^{-14}$ has been achieved. 
We decide to perform spectral normalization for $A$ and $C$, i.e., we divide the operators by their norms which are estimated using a power iteration scheme. 
An alternative to this approach is to avoid computing the spectral normalization. In the latter case, \eqref{eq:deq-gd} may not converge unless an appropriate value of $\tau$ is chosen after every optimization step. To guarantee the convergence to the fixed point, we check the gradient values before performing the optimization step, and we decrease $\tau$ by an arbitrary factor of $10$ if at least one component of the gradient becomes either $\infty$, or NaN; at the same time, we also increase the maximum number of iterations by the same factor to compute the fixed point \eqref{eq:deq-gd}.

\begin{figure}[p]
 \centering \includegraphics[trim={0, 0, 8.3cm, 0}, clip, width=0.66\linewidth]{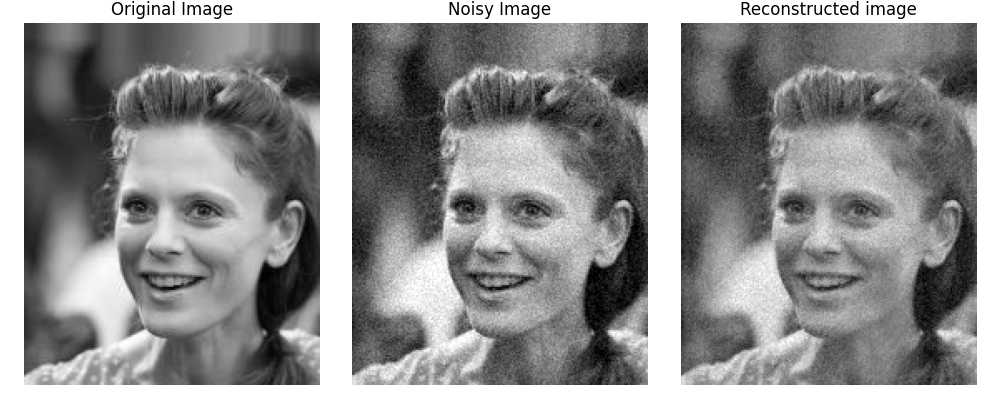}
 \\
 \subfloat[Reconstructed images, two output channels, kernels size $11\times 11$]{
 \begin{minipage}{0.33\textwidth}
  \centering
  \includegraphics[trim={17cm 0 0cm 0.505cm},clip, width=1.0\linewidth]{review-celeba/cuda_run_007_reconstruction__before_training_test.png}
 \end{minipage} \hfill
 \begin{minipage}{0.33\textwidth}
  \centering
  \includegraphics[trim={17cm 0 0cm 0.505cm},clip, width=1.0\linewidth]{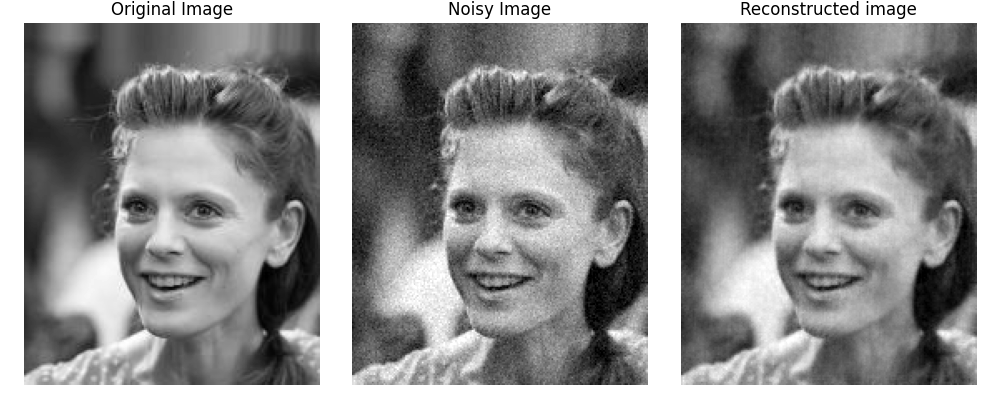}
 \end{minipage}
 \hfill
 \begin{minipage}{0.33\textwidth}
  \centering
  \includegraphics[trim={17cm 0 0cm 0.505cm}, clip, width=1.0\linewidth]{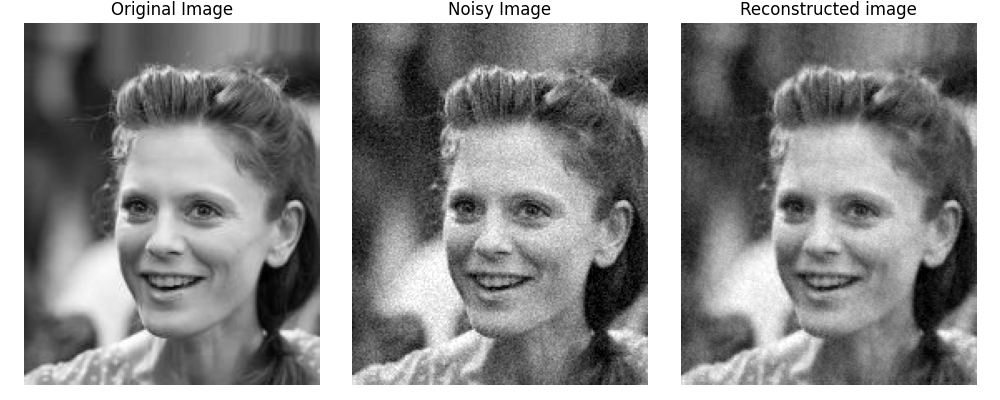}
 \end{minipage} \hfill
 }\\
 \subfloat[Reconstructed images, thirty output channels, kernels size $3\times 3$]{
 \begin{minipage}{0.33\textwidth}
  \centering
  \includegraphics[trim={17cm 0 0cm 0.505cm},clip, width=1.0\linewidth]{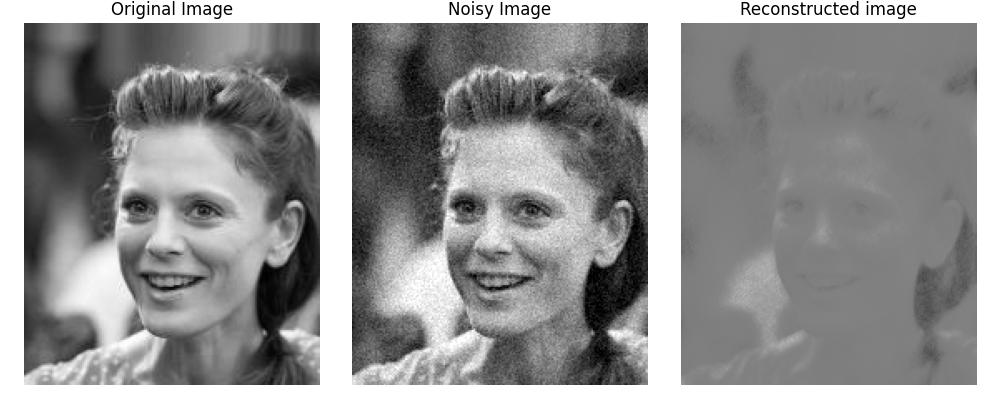}
 \end{minipage} \hfill
 \begin{minipage}{0.33\textwidth}
  \centering
  \includegraphics[trim={17cm 0 0cm 0.505cm},clip, width=1.0\linewidth]{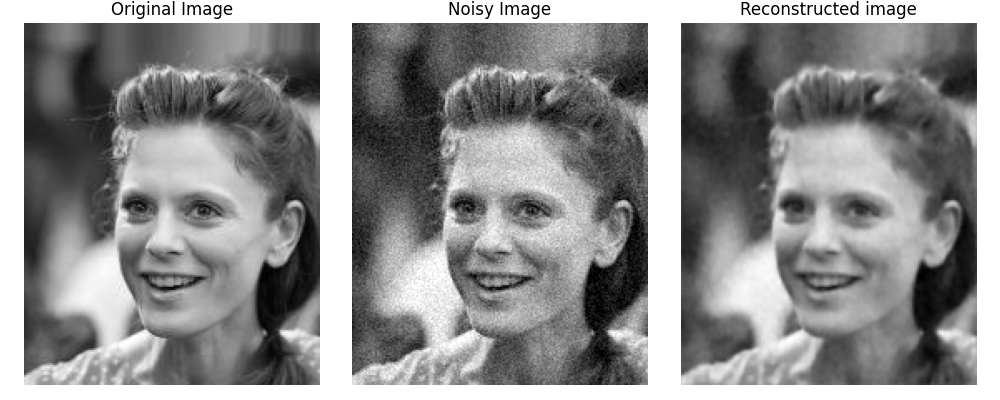}
 \end{minipage}
 \hfill
 \begin{minipage}{0.33\textwidth}
  \centering
  \includegraphics[trim={17cm 0 0cm 0.505cm}, clip, width=1.0\linewidth]{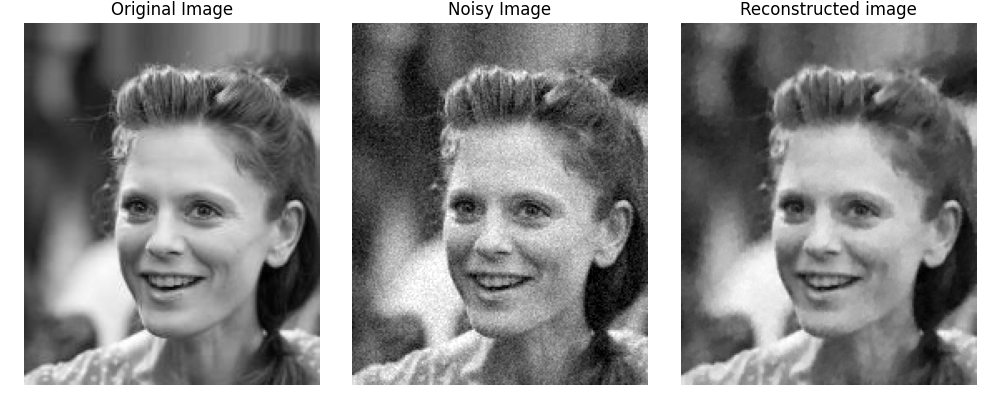}
 \end{minipage} \hfill
 }
 \caption{Denoising CelebA; sample from the test dataset. The first row contains the original image $u$ and the noisy image $f^\delta$. From left to right in the second and third rows: reconstructed image with random initializations of the kernels (left), with parameters learned using bilevel learning (center), and parameters learned using the DEQ model (right).}
 \label{fig:denoising__celeba_spectral_norm}
\end{figure}

\begin{figure}[htbp!]
 \subfloat[$11\times 11$ kernels before training.]{
  \centering
  \includegraphics[trim={0, 10.3cm, 0, 1cm}, clip, width=1.0\linewidth]{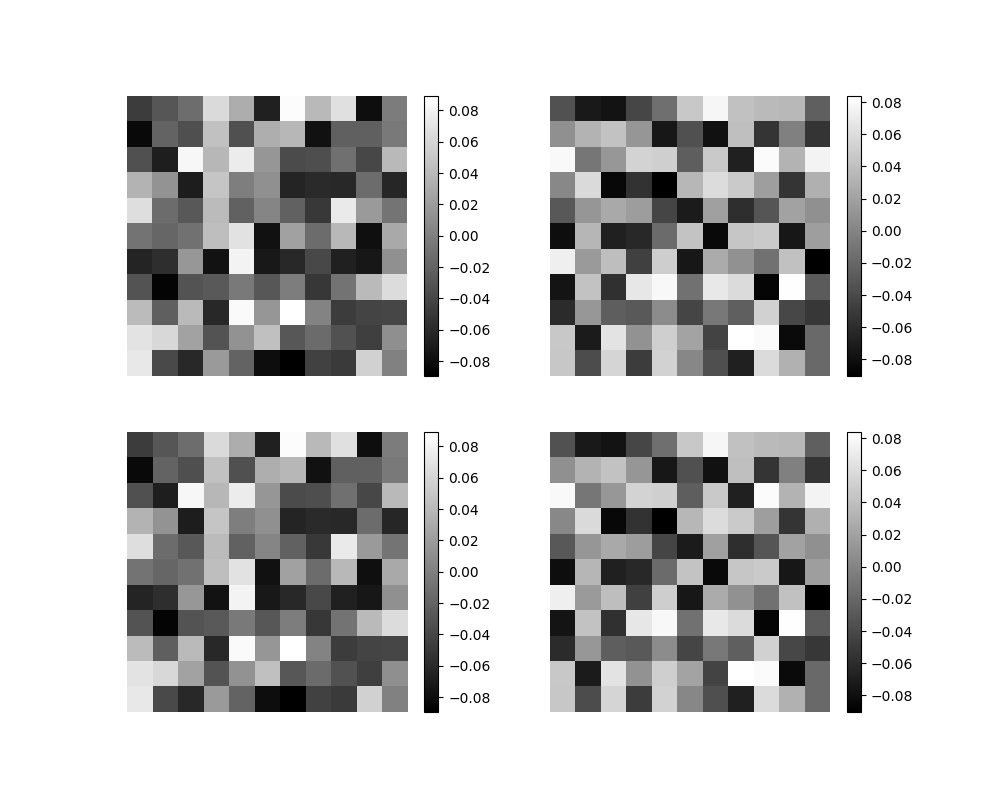}
  }\\
  \hfill
 \subfloat[$11\times 11$ kernels after training.]{
  \centering
  \includegraphics[trim={0, 10.2cm, 0, 2.1cm}, clip, width=1.0\linewidth]{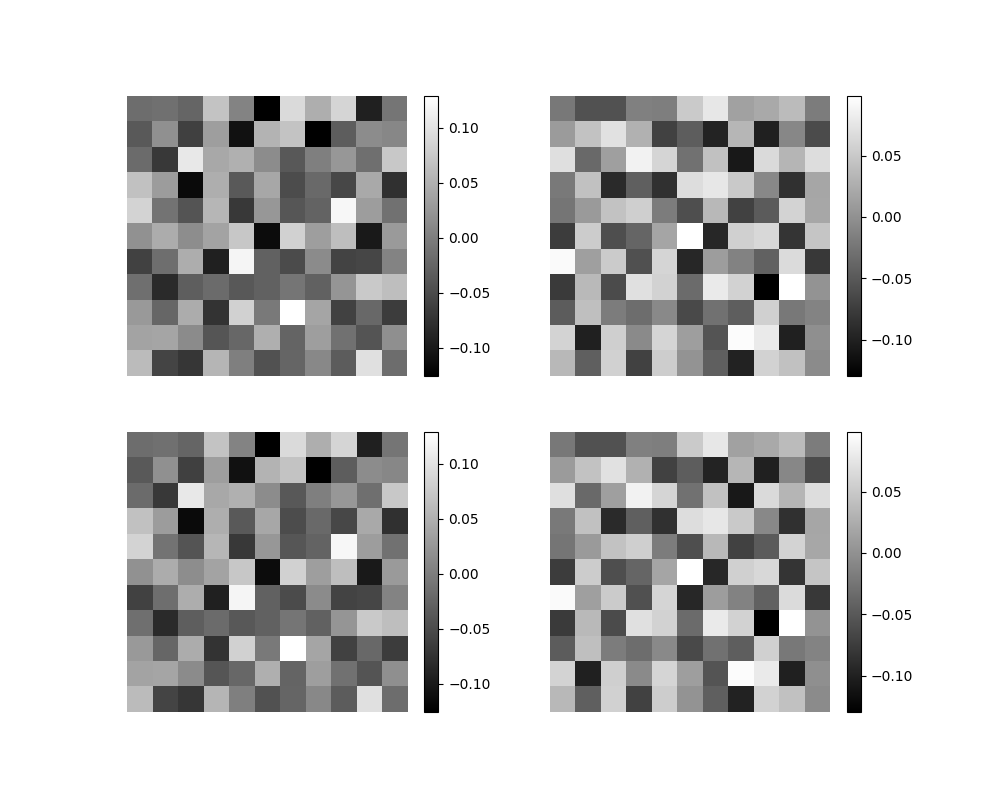}
  }
 \caption{Comparison of $11\times 11$ kernels of $A$, shown before and after training on the denoising task using the bilevel learning model. The weights of $C^\top$ are not shown here since $C^\top=A^\top$ for bilevel learning models.}
 \label{fig:comparison-kernels-11-11-denoising-BL}
\end{figure}

\begin{figure}[htbp!]
 \subfloat[$11\times 11$ kernels before training.]{
  \centering
  \includegraphics[trim={0.5cm, 10.2cm, 2.0cm, 2.2cm}, clip, width=0.5\linewidth]{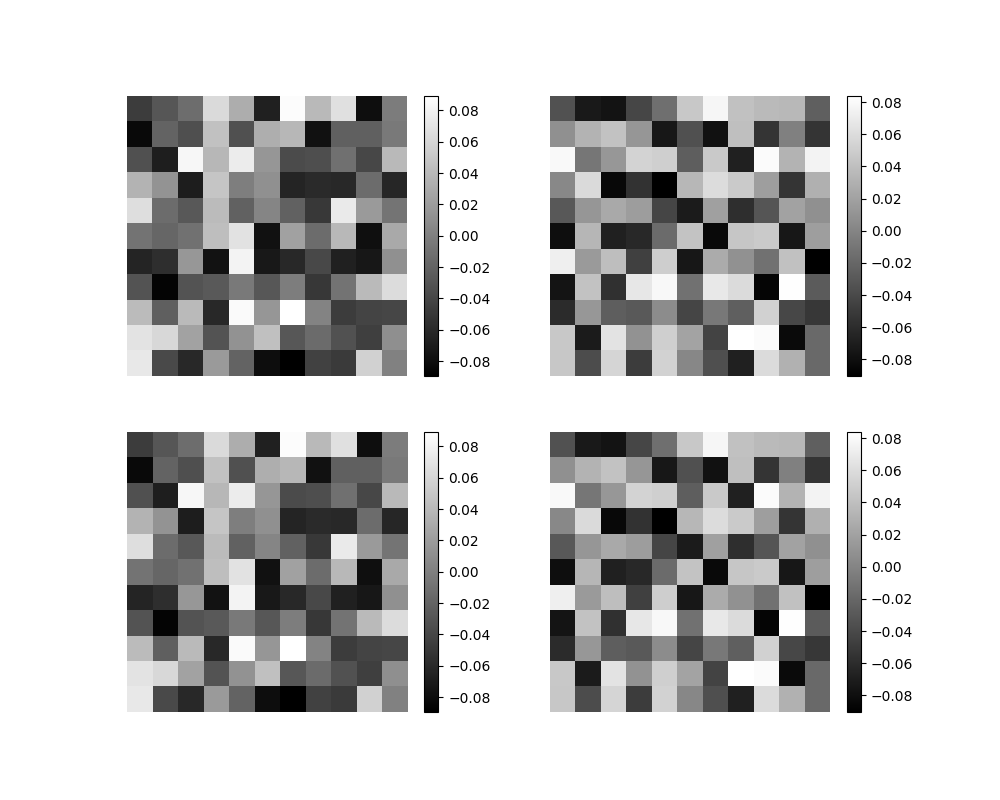}
  \includegraphics[trim={0.5cm, 1.67cm, 2.0cm, 10.7cm}, clip, width=0.5\linewidth]{review-celeba/cuda_run_007_kernels__before_training.png}
  }\\
  \hfill
 \subfloat[$11\times 11$ kernels after training.]{
  \centering
  \includegraphics[trim={0.5cm, 10.2cm, 2.0cm, 2.2cm}, clip, width=0.5\linewidth]{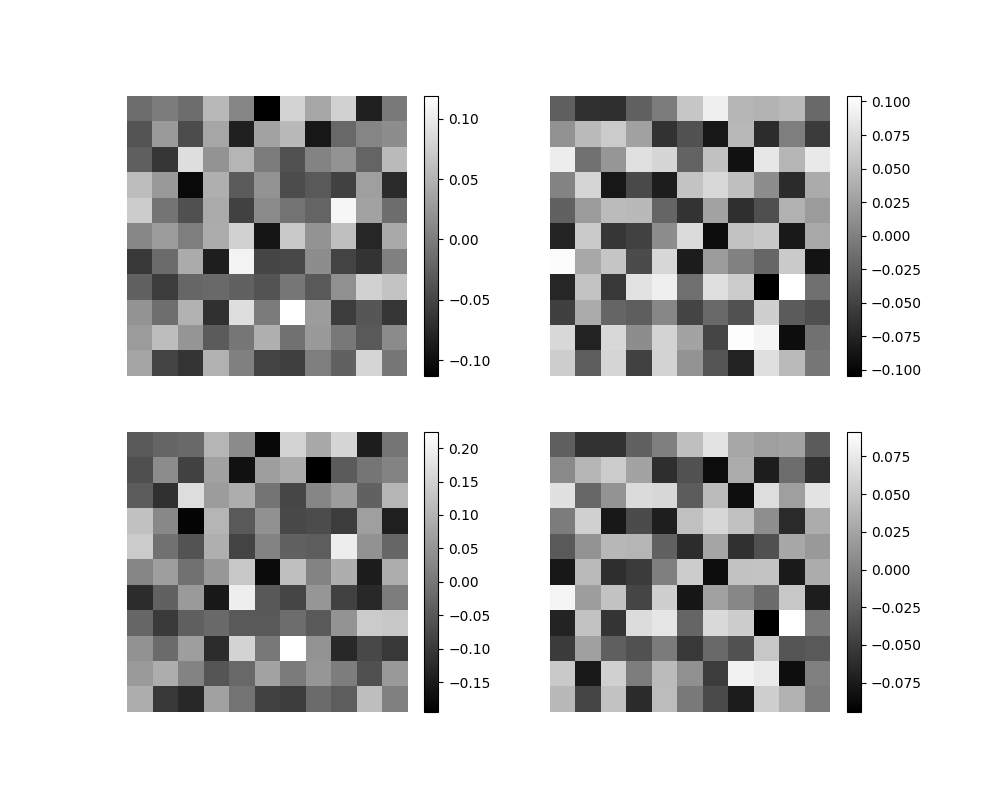}
  \includegraphics[trim={0.5cm, 1.67cm, 2.0cm, 10.7cm}, clip, width=0.5\linewidth]{review-celeba/cuda_run_007_kernels__after_training.png}
  }
 \caption{Comparison of $11\times 11$ kernels of $A$ (first two columns on the left) and of $C^\top$ (two columns on the right), shown before and after training on the denoising task using the DEQ model. 
 Note that, in the first row, the pairs of kernels in the first-third columns and the second-fourth are the same; this is because we initialize the kernels so that $C^\top=A^\top$ before training. After the training, they are different (second row).}
 \label{fig:comparison-kernels-11-11-denoising-DEQ}
\end{figure}

\begin{figure}[htbp!]
 \subfloat[$3\times 3$ kernels before training.]{
  \centering
  \includegraphics[trim={0, 11.8cm, 0, 0}, clip, width=1.0\linewidth]{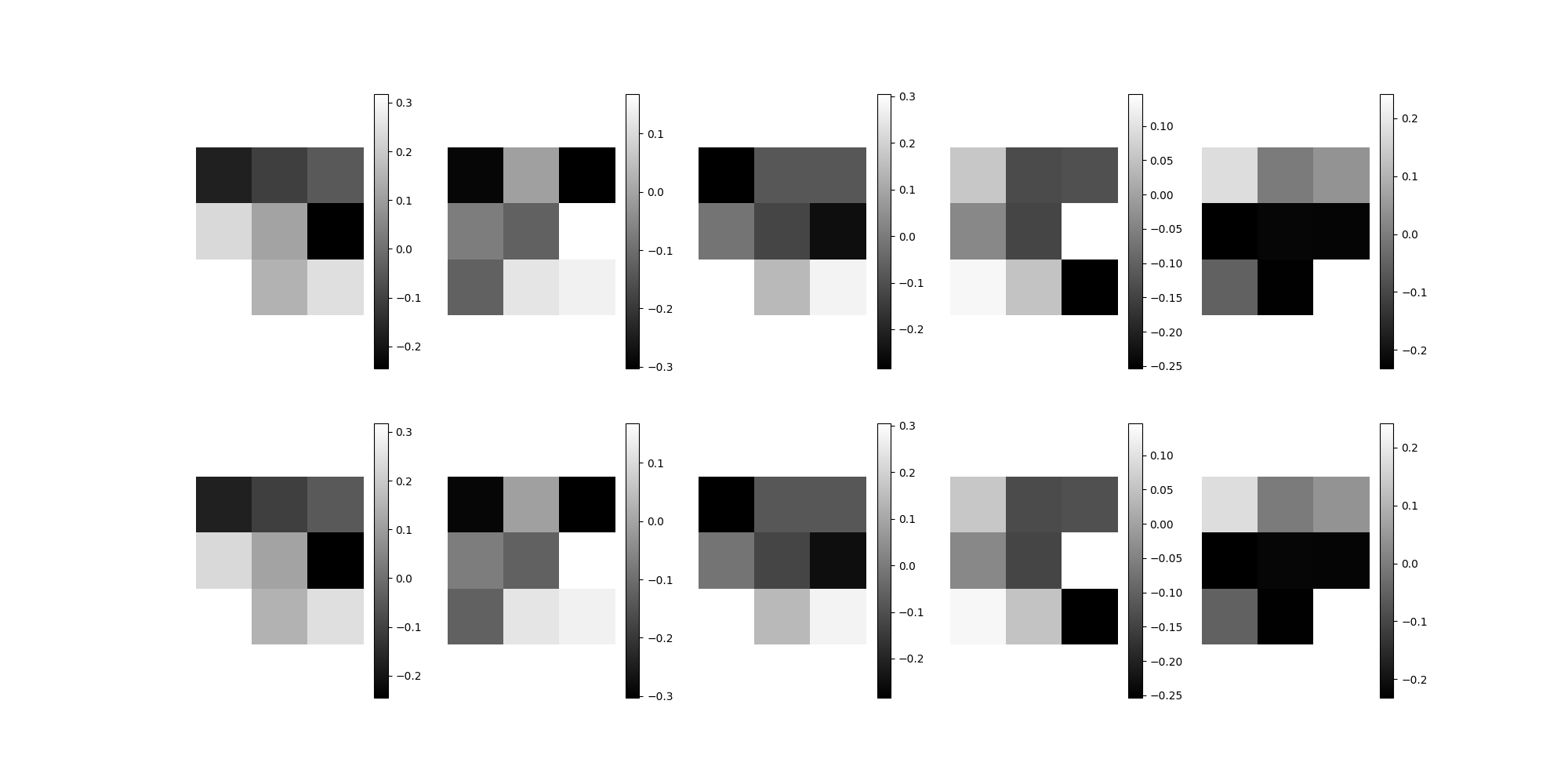}
  }\\
  \hfill
 \subfloat[$3\times 3$ kernels after training.]{
  \centering
  \includegraphics[trim={0, 12.0cm, 0, 2cm}, clip, width=1.0\linewidth]{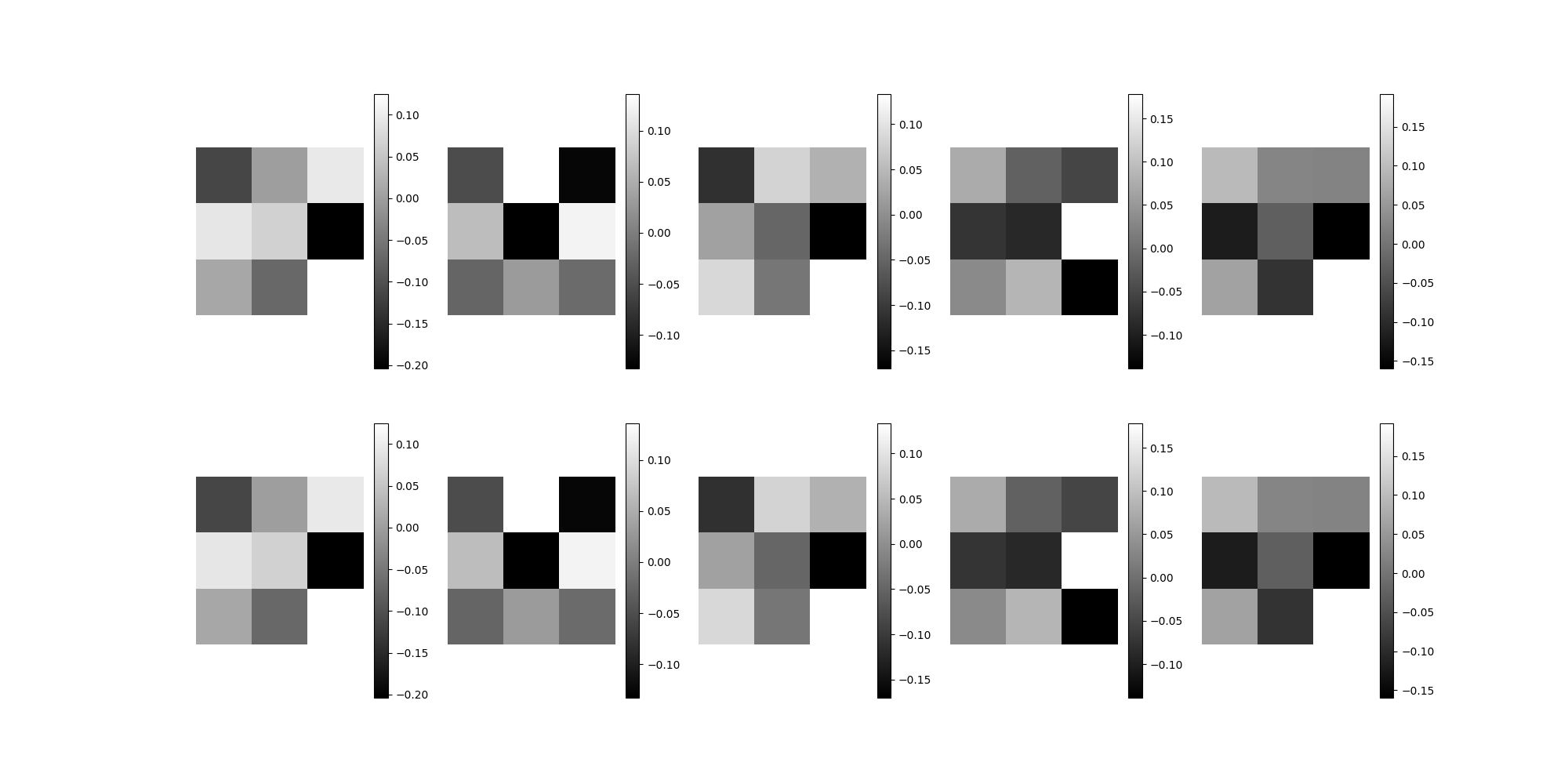}
  }
 \caption{Comparison of five (out of thirty) $3\times 3$ kernels of $A$, shown before and after training on the denoising task using the bilevel learning model. The weights of $C^\top$ are not shown here since $C^\top=A^\top$ for bilevel learning models.}
 \label{fig:comparison-kernels-3-3-denoising-BL}
\end{figure}

\begin{figure}[htbp!]
 \subfloat[$3\times 3$ kernels before training.]{
  \centering
  \includegraphics[trim={5cm, 11.8cm, 3cm, 2cm}, clip, width=0.5\linewidth]{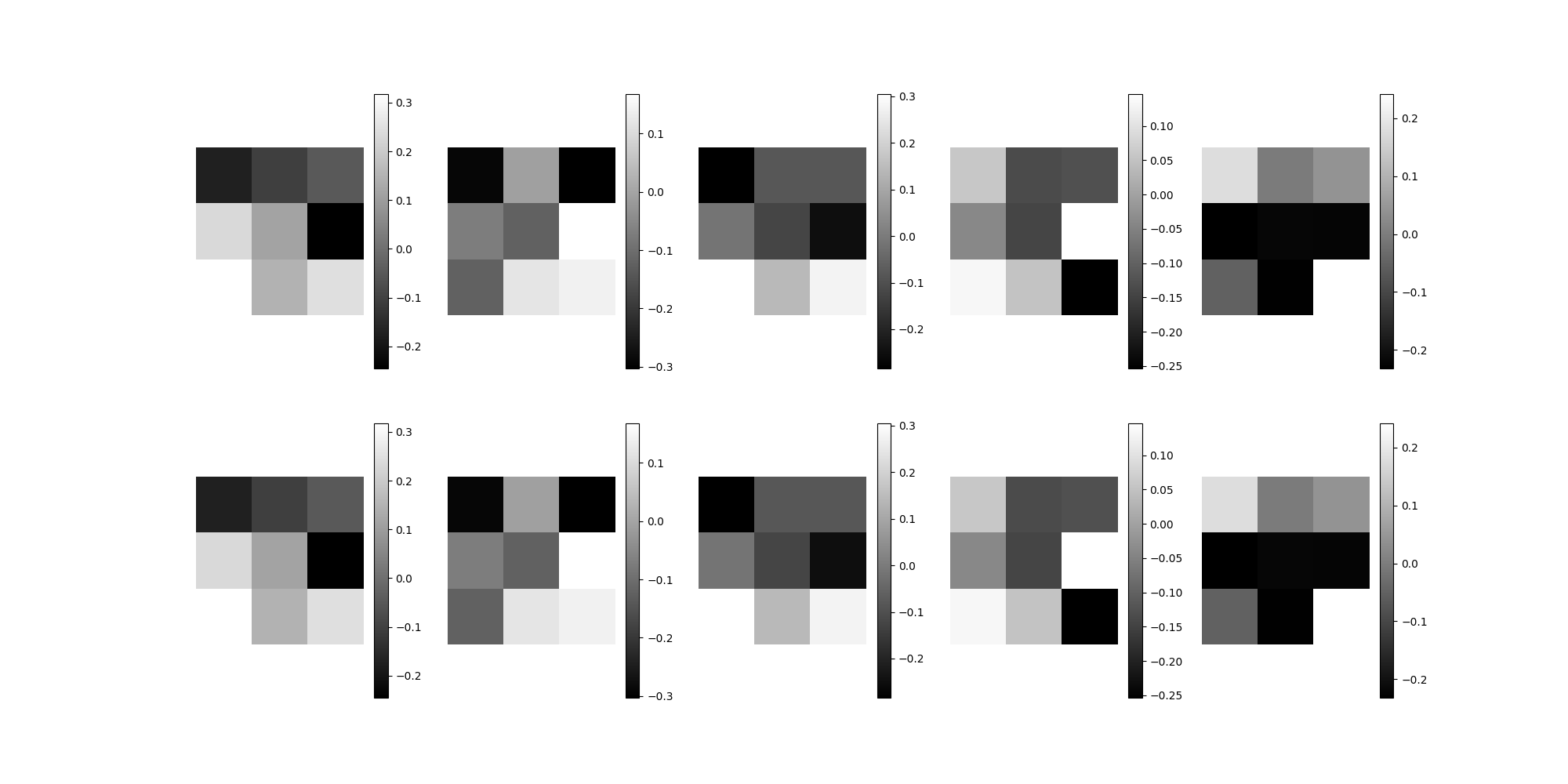}
  \includegraphics[trim={5cm, 11.8cm, 3cm, 2cm}, clip, width=0.5\linewidth]{review-celeba/cuda_run_008_kernels__before_training.png}
  }\\
  \hfill
 \subfloat[$3\times 3$ kernels after training.]{
  \centering
  \includegraphics[trim={5cm, 12.0cm, 3cm, 2cm}, clip, width=0.5\linewidth]{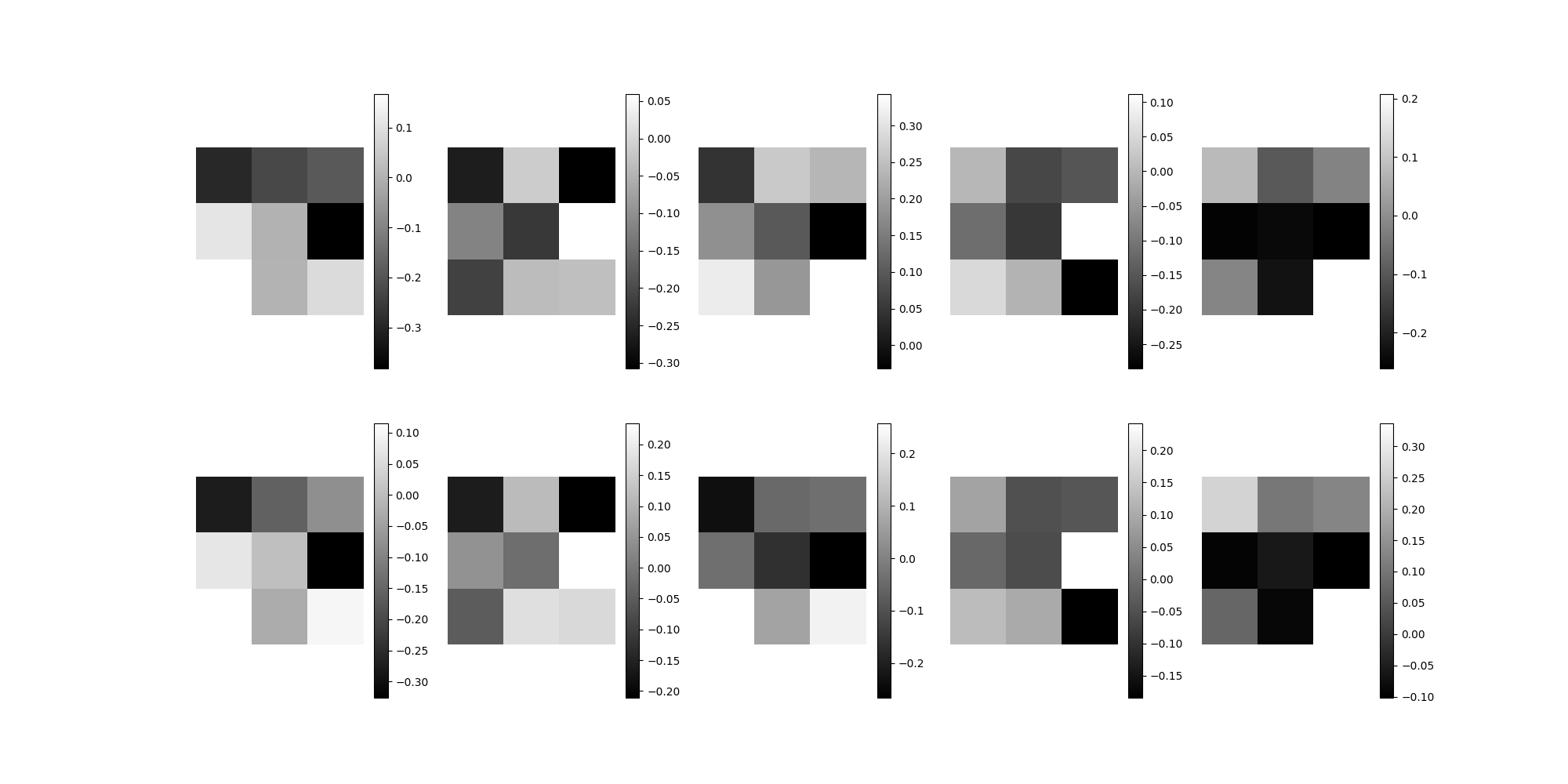}
  \includegraphics[trim={5cm, 12.0cm, 3cm, 2cm}, clip, width=0.5\linewidth]{review-celeba/cuda_run_008_kernels__after_training.png}
  }
 \caption{Comparison of $3\times 3$ kernels of $A$ (first five columns on the left) and of $C^\top$ (five columns on the right), shown before and after training on the denoising task using the DEQ model. 
 Note that, in the first row, kernels are pairwise equal; this is because we initialize the kernels so that $C^\top=A^\top$ before training. After the training, they are different (second row).}
 \label{fig:comparison-kernels-3-3-denoising-DEQ}
\end{figure}

The reconstruction of the images from noisy samples is shown in Fig.\ref{fig:denoising__celeba_spectral_norm}. The noise can be removed using either bilevel learning or DEQ models, with similar results in terms of MSE. Similarly, the quality of the reconstruction does not seem affected by the number of output channels and the size of the kernels, as long as the model is expressive enough. 
Interestingly, the kernels do not change too much even if the loss decreases. The major change is in the overall norm, whereas the relative values between the pairs of pixels remain similar (see Figs.\ref{fig:comparison-kernels-11-11-denoising-BL}-\ref{fig:comparison-kernels-11-11-denoising-DEQ} for a comparison of how much the kernels change in bilevel learning against DEQ methods on $11\times 11$ kernels, and Figs.\ref{fig:comparison-kernels-3-3-denoising-BL}-\ref{fig:comparison-kernels-3-3-denoising-DEQ} for $3 \times 3$ kernels). This is true for both the TV-like initialized kernels, where we could have expected to see the pattern change from the second-order to the first-order initialization (or vice-versa), and for the randomly initialized. Although this result seems counter-intuitive as it would be reasonable to expect to see some patterns appear in the kernels, a possible explanation is that the patterns are hidden in the kernels and are just not easily interpretable. Furthermore, it seems that any couple of randomly initialized kernels can perform denoising, provided that the weights are slightly changed.

\subsubsection{Deblurring}
\begin{figure}[hp]
 \centering \includegraphics[trim={0cm, 0, 8.7cm, 0}, clip, width=0.5\linewidth]{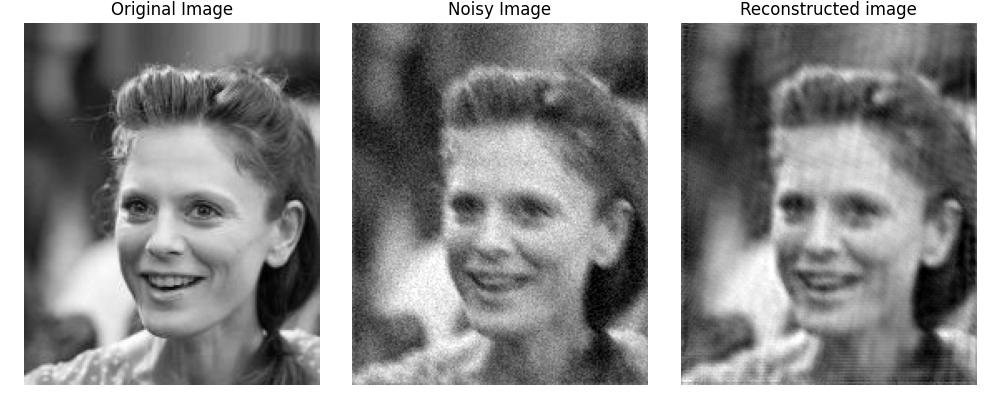}
 \\
 \subfloat[Reconstructed images, two output channels, kernels size $11\times 11$]{
 \begin{minipage}{0.25\textwidth}
  \centering
  \includegraphics[trim={17cm 0 0cm 0.505cm},clip, width=1.0\linewidth]{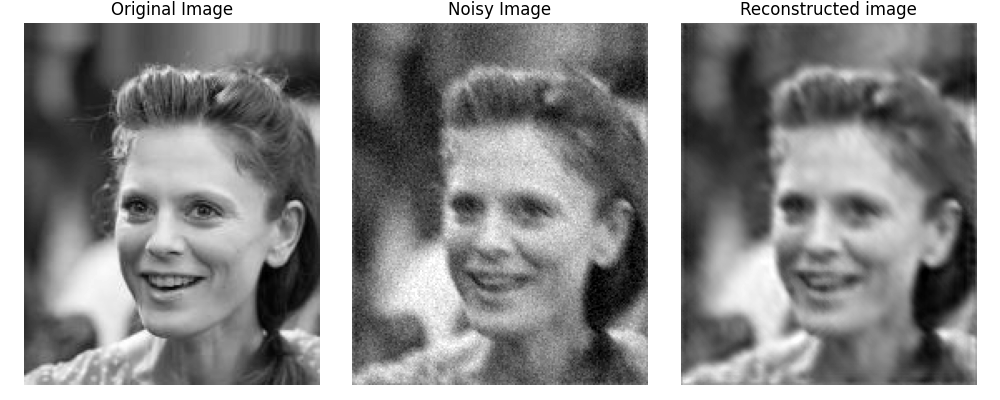}
 \end{minipage} \hfill
 \begin{minipage}{0.25\textwidth}
  \centering
  \includegraphics[trim={17cm 0 0cm 0.505cm},clip, width=1.0\linewidth]{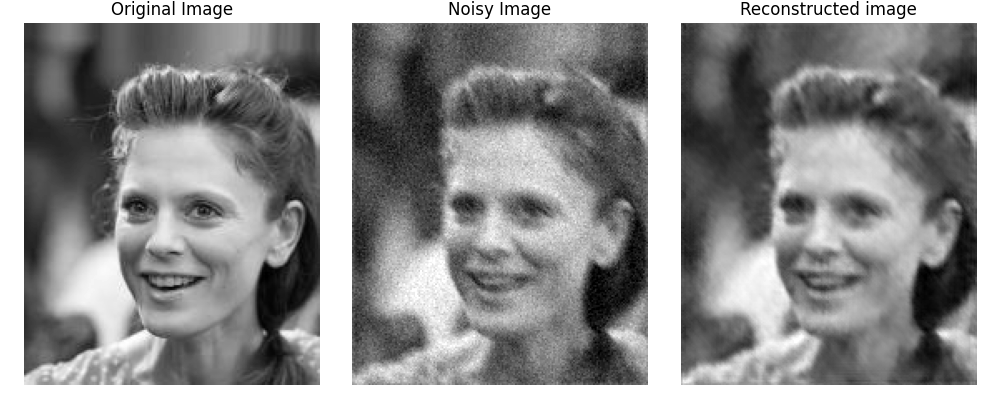}
 \end{minipage}
 \hfill
 \begin{minipage}{0.25\textwidth}
  \centering
  \includegraphics[trim={17cm 0 0cm 0.505cm}, clip, width=1.0\linewidth]{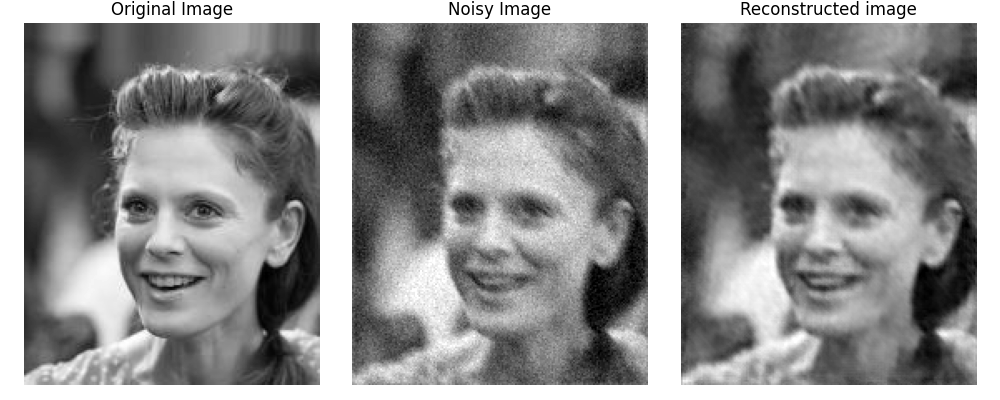}
 \end{minipage} \hfill 
 \begin{minipage}{0.25\textwidth}
  \centering
  \includegraphics[trim={17cm 0 0cm 0.505cm}, clip, width=1.0\linewidth]{review-celeba/cuda_run_023_reconstruction__before_training_test.png}
 \end{minipage} \hfill
 }\\
 \subfloat[Reconstructed images, thirty output channels, kernels size $3\times 3$]{
\begin{minipage}{0.25\textwidth}
  \centering
  \includegraphics[trim={17cm 0 0cm 0.505cm},clip, width=1.0\linewidth]{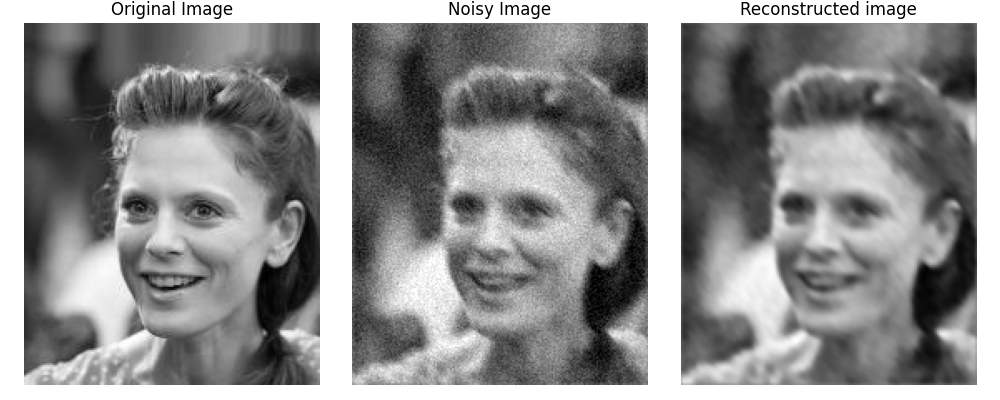}
 \end{minipage} \hfill
 \begin{minipage}{0.25\textwidth}
  \centering
  \includegraphics[trim={17cm 0 0cm 0.505cm},clip, width=1.0\linewidth]{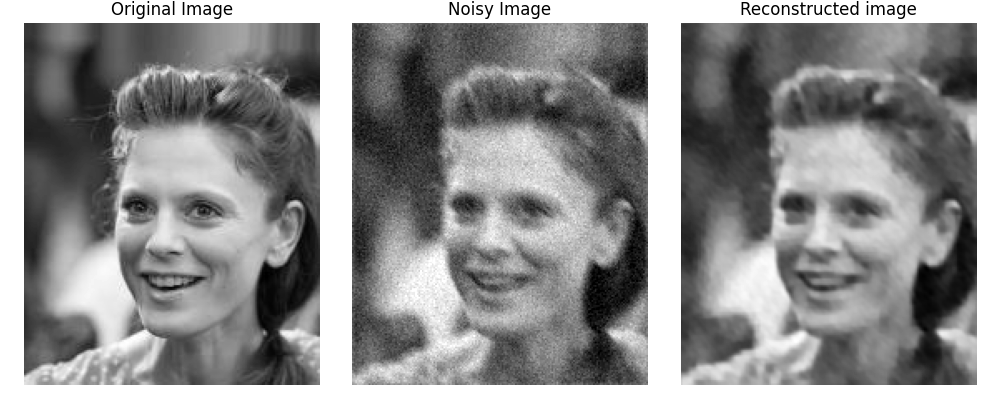}
 \end{minipage}
 \hfill
 \begin{minipage}{0.25\textwidth}
  \centering
  \includegraphics[trim={17cm 0 0cm 0.505cm}, clip, width=1.0\linewidth]{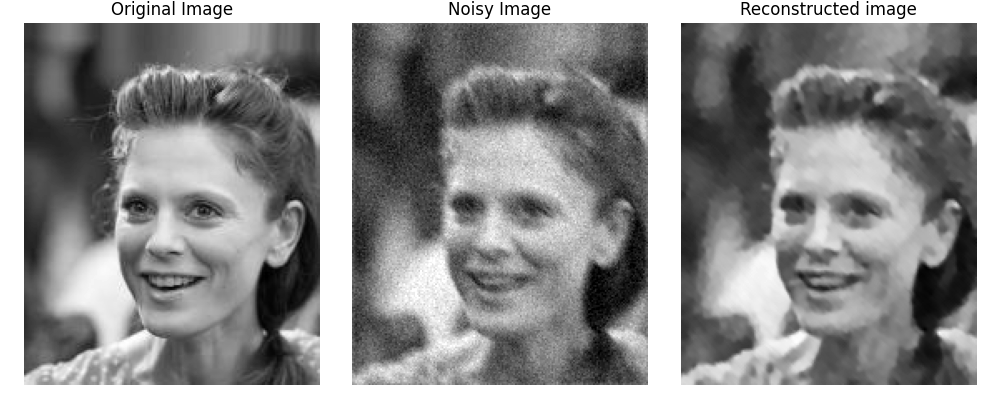}
 \end{minipage} \hfill 
 \begin{minipage}{0.25\textwidth}
  \centering
  \includegraphics[trim={17cm 0 0cm 0.505cm}, clip, width=1.0\linewidth]{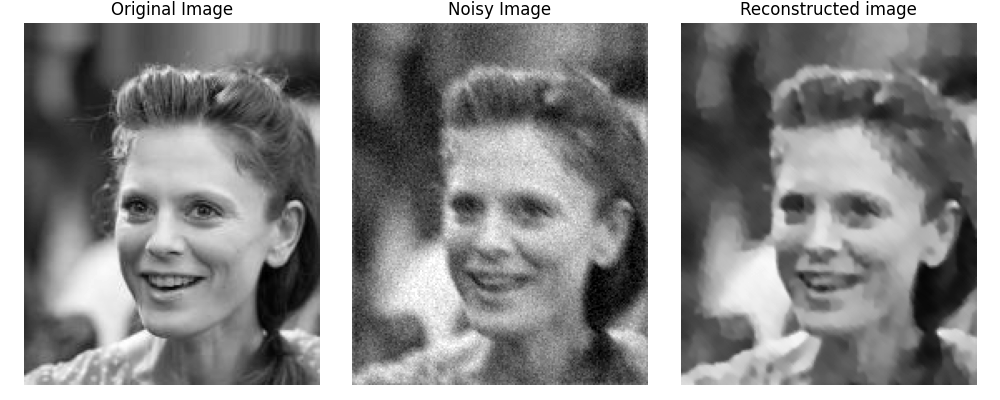}
 \end{minipage} \hfill
 }
 \caption{Deblurring CelebA, a sample from the test dataset. The first row contains the original image $u$ and the noisy blurred image $f^\delta$. From left to right in the second and third rows: reconstructed image with the optimal kernels found for the denoising task in the bilevel scenario (left), with parameters learned using bilevel learning (center-left), with parameters learned using DEQ model (center-right), and optimal kernels found for the denoising task in the DEQ scenario.}
 \label{fig:deblurring__celeba_spectral_norm}
\end{figure}

\begin{figure}[htbp!]
 \subfloat[$11\times 11$ kernels before training.]{
  \centering
  \includegraphics[trim={0, 10.3cm, 0, 1cm}, clip, width=1.0\linewidth]{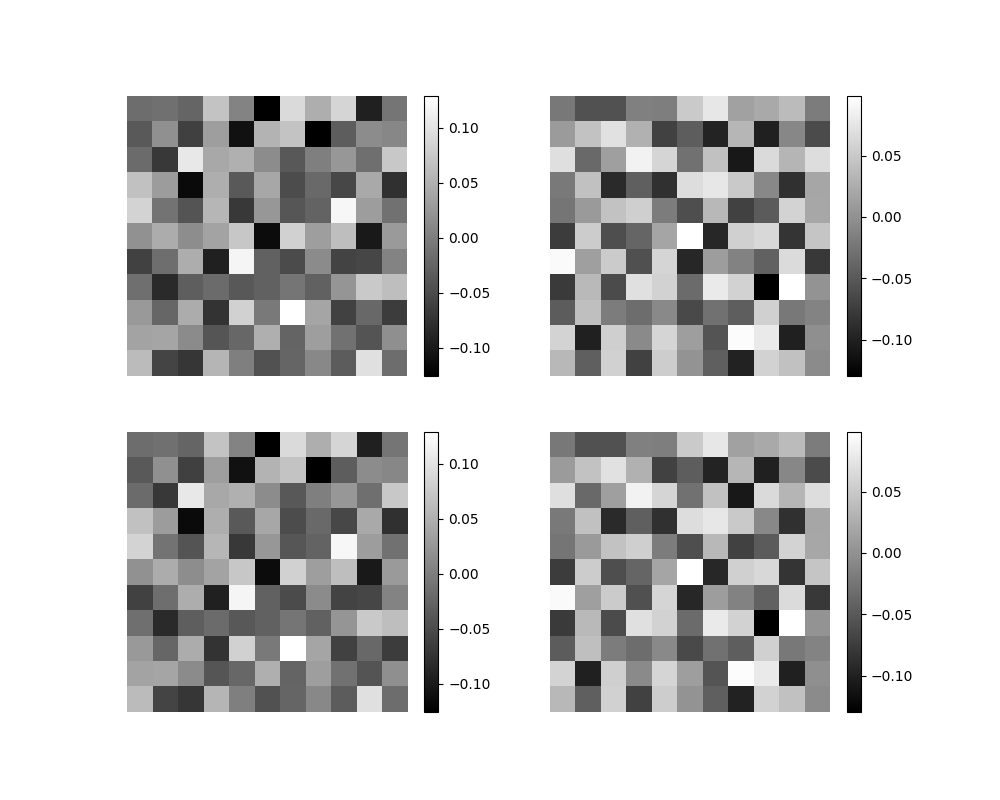}
  }\\
  \hfill
 \subfloat[$11\times 11$ kernels after training.]{
  \centering
  \includegraphics[trim={0, 10.2cm, 0, 2.1cm}, clip, width=1.0\linewidth]{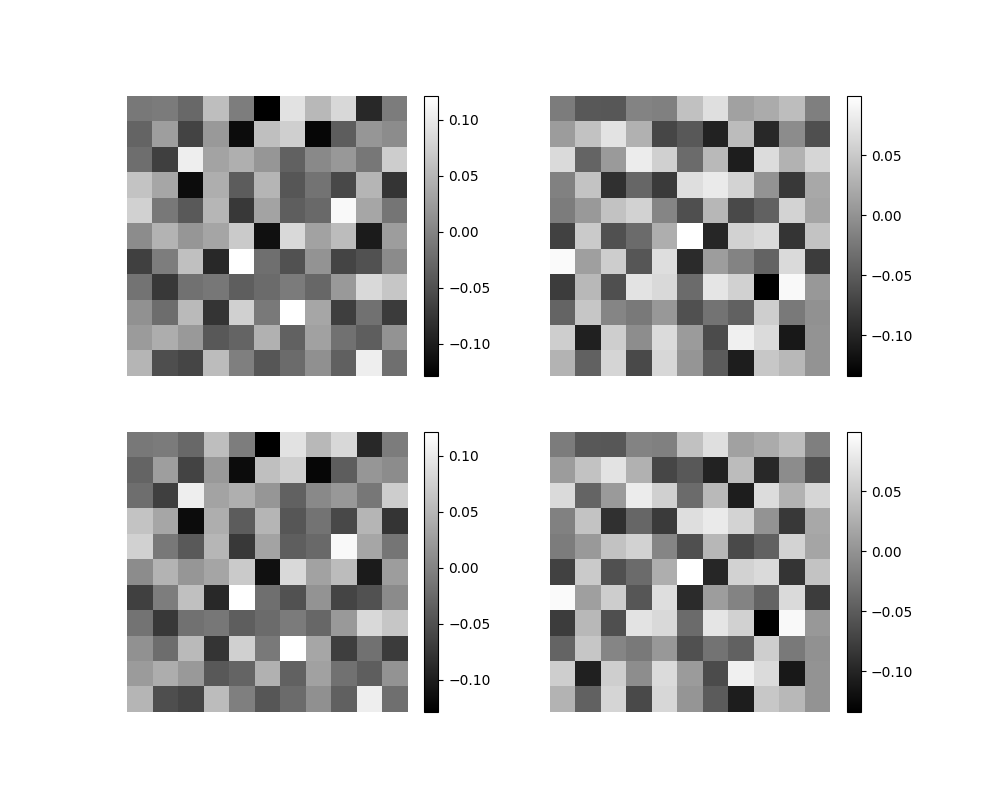}
  }
 \caption{Comparison of $11\times 11$ kernels of $A$, shown before and after training on the deblurring task using the bilevel learning model. The weights of $C^\top$ are not shown here since $C^\top=A^\top$ for bilevel learning models.}
 \label{fig:comparison-kernels-11-11-deblurring-BL}
\end{figure}

\begin{figure}[htbp!]
 \subfloat[$11\times 11$ kernels before training.]{
  \centering
  \includegraphics[trim={0.5cm, 10.2cm, 2.0cm, 2.2cm}, clip, width=0.5\linewidth]{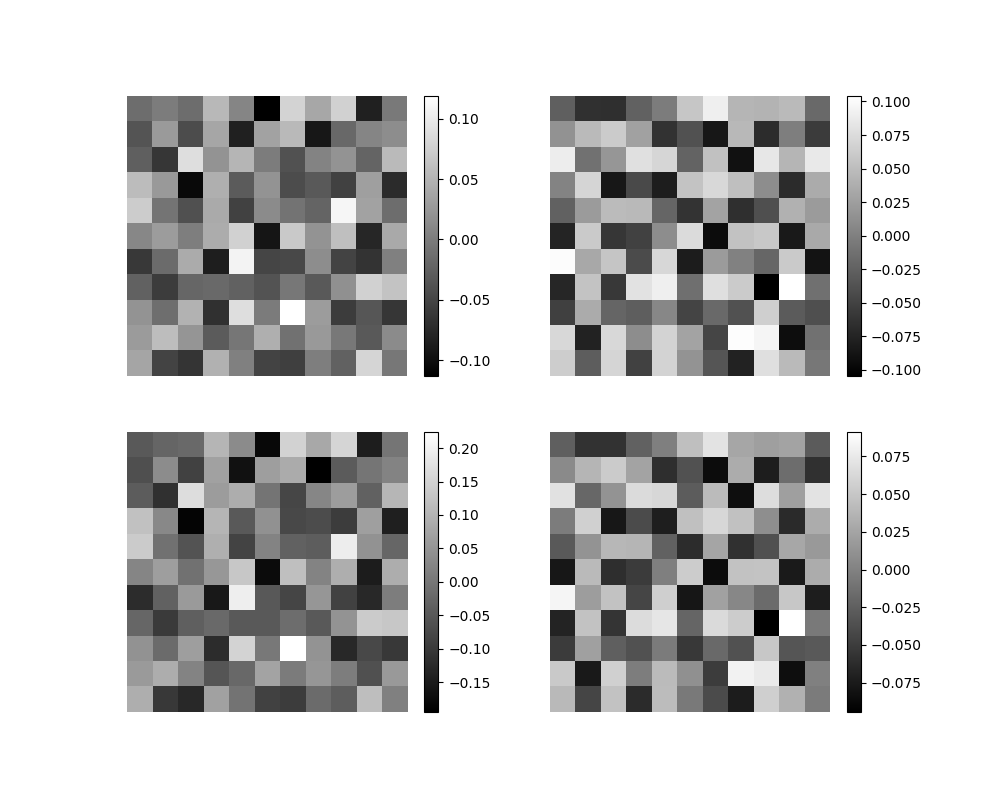}
  \includegraphics[trim={0.5cm, 1.67cm, 2.0cm, 10.7cm}, clip, width=0.5\linewidth]{review-celeba/cuda_run_023_kernels__before_training.png}
  }\\
  \hfill
 \subfloat[$11\times 11$ kernels after training.]{
  \centering
  \includegraphics[trim={0.5cm, 10.2cm, 2.0cm, 2.2cm}, clip, width=0.5\linewidth]{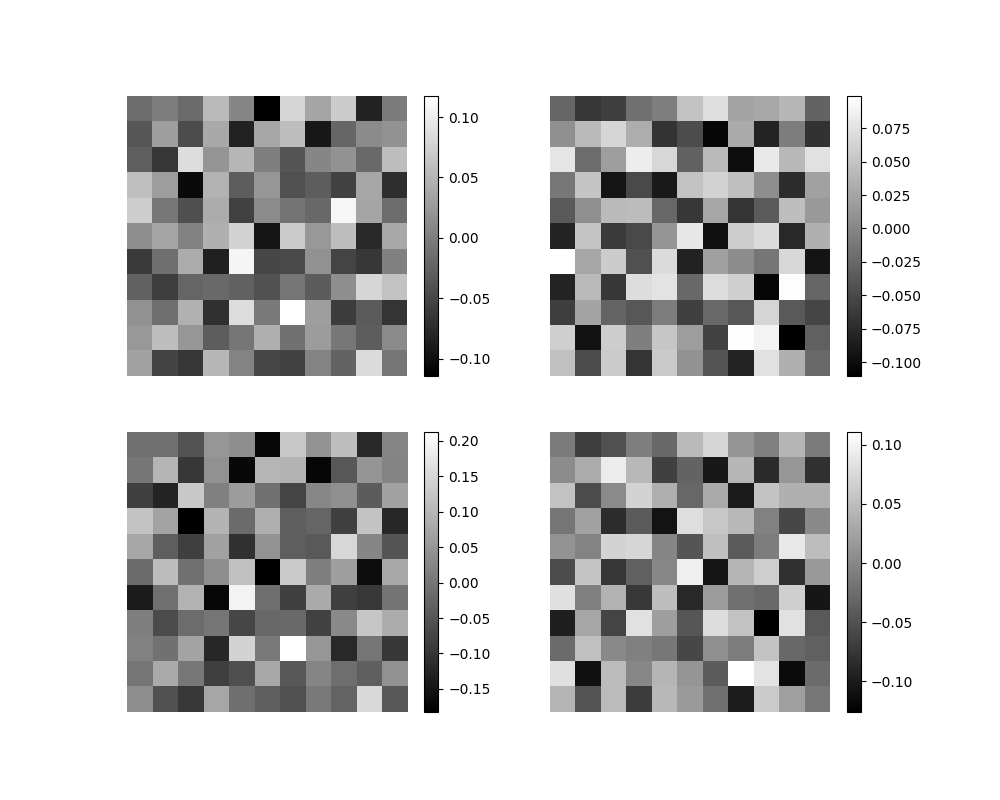}
  \includegraphics[trim={0.5cm, 1.67cm, 2.0cm, 10.7cm}, clip, width=0.5\linewidth]{review-celeba/cuda_run_023_kernels__after_training.png}
  }
 \caption{Comparison of $11\times 11$ kernels of $A$ (first two columns on the left) and of $C^\top$ (two columns on the right), shown before and after training on the deblurring task using the DEQ model. 
 Note that, in the first row, the pairs of kernels in the first-third columns and second-fourth are the same; this is because we initialize the kernels so that $C^\top=A^\top$ before training. After the training, they are different (second row).}
 \label{fig:comparison-kernels-11-11-deblurring-DEQ}
\end{figure}

\begin{figure}[htbp!]
 \subfloat[$3\times 3$ kernels before training.]{
  \centering
  \includegraphics[trim={0, 11.8cm, 0, 0}, clip, width=1.0\linewidth]{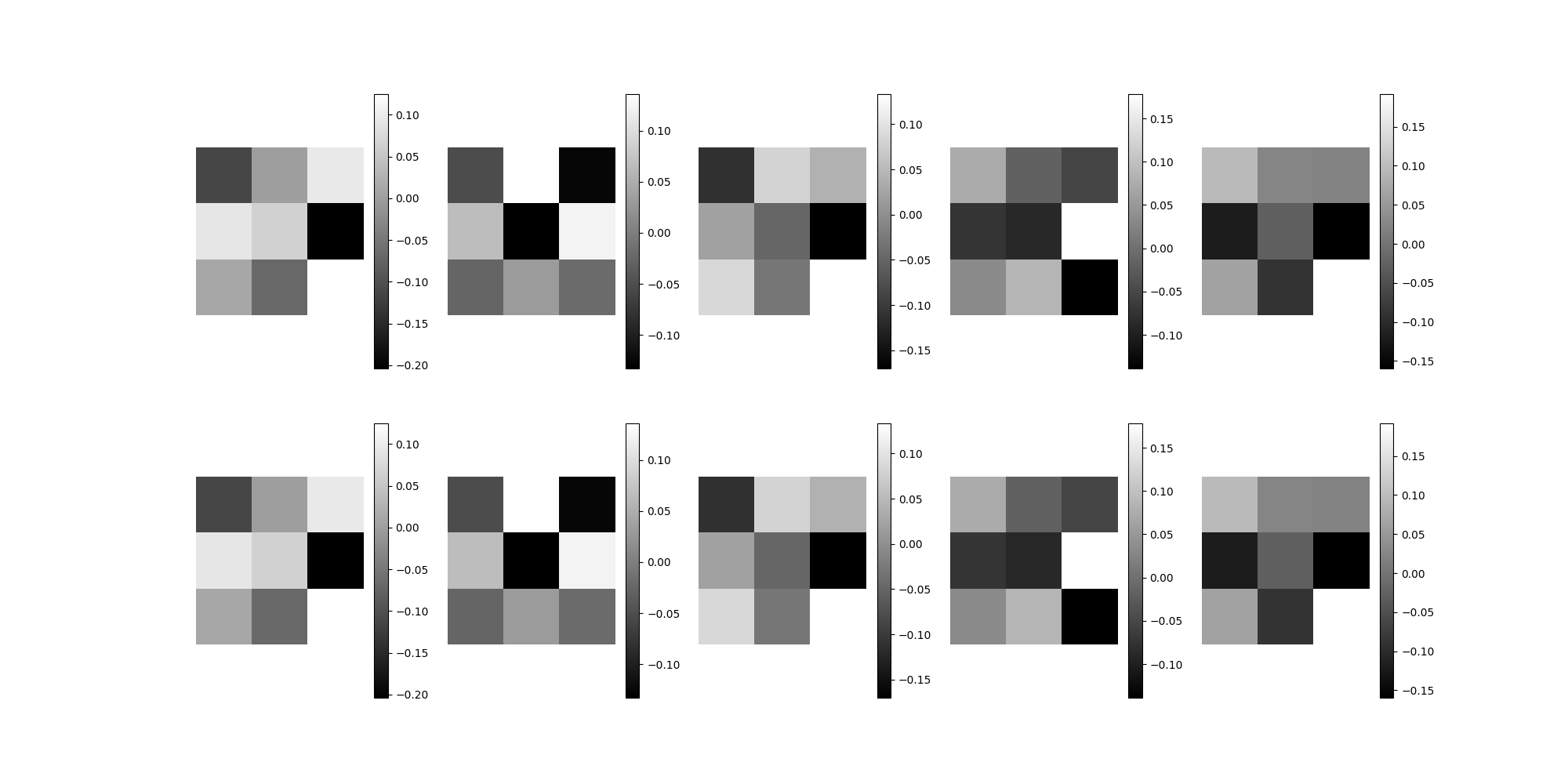}
  }\\
  \hfill
 \subfloat[$3\times 3$ kernels after training.]{
  \centering
  \includegraphics[trim={0, 12.0cm, 0, 2cm}, clip, width=1.0\linewidth]{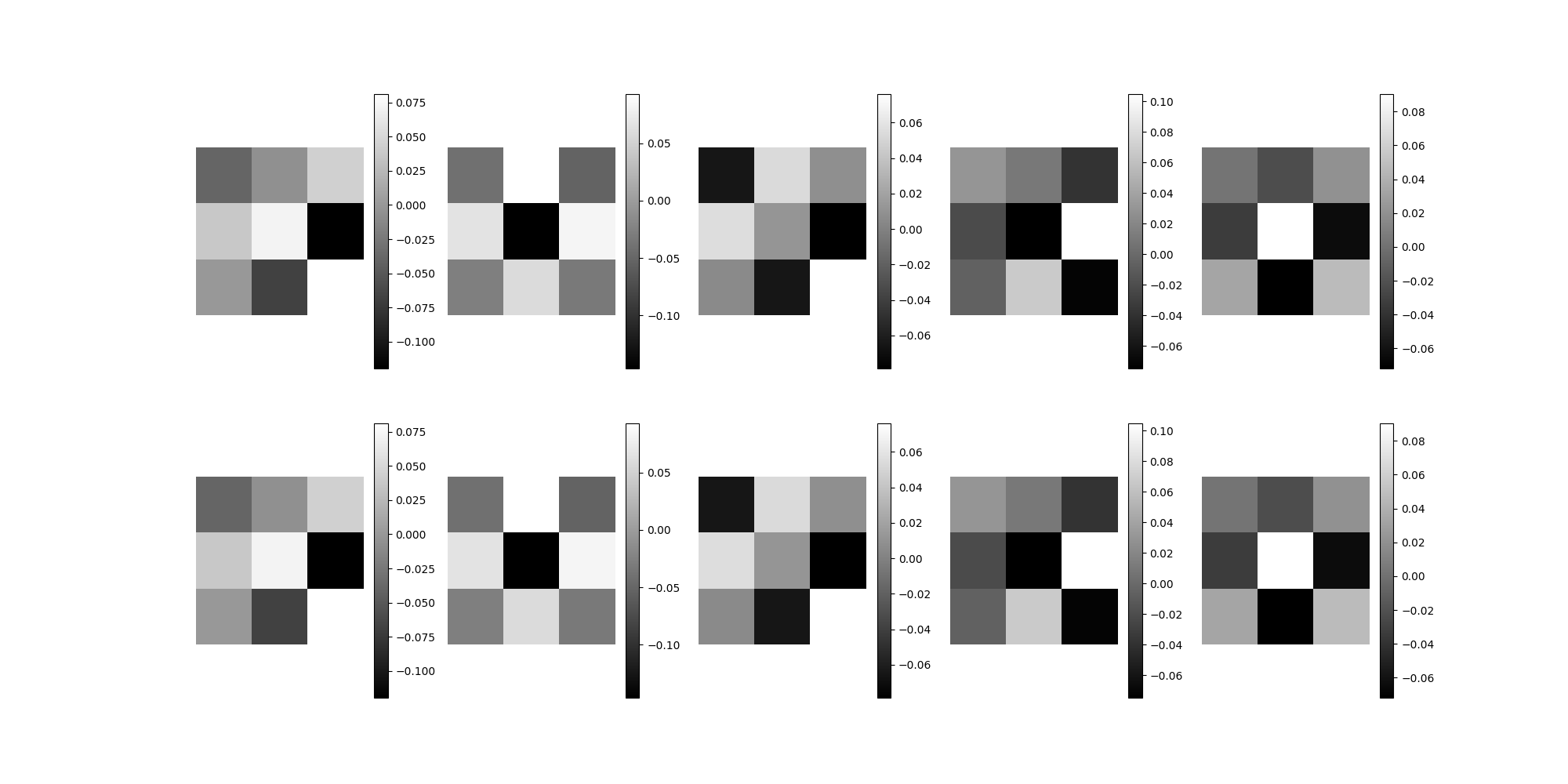}
  }
 \caption{Comparison of five (out of thirty) $3\times 3$ kernels of $A$, shown before and after training on the deblurring task using the bilevel learning model. The weights of $C^\top$ are not shown here since $C^\top=A^\top$ for bilevel learning models.}
 \label{fig:comparison-kernels-3-3-deblurring-BL}
\end{figure}

\begin{figure}[htbp!]
 \subfloat[$3\times 3$ kernels before training.]{
  \centering
  \includegraphics[trim={5cm, 11.8cm, 3cm, 2cm}, clip, width=0.5\linewidth]{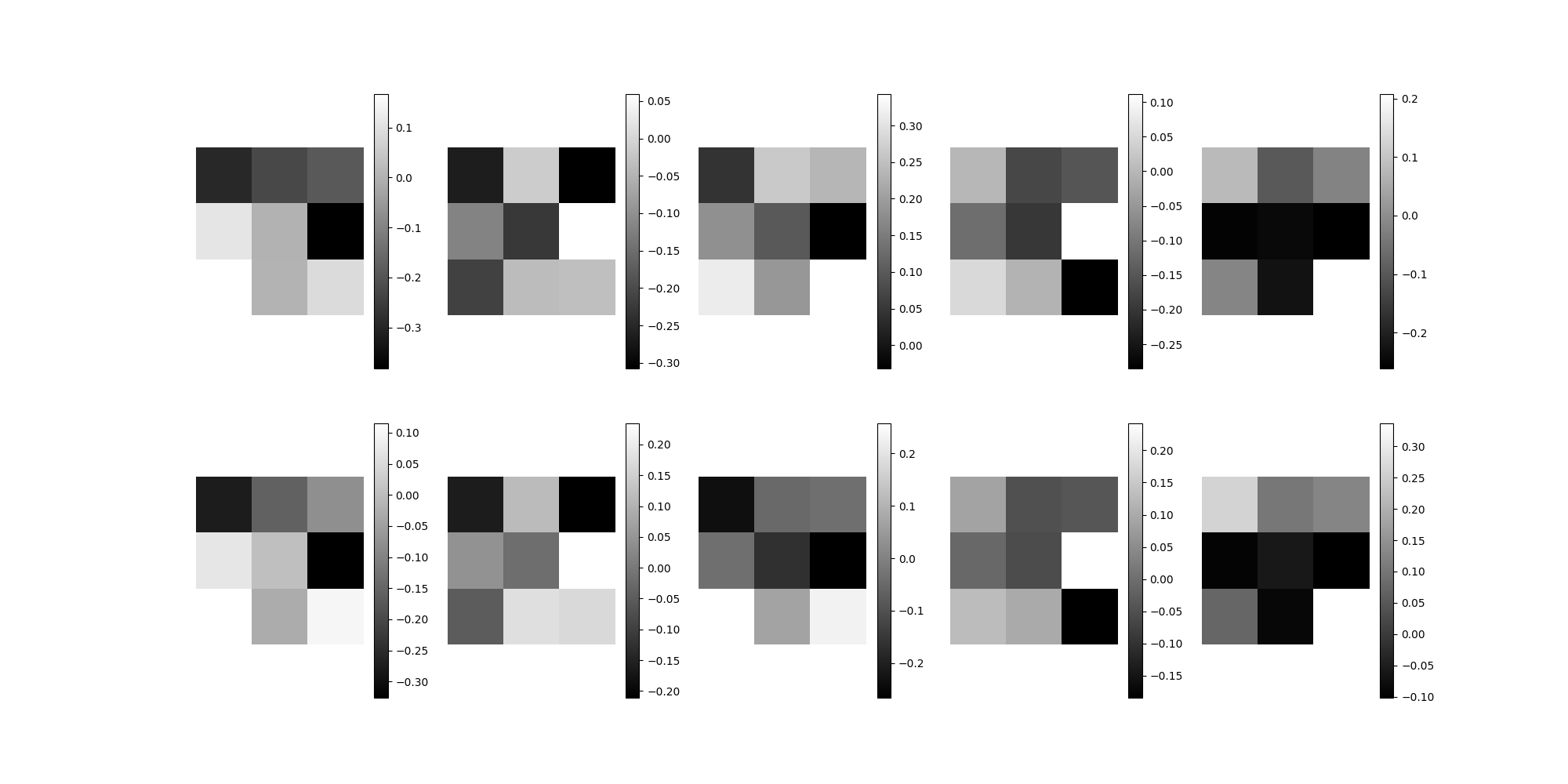}
  \includegraphics[trim={5cm, 11.8cm, 3cm, 2cm}, clip, width=0.5\linewidth]{review-celeba/cuda_run_024_kernels__before_training.png}
  }\\
  \hfill
 \subfloat[$3\times 3$ kernels after training.]{
  \centering
  \includegraphics[trim={5cm, 12.0cm, 3cm, 2cm}, clip, width=0.5\linewidth]{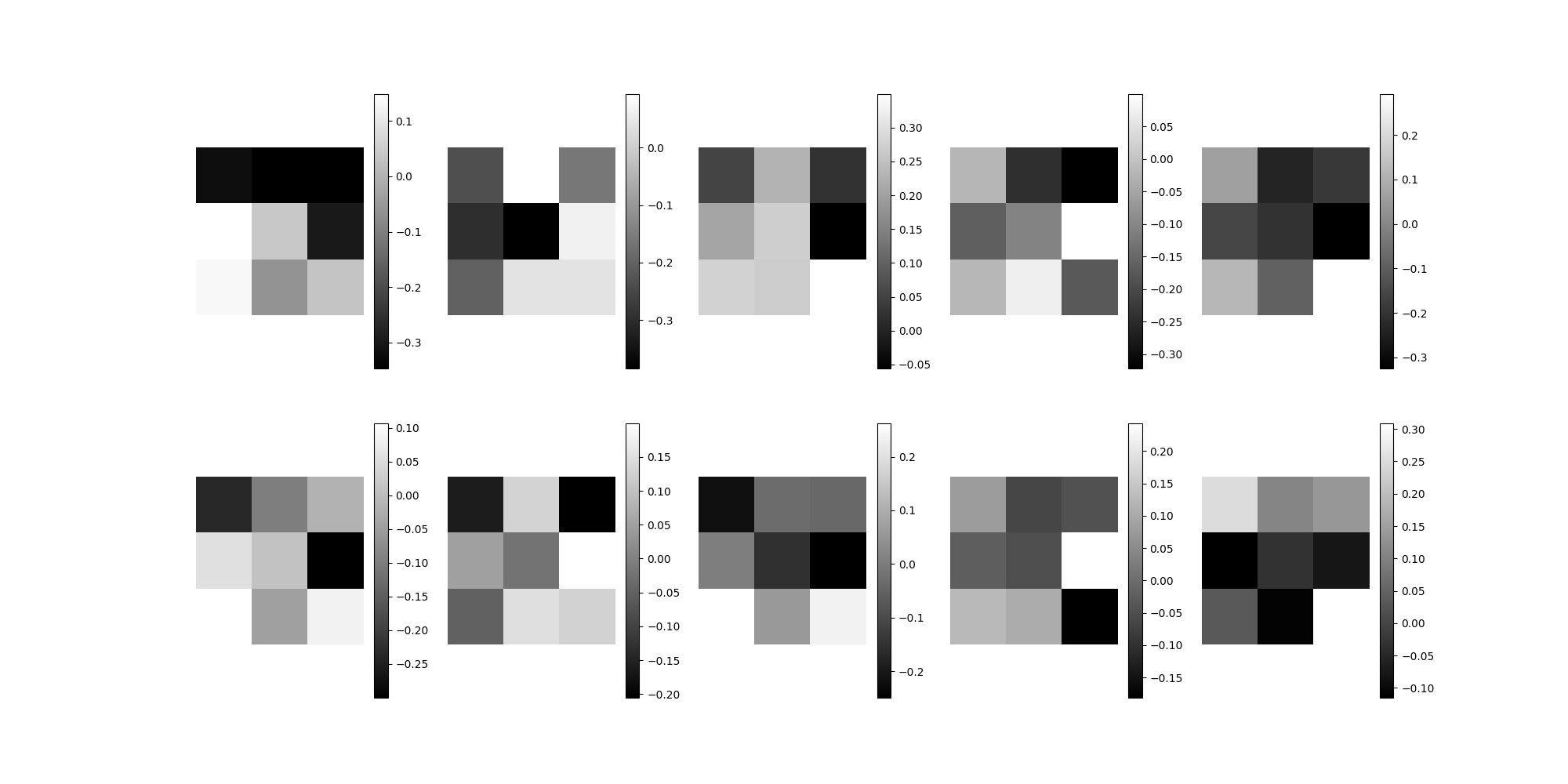}
  \includegraphics[trim={5cm, 12.0cm, 3cm, 2cm}, clip, width=0.5\linewidth]{review-celeba/cuda_run_024_kernels__after_training.png}
  }
 \caption{Comparison of $3\times 3$ kernels of $A$ (first five columns on the left) and of $C^\top$ (five columns on the right), shown before and after training on the deblurring task using the DEQ model. 
 Note that, in the first row, kernels are pairwise equal; this is because we initialize the kernels so that $C^\top=A^\top$ before training.}
 \label{fig:comparison-kernels-3-3-deblurring-DEQ}
\end{figure}

The second inverse problem we consider for the CelebA dataset is deblurring. We perform a warm-start, loading the optimal kernels' weights learned in the denoising task.
The reconstructed images are shown in Fig.\ref{fig:deblurring__celeba_spectral_norm}. Both bilevel learning and DEQ models achieve a similar quality in the reconstructed images. Considerations similar to the ones drawn in the denoising task can be done on the learned kernels for the deblurring task
(Figs.\ref{fig:comparison-kernels-11-11-deblurring-BL}-\ref{fig:comparison-kernels-11-11-deblurring-DEQ} show a comparison of how much the kernels change in bilevel learning against DEQ methods on $11\times 11$ kernels, and Figs.\ref{fig:comparison-kernels-3-3-deblurring-BL}-\ref{fig:comparison-kernels-3-3-deblurring-DEQ} for $3 \times 3$ kernels).

\section{Conclusions \& Outlook}\label{sec:conclusions-and-outlook}
We have framed bilevel learning methods as deep equilibrium methods and have compared several models of each class of similar complexity. From those numerical results, we have observed that the two methods behave similarly in terms of average loss when they are used as regularizers. We can even argue that regularizers trained by bilevel learning were performing slightly better than their deep equilibrium counterparts, which can have different reasons. The average loss of the trained models with bilevel learning is smaller, with a lower interquantile range than deep equilibrium models. The results also suggest bilevel learning is less sensitive to the choice of hyperparameters than deep equilibrium method, making it more robust in terms of hyperparameters selection, and easier to train.
This is in contrast to the a-priori assumption that deep equilibrium models might perform better because they are more general than their bilevel learning counterparts. The number of parameters in the deep equilibrium models we considered is around twice than the ones of the bilevel learning models. 
The experiments suggest that the extra constraints imposed by requiring the regularizing network to be the gradient of a regularizer does not hinder the performance, and may even enhance it.
We emphasize that these observations are limited to the gradient descent deep equilibrium architecture. Further comparisons with other deep equilibrium architectures are subject to future work. %
It would also be interesting to study deep equilibrium models mathematically in greater detail, and to identify conditions that can lead to theoretical guarantees for them. 

\section*{Acknowledgements}
The authors would like to thank the Isaac Newton Institute for Mathematical Sciences, Cambridge, for support and hospitality during the programme \emph{Mathematics of Deep Learning} where work on this paper was undertaken. This work was supported by EPSRC grant no EP/R014604/1. This research utilized Queen Mary's Apocrita and Andrena HPC facilities, supported by QMUL Research-IT \url{http://doi.org/10.5281/zenodo.438045}. DR acknowledges support from EPSRC grant EP/R513106/1. MB acknowledges support from the Alan Turing Institute. MJE acknowledges support from the EPSRC (EP/S026045/1, EP/T026693/1, EP/V026259/1) and the
Leverhulme Trust (ECF-2019-478).

\clearpage
\printbibliography


\end{document}